\documentclass{elsarticle}

\usepackage{lineno,hyperref}
\modulolinenumbers[5]

\journal{}

\usepackage{amsmath}
\usepackage{lineno,hyperref}
\usepackage{subfigure}
\usepackage{latexsym}
\usepackage{amsmath}
\usepackage{amssymb}
\usepackage{amsfonts}
\usepackage{verbatim}
\usepackage[latin1]{inputenc}
\usepackage{mathrsfs}
\usepackage{amsfonts}
\usepackage{graphicx}
\usepackage{colortbl,dcolumn}
\usepackage{amsmath}
\usepackage{psfrag}
\usepackage{booktabs}
\usepackage{lipsum}
\usepackage{amsfonts}
\usepackage{graphicx}
\usepackage{epstopdf}
\usepackage{algorithmic}
\ifpdf
  \DeclareGraphicsExtensions{.eps,.pdf,.png,.jpg}
\else
  \DeclareGraphicsExtensions{.eps}
\fi
\usepackage[misc]{ifsym}
\usepackage{enumerate}

\newtheorem{as}{Assumption}
\newtheorem{tm}{Theorem}
\newtheorem{rk}{Remark}

\newtheorem{prop}{Proposition}
\newtheorem{lm}{Lemma}
\newtheorem{cor}{Corollary}

\newcommand{\E}{\mathbb E}
\newcommand{\PP}{\mathbb P}
\newcommand{\N}{\mathbb N}
\newcommand{\R}{\mathbb R}

\newcommand{\LL}{\mathcal L}

\newcommand{\HH}{\mathbb H}

\newcommand{\FFF}{\mathscr F}
\newcommand{\<}{\langle}
\renewcommand{\>}{\rangle}
\usepackage{graphicx}
\allowdisplaybreaks[4]

\begin{document}

\begin{frontmatter}

\title{Weak convergence and invariant measure
of a full discretization for parabolic SPDEs with non-globally Lipschitz coefficients}
\tnotetext[mytitlenote]{This work is partially  supported by National Natural Science Foundation of China (NO. 91530118, NO. 91130003, NO. 11021101, NO. 91630312 and NO. 11290142).}

\author[cas]{Jianbo Cui\corref{cor}}
\ead{jianbocui@lsec.cc.ac.cn}

\author[cas]{Jialin Hong}
\ead{hjl@lsec.cc.ac.cn}

\author[cas]{Liying Sun}
\ead{liyingsun@lsec.cc.ac.cn}

\cortext[cor]{Corresponding author.}

\address[cas]{1. LSEC, ICMSEC, 
			Academy of Mathematics and Systems Science, Chinese Academy of Sciences, Beijing,  100190, China\\
			2. School of Mathematical Science, University of Chinese Academy of Sciences, Beijing, 100049, China}

\begin{abstract}

Approximating the invariant measure and the expectation of the functionals for parabolic stochastic partial differential equations (SPDEs)  with non-globally Lipschitz coefficients is an active research area and is far from being well understood. 
In this article, we study such problem in terms of a  full discretization based on the spectral Galerkin method and the temporal implicit Euler scheme.
By deriving the a priori estimates and regularity estimates of the numerical solution via a variational approach and Malliavin calculus, we establish the sharp weak convergence rate of the full discretization.
When the SPDE admits a unique $V$-uniformly ergodic invariant measure, we prove that the invariant measure can be approximated by the full discretization.
The key ingredients lie on the time-independent weak convergence analysis and time-independent regularity estimates of the corresponding Kolmogorov equation.
Finally, numerical experiments confirm the theoretical findings.
\end{abstract}


\end{frontmatter}

\linenumbers


\section{Introduction}

Numerical approximations of stochastic partial differential equations (SPDEs) with local Lipschitz continuous coefficients, as an active area of research, has been widely concerned in the recent years. Weak approximation of such SPDE is still far from well understood, although some progress has been made.
We are only aware that there are several results on weak convergent semi-discretizations, such as some temporal splitting methods in \cite{BG18b} and the spatial finite element methods in \cite{CH18} for parabolic SPDEs with non-globally Lipschitz coefficients, and a temporal splitting method in \cite{CH17} for the stochastic Schr\"odinger equation with a cubic nonlinearity. 
To the best of our knowledge, there has been no essentially sharp weak convergence rate result of full discretization for parabolic SPDEs with non-globally Lipschitz coefficients.

One motivation of this present work is considering this direction and studying the numerical weak 
approximation for the
following parabolic SPDE 
\begin{equation}
\begin{split}
\label{spde}
dX(t)&=(AX(t)+F(X(t)))dt+dW(t),\; t>0, \\
X(0)&=X_0,\; 
\end{split}
\end{equation}
where $A$ is the Laplacian operator on  $\mathcal O:=[0,L]^d, d\le 3, L>0$ under homogenous Dirichlet boundary condition, $F$ is the Nemytskii operator of a real-valued one-sided Lipschitz function $f$ and $\{W(t)\}_{t\ge 0}$ is a generalized $Q$-Wiener process on a filtered probability space $(\Omega,\mathcal F,\PP,\{\mathcal F_t\}_{t\ge 0})$ (see Section 2 for details).
After discretizing Eq. \eqref{spde} by the spatial spectral Galerkin method, we apply the temporal backward Euler method to propose the full discretization \eqref{full-spde}.
Let $N$ be the dimension of the spectral Galerkin projection space and $\delta t$ be the time stepsize.
Let $\lambda_1$ be the smallest eigenvalue of $-A$,
and 
$\lambda_F$ be the one-sided Lipschitz constant of $F$.
Denote $\HH:=L^2(\mathcal O)$, and $\mathcal C_b^2(\HH,\R)$ 
the space of twice continuous differentiable functionals from $\HH$ to $\R$ with bounded first and second derivatives. 
By denoting $\{X^N_k\}_{k\in \N^+}$ the numerical solution of \eqref{full-spde}, the essentially sharp weak convergence rate of \eqref{full-spde} is shown in the following theorem.

\begin{tm}\label{weak0}
Let Assumptions \ref{as-lap}-\ref{as-dri} hold with $\beta \in (0,1]$, $\gamma\in (0,\beta)$,
$X_0\in \HH^{\frac d2+\epsilon}$ with a sufficient small positive constant $\epsilon$, $T>0$ and $\delta t_0\in 
(0,1\land \frac 1{(2\lambda_F-2\lambda_1)\lor 0})$. Then for any  $\phi\in \mathcal C_b^2(\HH)$, 
there exists  $C(T,X_0,Q,\phi)>0$ such that for any $\delta t\in (0, \delta t_0]$, $K\delta t=T$, $K\in\N^+$ and $N\in \N^+$,
\begin{align*}
\Big|\E\Big[\phi(X(T))-\phi(X^N_{K})\Big]\Big|
\le C(T,X_0,Q,\phi)\Big(\delta t^{\gamma}+\lambda_N^{-\gamma}\Big).
\end{align*}
\end{tm}

For SPDEs with smooth and regular coefficients, 
there have already been different approaches to studying the weak convergence rate of numerical methods 
(see e.g. \cite{AL16,BD18,CJK14,Deb11,KLL13,WG13}). However,
for the full discretization of parobolic SPDEs with superlinearly growing nonlinearities, it is still unclear  how to analyze its sharp weak convergence rate.
The key points to gain the error estimate in Theorem \ref{weak0} are applications of  the regularity estimates of the regularized Kolmogorov equation, and the a priori estimates of the numerical solution  and its Malliavin derivative.
We would like to mention that proving this result confronts at two main difficulties, one being the full implicity of the proposed method and another being to get the a priori estimates independent of both $N$ and $\delta t$ for \eqref{full-spde}.
These a priori estimates are not trivial due to the loss of the maximum principle for the analytic semigroup.
To overcome these difficulties, 
we make use of 
some Sobolev--Gagliard--Nirenberg inequalities, the It\^o formula for Skorohod integrals 
and the equivalence between a random PDE and Eq. \eqref{spde}.
Meanwhile, the approach to the weak convergence analysis is also available for other numerical methods at any finite time.

Based on the weak error analysis, we further study whether the proposed method \eqref{full-spde} can be applied to approximating  the invariant measure of the considered SPDE.
In many physical applications, the approximation of  the invariant measure is of fundamental importance,  especially when the invariant measure of the original system is unknown. 
For the results on numerically approximating the invariant measures of SPDEs, we refer to \cite{Bre14,HWZ17} and references therein.
For instance, the authors in \cite{HWZ17} consider the invariant measure of a full discretization and study the error 
between the invariant measure of the semi-discretization and that of the full discretization for the stochastic nonlinear Schr\"odinger equation.
The author in \cite{Bre14} investigates the error between the invariant measure for the temporal  semi-implicit method and the invariant measure of parabolic SPDE with Lipschitz and regular coefficients. 
Nevertheless, it is still not well known how to numerically approximate the invariant measures of  parabolic SPDEs with non-global Lipschitz coefficients and  how
to estimate the error  between these invariant measures.

To  solve these problems, we present the time-independent weak convergence analysis of
the proposed full discretization \eqref{full-spde},
which is much more involved than the time-dependent case. 
The main difficulties lie on deducing time-independent a priori estimations of numerical solution and showing the time-independent regularity estimates of Kolmogorov equation with respect to the spectral Galerkin approximation.
To this end, we introduce  the strong dissipative condition and the non-degenerate condition (see Section 4 for details).
Under the strong dissipative condition,  the time-independent regularity estimates of Kolmogorov equation are obtained by using  decay estimates.
Under the non-degenerate condition, 
we first study the $V$-uniform ergodicity of the invariant measure of the spectral Galerkin approximation. Then  the time-independent  regularity estimates of the corresponding Kolmogorov equation is established by  means of  the
Bismut--Elworthy--Li formula.
Finally,  the following result is proven.

\begin{tm}\label{weak-dis0}
Let Assumptions \ref{as-lap}-\ref{as-dri} hold with $\beta \in (0,1]$, $\gamma\in (0,\beta)$,
$X_0\in \HH^{\frac d2+\epsilon}$ with a sufficient small positive constant $\epsilon$, and $\delta t_0\in 
(0,1\land \frac 1{(2\lambda_F-2\lambda_1)\lor 0})$.
In addition, under Assumption \ref{erg-as} or \ref{erg-as1}, for any $\phi\in\mathcal C_b^2(\HH)$, 
there exists $C(X_0,Q,\phi)>0$ such that 
for 
$\delta t\in (0, \delta t_0]$, $K\ge 2$ and $N\in \N^+$,
\begin{align*}
|\E[\phi(X(K\delta t,X_0))-\phi(X_K^N(X_0^N))]|
&\le C(X_0,Q,\phi)(1+(K\delta t)^{-\gamma})(\delta t^{\gamma}+\lambda_N^{-\gamma}).
\end{align*}
\end{tm}

Then, as a consequence of Theorem \ref{weak-dis0}, the approximate error of the invariant measure for the proposed full discretization \eqref{full-spde} is obtained through the weak convergence analysis and the exponential ergodicity of 
Eq. \eqref{spde}. 

\begin{cor}\label{erg-cor}
Under the same conditions of Theorem \ref{weak-dis0}, for any  $\phi\in\mathcal C_b^2(\HH),$
there exist constants $c>0$, $C(X_0,Q,\phi)>0$ such that for any large $K$, $\delta t\in (0, \delta t_0]$
and $N\in \N^+$,
\begin{align*}
\Big|\E\Big[\phi(X_K^N(X_0^N))-\int_{\HH}\phi d\mu\Big]\Big|
&\le C(X_0,Q,\phi)(\delta t^{\gamma}+
\lambda_N^{-\gamma}+e^{-cK\delta t}),
\end{align*}
where $\mu$ is the unique invariant measure of 
Eq. \eqref{spde}.
Furthermore, if $\mu^{N,\delta t}$ is an ergodic invariant measure of the numerical  
solution $\{X^N_k\}_{k\in \N^+}$,
we have 
\begin{align*}
\Big|\E\Big[ \int_{P^N(\HH)}\phi d \mu^{N,\delta t}-\int_{\HH}\phi d\mu\Big]\Big|
&\le C(X_0,Q,\phi)(\delta t^{\gamma}+\lambda_N^{-\gamma}).
\end{align*}
\end{cor}
To the best of our knowledge, these are the first results on the time-independent weak error analysis and the convergence rates of the invariant measures of numerical methods for parabolic SPDEs with non-globally Lipschitz coefficients, especially for Eq. \eqref{spde}.

The remainder of this paper is organized as follows. 
In Section \ref{sec-3} we first introduce  
 some notations and assumptions. Then we propose the full discretization \eqref{full-spde} and  present both the regularity estimates and a priori estimates of the numerical solution, as well as the a priori estimates of semi-discretized stochastic convolution.
In Section \ref{sec-4}, we use the splitting based regularizing procedure and Malliavin calculus  to study the weak convergence rate of the proposed full discretization \eqref{full-spde}.
In Section \ref{sec-5} we show the regularity estimates of Kolmogorov equation with respect to the spectral Galerkin approximation, and use \eqref{full-spde}  to approximate the invariant measure of Eq. \eqref{spde} based on the time-independent weak error analysis.
Finally,  numerical tests are presented to
verify our theoretical results.

\section{ Preliminaries and full discretization}\label{sec-3}

In this section, we  give some basic assumptions and preliminaries, and introduce 
the spatial spectral Galerkin method 
and the implicit Euler type full discretization.
Furthermore, we present both the strong convergence and some a priori estimates for the proposed method.

\subsection{Preliminaries and assumptions}
Let $(\mathcal H,|\cdot|_{\mathcal H})$ and $(\widetilde H, \|\cdot\|_{\widetilde H})$  be separable Hilbert spaces.
We denote $\mathcal C_b^k(\mathcal H,\R),$ $k\in \N^+,$
the space of $k$ times continuous differentiable functionals from $\mathcal H$ to $\R$ with bounded derivatives up to order $k$, and
$B_b(\mathcal H,\R)$ 
the space of measurable and bounded functionals.
 Define 
\begin{align*}
\|\phi\|_{0}:=\sup\limits_{x\in \mathcal H}|\phi(x)|,\quad
|\phi|_{1}:=\sup\limits_{x\in \mathcal H}|D\phi(x)|_{\mathcal H},\quad
|\phi|_{2}:=\sup\limits_{x\in \mathcal H}|D^2\phi(x)|_{\mathcal L(\mathcal H)}
\end{align*}
with $D^k\phi$, $k=1,2$, being the $k$-th derivative of $\phi\in \mathcal C_b^k(\mathcal H,\R),$ and $\mathcal L(\mathcal H)$ being the space of linear operators from $\mathcal H$ into itself.
Denote by $\LL_2(\mathcal H,\widetilde H)$ the space 
of Hilbert--Schmidt operators from $\mathcal H$ into $\widetilde H$, equipped with the norm $\|\cdot\|_{\LL_2(\mathcal H,\widetilde H)}=(\sum_{k\in \N^+}\|(\cdot) f_k\|^2_{\widetilde H})^{\frac{1}{2}}$, where $\{f_k\}_{k\in \N^+},$ is an any 
orthonormal basis of $\mathcal H$. 
Given a Banach space $(\mathcal E,\|\cdot\|_{\mathcal E})$, we denote by $\gamma( \mathcal H, \mathcal E)$ the space of $\gamma$-radonifying operators endowed with the norm
$\|\cdot \|_{\gamma(\mathcal  H, \mathcal E)}=(\widetilde \E\|\sum_{k\in\N^+} (\cdot) f_k \gamma_k \|^2_{\mathcal E})^{\frac 12}$,
where $(\gamma_k)_{k\in\N^+}$ is a Rademacher sequence on a
probability space $(\widetilde \Omega,\widetilde \FFF, \widetilde \PP)$.

Moreover, we define $\HH:=L^2(\mathcal O)$ endowed with the norm $\|\cdot \|$ and the inner product $\<\cdot,\cdot\>,$ and denote $\mathcal L_2^0:=\mathcal L_2(\HH, U_0)$ with $U_0:=Q^{\frac 12}(\HH),$
where $Q\in \mathcal L(\HH)$ is self-adjoint and positive.
We also use the notation $\mathcal C_b^k(\HH):=\mathcal C_b^k(\HH,\R)$, $k\in \N^+$.
Meanwhile, let $\mathcal L:=\mathcal L(\HH),$
$E:=\mathcal C(\mathcal O;\R),$  
 $L^q:=L^q(\mathcal O; \R), q\ge 1$ equipped with the norm $\|\cdot \|_{L^q},$ and $H^{k}$, $H^{k}_0, k\in \N^+$ be the usual Sobolev spaces equipped with usual norms.

Let $\mathcal  I: L^2([0,T]; U_0) \to L^2(\Omega)$ be an isonormal process, i.e, $\mathcal I(\psi)$ is the centered Gaussian random variable, for any $\psi \in L^2([0,T]; U_0),$
and $\E[\mathcal I(\psi_1)\mathcal I(\psi_2)]=\<\psi_1,\psi_2\>_{L^2([0,T]; U_0)},$ for any $\psi_1,\psi_2\in L^2([0,T]; U_0).$
We denote the family of smooth real-valued cylindrical random variables by
\begin{align*}
\mathcal S=\Big\{\mathcal X=g(\mathcal I(\psi_1),\cdots,\mathcal I(\psi_n)): g\in \mathcal C_p^{\infty}(\R^n), \psi_j \in L^2([0,T]; U_0), j=1,\cdots,n \Big\},
\end{align*}
where $\mathcal C_p^{\infty}(\R^n)$ is the space of all real-valued $\mathcal C^{\infty}$ functions on $\R^n$ with polynomial growth,
and the family of smooth cylindrical $\HH$-valued 
random variables by
\begin{align*}
\mathcal S(\mathbb H)=\Big\{G=\sum_{i=1}^M\mathcal X_i\otimes h_i: \mathcal X_i \in \mathcal S, h_i\in \HH, M\ge 1\Big\},
\end{align*} 
where $\otimes$ denotes the tensor product.
For $G=\sum\limits_{i=1}^Mg_i(\mathcal I(\psi_1),\cdots,\mathcal I(\psi_n))\otimes h_i,$ define its Malliavin derivative 
\begin{align*}
\mathcal D_sG=\sum_{i= 1}^M\sum_{j=1}^n\partial_j g_i(\mathcal I(\psi_1),\cdots,\mathcal I(\psi_n))\otimes (h_i\otimes \psi_j(s)),
\end{align*}
where $s\in [0,T]$.
 Let $\mathbb D^{1,2}(\HH)$  be the closure of $\mathcal S(\HH)$ under the norm 
\begin{align*}
\|G\|_{\mathbb D^{1,2}(\HH)}=\Big(\E[\|G\|^2]+\E[\int_0^T\|\mathcal D_sG\|^2ds]\Big)^{\frac 12}.
\end{align*}
Then the Malliavin integration by parts formula holds (see, e.g., \cite[Section 2]{Deb11}),
namely, for any random variable $G\in \mathbb D^{1,2}(\HH)$
and any predictable  process $\Theta \in L^2([0,T];\mathcal L_2^0)$, we have 
\begin{align}
\label{int-by}
\E\Big[\Big\<\int_0^T\Theta (t)dW(t),G\Big\>\Big]
=\E\Big[ \int_0^T \Big\<\Theta (t),\mathcal D_tG\Big\>_{\mathcal L_2^0}dt\Big].
\end{align}
This property is the key to analyzing the weak convergence rate of numerical method in Sections \ref{sec-4} and \ref{sec-5}.
Additionally, the Malliavin derivative satisfies the chain rule, that is, for $\sigma(G)\in \mathbb D^{1,2}(\mathcal H),$ 
\begin{align*}
\mathcal D_t^y(\sigma(G))&=\mathcal D\sigma(G)\cdot\mathcal D_t^yG, \quad y\in U_0, \quad G\in \mathbb D^{1,2}(\HH),\\
\mathcal D_t(\sigma(G))&=\mathcal D\sigma(G)\mathcal D_tG, \quad G\in \mathbb D^{1,2}(\HH),
\end{align*} 
where $\mathcal D_t^yG:=\mathcal D_tG y$ is the derivative of  $G$ in the direction of $y\in U_0$, $\sigma \in \mathcal C_b^1(\HH,\mathcal H).$

Throughout this article, we use $c$, $C$ to denote generic constants, independent of  $N$ and $\delta t$, which may differ from one place to another. 
Unless otherwise specified, we always assume that $X_0$ is a deterministic function in $\HH^{\frac d2+\epsilon}$, where $\epsilon$ is sufficiently small positive number.
In the following, we introduce some assumptions on  both the coefficients and driving noises for Eq. \eqref{spde}.

\begin{as}\label{as-lap}
Let $\mathcal O:=[0,L]^d, d\le 3, L>0$.
Let $A: D(A)\subset \HH \to \HH$ be the Laplacian operator on $\mathcal O$ with the homogenous Dirichlet boundary condition, i.e., $Au=\Delta u, u\in D(A)$.
\end{as}

This assumption implies that the operator
$A$ generates 
an analytic and contraction $C_0$-semigroup $S(t), t\ge 0$ in $\HH$ and $L^q$, $q\ge 1$ and  that the existence of the eigensystem $\{\lambda_k, e_k\}_{k\in\N^+}$ of 
$\HH$, such that $\{\lambda_k\}_{k\ge1}$ is an increasing sequence, $-Ae_k=\lambda_k e_k$, $\lim\limits_{k\to \infty}\lambda_k=\infty$ and $\sup\limits_{k\in \N^+}\|e_k\|_E\le C$.
Let $\mathbb H^r$ be the Banach space equipped with the norm
 $\|\cdot \|_{\HH^r}:= \|(-A)^{\frac {r}2}\cdot \|$ 
for  the fractional power $(-A)^{\frac{r}{2}}, r\ge 0$.
We also remark that Assumption \ref{as-lap} can be 
extended to the case that $A$ is a second order elliptic operator on a regular domain and a part of A in $E$ generates an analytic semigroup in $E$.

\begin{as}\label{as-noi}
Let $W(t)$ be a Wiener process with covariance operator $Q$, where $Q$ is a bounded, linear, self-adjoint and positive definite operator on $\HH$ and satisfies
$\|(-A)^{\frac {\beta-1}2 }\|_{\LL_2^0}<\infty$ with $0< \beta\le 2$. 
Assume in addition that $\beta>\frac d2$ or $A$ commutes with $Q$.
\end{as}

In the case of investigating the strong error estimate,  the additional condition that $\beta>\frac d2$ or $A$ commutes with $Q$ can be weakened. 
The additional assumption is used to ensure a priori estimates of exact and numerical solutions when studying the weak convergence rates of numerical methods.
In order to get the time-independent error estimate and to approximate the invariant measure, some dissipative condition and non-degenerate condition 
are presented in Section \ref{sec-5}.

\begin{as}\label{as-dri}
Let $f$ be a cubic polynomial with $f(\xi)=-a_3\xi^3+a_2\xi^2+a_1\xi+a_0$,  $a_i\in \R$, $i=0,1,2,3$, $a_3>0$ and 
let $F: L^{6} \to \HH$ be the Nemytskii operator defined by $F(X)(\xi)=f(X(\xi))$.
\end{as}

Notice that $\lambda_F=\sup\limits_{\xi \in \R}f'(\xi).$ 
The above assumption ensures that $F$ satisfies 
\begin{align*}
\<F(u)-F(v),u-v\>&\le \lambda_F \|u-v\|^2,\\
\|F(u)-F(v)\|&\le C_f(1+\|u\|_{E}^{2}+\|v\|_{E}^2)\|u-v\|,
\end{align*}
for $C_f>0$.
In the case that $f(\xi)=-\xi^3+\xi$, Eq. \eqref{spde} corresponds to the  stochastic Allen--Cahn equation or stochastic Ginzburg--Landau equation.

\subsection{Full discretization}
Now we are in a position to give both the semi-discretization and the full discretization for Eq. \eqref{spde}.
In the sequel, we let $\delta t\in (0,\delta t_0]$, $\delta t_0\in \R^+$, $t_k$=$k\delta t$, $k\in \N^+$ and $N\in \N^+$. The notations $[t]_{\delta t}:=\max \big(\{0,\delta t,\cdots,k\delta t,\cdots \}\cap [0,t]\big)$, $\lfloor t \rfloor=\frac {[t]_{\delta t}}{\delta t}$ are used frequently. 

Denote the spectral Galerkin projection by $P^N$. Using spectral Galerkin method in 
space, we get the following semi-discretization
\begin{align}
\label{semi-spde}
dX^N(t)=AX^N(t)dt+P^NF(X^N(t))dt+P^NdW(t).
\end{align}
For the  weak convergence analysis in a finite interval $[0,T]$, we choose a proper $\delta t$ such that $K\delta t=T$ for some $K \in \N^+$.
For the time-independent weak convergence analysis, we fixed the stepsize $\delta t$ and let $K\in \N^+$.
By applying implicit Euler method to discretize Eq. \eqref{semi-spde} further, we get the
full discretization  
\begin{align*}
X^N_{k+1}=X^N_k+\delta tAX^N_{k+1}
+\delta tP^NF(X^N_{k+1})+ P^N\delta W_{k},
\end{align*}
where $X^N_0=P^NX_0$ and $P^N\delta W_{k}=P^N(W((k+1)\delta t)-W(k\delta t)).$
Here, for the sake of simplicity, we omit the dependence on the initial data and denote  $X^N_k:=X^N_k(X_0^N)$, $k \in \N^+$.
Denoting $S_{\delta t}:=(I-A\delta t)^{-1}$, then the
full discretization can be rewritten as
\begin{align}\label{full-spde}
X^N_{k+1}=S_{\delta t}X^N_k
+\delta tS_{\delta t}P^NF(X^N_{k+1})
+S_{\delta t}P^N\delta W_{k},
\end{align}
which is equivalent to $X^N_{k+1}=Y^N_{k+1}+Z^N_{k+1}$ with
\begin{align}\label{eq-ran}
Y^N_{k+1}&=Y^N_k+\delta t AY^N_{k+1}
+\delta tP^NF(Y^N_{k+1}+Z^N_{k+1}),\\\nonumber
Z^N_{k+1}&=\sum_{j=0}^k S_{\delta t}^{k+1-j}P^N \delta W_j.
\end{align}

Moreover, based on Assumption \ref{as-dri}, the solvability of the proposed method \eqref{full-spde} is obtained if the time stepsize $\delta t$ is small. 
Indeed,  if $\delta t_0< 1\land \frac 1{(2\lambda_F-2\lambda_1)\lor 0},$ then the proposed method has a unique solution when $\delta t\in (0,\delta t_0]$.

Let Assumptions \ref{as-lap}-\ref{as-dri} hold 
with $\beta >\frac d2$ and $X_0\in \HH^{\beta}$, 
or with $A$ commuting with $Q$ and $X_0\in \HH^{\beta}\cap E$.  It can be shown (see e.g. \cite[Chapter 6]{Cer01}) that for any $T>0$, $$\sup\limits_{s\in[0,T]}\E\Big[\|X(s)\|_{E}^p\Big]+\sup\limits_{s\in[0,T]}\E\Big[\|X(s)\|^p_{\HH^{\beta}}\Big]\le C(X_0,Q,p)$$ and $$\E\Big[\|X(t)-X(s)\|^p\Big]\le C(X_0,Q,p)(t-s)^{\frac {\beta p}2},$$ where $p\ge 1$, $0 \le s<t\le T$.
Furthermore, the following strong error estimate holds. Its proof is similar to the proofs of \cite[Theorem 4.1]{QW18}  and \cite[Theorem 3.1]{CH18}.
\begin{lm}\label{str}
Under  Assumptions \ref{as-lap}-\ref{as-dri},
let $X_0\in \HH^{\beta}$, $T=K\delta t$, $K\in \N^+$ and $N\in \N^+$. Then 
the  full discretization
is strongly convergent and satisfies
\begin{align*}
\sup_{k\le K}\Big\|X_k^N-X(t_k)\Big\|_{L^p(\Omega;\HH)}
&\le C(T,X_0,Q)(\delta t^{\frac \beta 2}+\lambda_N^{-\frac \beta 2}).
\end{align*}
\end{lm}
\begin{rk}\label{rk-str}
Similar to Lemma \ref{str}, for $T=K\delta t>0$, $k\le K$, $k\in \N$ and $N\in \N^+$, we have
\begin{align*}
\sup_{k\le K}\Big\|X_k^N-X^N(t_k)\Big\|_{L^p(\Omega;\HH)}
&\le C(T,X_0,Q)\delta t^{\frac \beta 2},
\end{align*}
and 
\begin{align*}
\sup_{t\in [0,T]}\Big\|X^N(t)-X(t)\Big\|_{L^p(\Omega;\HH)}
&\le C(T,X_0,Q)\lambda_N^{-\frac \beta 2}.
\end{align*}
\end{rk}

Lemma \ref{str} yields a upper bound on the weak convergence error estimate. 
Combining with weak convergence result in Theorem \ref{weak0}, we immediately obtain that  
the weak convergence rate is
$\mathcal O\left(\delta t^{(1\land \beta-\epsilon_1)\lor \frac \beta 2}+\lambda_N^{-(1\land\beta-\epsilon_1)\lor  \frac \beta 2)}\right)$ where $\epsilon_1$ is a sufficient small positive number.
Thus in Sections  \ref{sec-3}-\ref{sec-5}, we mainly focus on weak convergence rates of numerical methods in the case that $\beta\in (0,1]$. We would like to mention that for SPDEs with non-globally Lipschitz coefficients, there already exist a lot of results on the strong convergence and strong convergence rate of numerical approximations, see \cite{AC17,BGJK17,BJ16,BCH18,BG18,CHL16b,CHLZ17,CHS18a,KLL18,LQ18,Wang18} and references therein.

\subsection{A priori estimate of the  full discretization}
In this subsection, we present the  time-independent a priori estimate of the proposed
numerical method. The following lemma is about the a priori estimate of
the solution $X^N$ for the  
spectral Galerkin method, which is very useful in Section \ref{sec-5}. Its proof is similar to  that of 
the numerical solution, see  
Lemmas \ref{pri-z} and \ref{pri-y}.

\begin{lm}\label{pri-xn}
Let Assumptions \ref{as-lap}-\ref{as-dri} hold, $p\ge1$. For $\gamma\in (0,\beta]$, 
there exist  $C(X_0,Q,p)>0$ and $C( X_0,Q,p,\gamma)>0$ such that 
\begin{align*}
\sup_{t\ge0}\E\Big[\|X^N(t,X_0)\|_{\HH^{\gamma}}^p\Big]\le  C(X_0,Q,p,\gamma) \;\;\text{and} \;\;\sup_{t\ge 0}\E\Big[\|X^N(t,X_0)\|_E^p\Big]\le C( X_0,Q,p).
\end{align*}
\end{lm}

We also need the uniform bound of $X^N_k$, $k\le K$, for the full discretization \eqref{full-spde}. To this end, it suffices to show  a priori  estimates of $Y^N_k$ and $Z^N_k$, $k\in \N^+$.

\begin{lm}\label{pri-z}
Let  Assumptions \ref{as-lap}-\ref{as-noi} hold and  $p\ge 1$. There exists $C(Q,p)>0$ such that the discretized stochastic convolution $\{Z^N_k\}_{k\in \N^+}$ satisfies 
\begin{align*}
\sup_{k\in\N^+}\|Z^N_k\|_{L^p(\Omega;E)}\le C(Q,p).
\end{align*}
\end{lm}
\textbf{Proof}
Under the condition  that $A$ commutes with $Q$, we apply the  Burkholder inequality and get 
\begin{align*}
\|Z^N_{k+1}\|_{L^p(\Omega;E)}&=\Big\|\sum_{j=0}^k S_{\delta t}^{k+1-j}P^N \delta W_j\Big\|_{L^p(\Omega;E)}\\
&\le 
\|S_{\delta t}^{k+1-\lfloor \cdot \rfloor}\|_{L^p(\Omega;L^2([0,t_{k+1}];\gamma(\HH;E)))}\\
&\le C\sqrt{\sum_{l\in \N^+}\sum_{j=0}^{k}\Bigg(\frac {1}{1+\lambda_l\delta t}\Bigg)^{2(k+1-j)}q_l\delta t}\\
&\le  C\sqrt{\sum_{l\in \N^+}\frac {1}{\lambda_l(2+\lambda_l\delta t)}\Bigg(1-(\frac {1}{1+\lambda_l\delta t})^{2(k+1)}\Bigg)q_l}\\
&\le C\sqrt{\sum_{l\in \N^+}\lambda_l^{\beta-1}q_l
\sup_{l\in\N^+}
\frac 1{\lambda_l^\beta(2+\lambda_l\delta t
)}}\le C(Q,p),
\end{align*}
where $\{q_l\}_{l\ge 1}$ is the sequence of  eigenvalues of $Q$.
In another case that $\beta >\frac d2$, it follows  from the Sobolev embedding theorem and the Burkholder--Davis--Gundy inequality that  for a sufficient small number $\epsilon>0$,
\begin{align*}
\|Z^N_{k+1}\|_{L^p(\Omega;E)}&\le \Big\|\sum_{j=0}^k S_{\delta t}^{k+1-j}P^N \delta W_j\Big\|_{L^p(\Omega; \HH^{\frac d2+\epsilon})}\\
&\le C\sqrt{\sum_{j=0}^k\left\|(-A)^{\frac d4+\frac \epsilon2}S_{\delta t}^{k+1-j}Q^{\frac 12}\right\|_{\LL_2}^2\delta t}\\
&\le C\sqrt{\sum_{j=0}^k\left\|(-A)^{\frac {\frac d2+\epsilon+1-\beta}2}S_{\delta t}^{k+1-j}\right\|_{\LL}^2\|(-A)^{\frac {\beta-1}2}\|_{\LL_2^0}^2 \delta t}
\\
&\le C\sqrt{\sum_{j=0}^k \frac 1{((k+1-j)\delta t)^{\frac d2+\epsilon+1-\beta}}\frac 1{(1+\lambda_1\delta t)^{(k+1-j)(1+\beta-\frac d2 -\epsilon)}}\delta t}\\
&\le C\sqrt{\int_0^{\infty}t^{-\frac d2+\beta-\epsilon-1}\frac 1{(1+\lambda_1\delta t)^{(1+\beta-\frac d2-\epsilon)\lfloor t\rfloor}}dt} 
\le C(Q,p).
\end{align*}
\qed

\begin{lm}\label{pri-y}
Under Assumptions \ref{as-lap}-\ref{as-dri}, there exists  $C(Q,X_0,p)>0$ such that  the solution $\{Y^N_k\}_{k\in \N^+}$ of Eq. \eqref{eq-ran}  satisfies 
\begin{align*}
\sup_{k\in \N^+}\|Y^N_k\|_{L^p(\Omega;E)}\le C(Q,X_0,p).
\end{align*}
\end{lm}
\textbf{Proof}
By multiplying $Y^N_{k+1}$ on both sides of Eq. \eqref{eq-ran}, integrating over $\mathcal O$ and the Young inequality, we have 
\begin{align*}
\frac 12\|Y_{k+1}^N\|^2
&\le \frac 12\|Y_k^{N}\|^2-\delta t\|\nabla Y_{k+1}^N\|^2
- (a_3-\epsilon) \delta t\|Y_{k+1}^N\|_{L^4}^4
+C(\epsilon)\delta t(1+\|Z_{k+1}^N\|_E^4)\\
&\le \frac 12\|Y_k^{N}\|^2-
\lambda_1\delta t\|Y_{k+1}^N\|^2
+C(\epsilon)\delta t(1+\|Z_{k+1}^N\|_E^4).
\end{align*}
The Gronwall inequality, together with the a priori estimate of $Z_k^N$ in Lemma \ref{pri-z}, leads to 
\begin{align*}
\E\Big[\|Y_{k+1}^N\|^2\Big]&\le C(Q,X_0)\sum_{j=0}^{k}\frac 1{(1+2\lambda_1\delta t)^{k+1-j}}\delta t\le C(Q,X_0).
\end{align*}
Next, we turn to estimate the a priori estimate in $E$ by the mild form of $$Y^N_{k+1}=S_{\delta t}^{k+1}Y_0^N+\sum_{j=0}^{k}S_{\delta t}^{k+1-j}P^NF(Y_{j+1}^N+Z_{j+1}^N)\delta t.$$
The Sobolev embedding theorem and the  smoothing effect of $S_{\delta t}$  yield that 
\begin{align*}
&\|Y^N_{k+1}\|_E\\
&\le \|Y_0^N\|_E+C\sum_{j=0}^{k}\|(-A)^{\frac d4+\epsilon}S_{\delta t}^{k+1-j}P^NF(Y_{j+1}^N+Z_{j+1}^N)\|\delta t\\
&\le \|Y_0^N\|_E+C\sum_{j=0}^k\frac 1{((k+1-j)\delta t)^{\frac d4+\epsilon}}\frac 1{(1+\lambda_1\delta t)^{(k+1-j)(1-\frac d4-\epsilon)}}\delta t
\|F(Y_{j+1}^N+Z_{j+1}^N)\|\\
&\le \|Y_0\|_{\HH^{\frac d2+\epsilon}}+C\sum_{j=0}^k\frac 1{((k+1-j)\delta t)^{\frac d4+\epsilon}}\frac 1{(1+\lambda_1\delta t)^{(k+1-j)(1-\frac d4-\epsilon)}}\delta t\\
&\qquad(1+\|Y_{j+1}^N\|_{\HH^1}^3+\|Z_{j+1}^N\|_{L^6}^3).
\end{align*}
The above estimate can be improved in $d=1$ by using Gagliard--Nirenberg inequality $\|u\|_{L^6}\le C\|\nabla u\|^{\frac 13}\|u\|^{\frac 23}$ and 
the estimation of $\sum_{j=0}^k\|\nabla Y_{j+1}\|^2$ (see the proof \cite[Proposition 3.1]{CH18}).
For  $d=2,3$,
we need to give the a priori estimate of $\|\nabla Y_{k+1}^N\|$.

Multiplying  the term $-A Y_{k+1}^N$ on both sides of  Eq. \eqref{eq-ran} and  integrating over $\mathcal O,$ we obtain 
\begin{align*}
\|\nabla Y_{k+1}^N\|^2
&\le \|\nabla Y_{k}^N\|^2 -2\delta t\|A Y_{k+1}^N\|^2+2\delta t\<\nabla F(Y_{k+1}^N+Z_{k+1}^N),\nabla Y_{k+1}^N\>
\\
&\le \|\nabla Y_{k}^N\|^2-(2-\epsilon)\delta t\|A Y_{k+1}^N\|^2\\
&\quad +C(\epsilon)\delta t(\|Y_{k+1}^N\|^2+\|Z_{k+1}^N\|^2+\|Z_{k+1}^N\|_{L^6}^6)\\
&\quad -2a_3\delta t\<Y_{k+1}^N\nabla Y_{k+1}^N,Y_{k+1}^N\nabla Y_{k+1}^N\> 
+C\delta t\<Y_{k+1}^N(Z_{k+1}^N)^2,A Y_{k+1}^N\>\\
&\quad+C\delta t\<(Y_{k+1}^N)^2Z_{k+1}^N,A Y_{k+1}^N\>
+C\delta t\<(Y_{k+1}^N+Z_{k+1}^N)^2,A Y_{k+1}^N\>\\
&\le  \|\nabla Y_{k}^N\|^2+C(\epsilon)\delta t(1+\|Y_{k+1}^N\|^2+\|Z_{k+1}^N\|_{E}^6+\|Y_{k+1}^N\|^2\|Z_{k+1}^N\|_E^4)\\
&\quad
-(2-2\epsilon)\delta t\|A Y_{k+1}^N\|^2
+C\delta t\|Z_{k+1}^N\|_E^2\|Y_{k+1}^N\|_{L^4}^4.
\end{align*}
The Gagliardo--Nirenberg--Sobolev inequality yields that 
\begin{align*}
\|\nabla Y_{k+1}^N\|^2
&\le \|\nabla Y_{k}^N\|^2+C(\epsilon)\delta t(1+\|Y_{k+1}^N\|^2+\|Z_{k+1}^N\|_{E}^6+\|Y_{k+1}^N\|^2\|Z_{k+1}^N\|_E^4)\\
&\quad
-(1-2\epsilon)\delta t\|A Y_{k+1}^N\|^2
+C\delta t\|Z_{k+1}^N\|_E^2\|A Y_{k+1}^N\|^{\frac d2}\|Y_{k+1}^N\|^{4-\frac d2}\\
&\le  \|\nabla Y_{k}^N\|^2+C(\epsilon)\delta t(1+\|Y_{k+1}^N\|^2+\|Z_{k+1}^N\|_{E}^6+\|Y_{k+1}^N\|^2\|Z_{k+1}^N\|_E^4)\\
&\quad
-(2-3\epsilon)\delta t\|A Y_{k+1}^N\|^2
+C(\epsilon)\delta t\|Z_{k+1}^N\|_E^{\frac 8{4-d}}\|Y_{k+1}^N\|^{\frac {16-2d}{4-d}}.
\end{align*}
Combining the estimations of $Y_{k+1}^N$ and $ \nabla Y_{k+1}^N$ with the equivalence of the norm in $\HH^2$ and the norm in $H^1_0\cap H^2$, we obtain
\begin{align*}
\|Y_{k+1}^N\|_{\HH^1}^2
&\le \|Y_{k}^N\|_{\HH^1}^2
-c(2-3\epsilon)\delta t \|Y_{k+1}^N\|_{\HH^2}^2+
C(\epsilon)\delta t(1+\|Y_{k+1}^N\|^2+\|Z_{k+1}^N\|_{E}^6\\
&\quad+\|Y_{k+1}^N\|^2\|Z_{k+1}^N\|_E^4)
+C\delta t\|Z_{k+1}^N\|_E^{\frac 8{4-d}}\|Y_{k+1}^N\|^{\frac {16-2d}{4-d}}\\
&\le 
 \|Y_{k}^N\|_{\HH^1}^2
-c(2-3\epsilon)\delta t \|Y_{k+1}^N\|_{\HH^1}^2+
C(\epsilon)\delta t(1+\|Y_{k+1}^N\|^2+\|Z_{k+1}^N\|^2\\
&\quad+\|Z_{k+1}^N\|_{E}^6+\|Y_{k+1}^N\|^2\|Z_{k+1}^N\|_E^4)
+C\delta t\|Z_{k+1}^N\|_E^{\frac 8{4-d}}\|Y_{k+1}^N\|^{\frac {16-2d}{4-d}}.
\end{align*}
By using Gronwall's inequality and then taking the $p$-th moment on both sides, combining with the a priori estimate of $\|Z_{k}^N\|_E$ in Lemma \ref{pri-z} and $\|Y_{k}^N\|$, we complete the proof.
\qed

By a more refined estimate, the following result holds.

\begin{cor}\label{pri-xnk}
Under Assumptions \ref{as-lap}-\ref{as-dri}, for $p\ge 1$,
there exists $C(Q,X_0,p)>0$ such that
\begin{align*}
\left\|\sup_{k\in \N^+}\|Y^N_k\|_E\right\|_{L^p(\Omega;\R)}\le C(Q,X_0,p),
\end{align*}
Moreover, it holds that
\begin{align*}
\sup_{k\in \N^+}\left\|X^N_k\right\|_{L^p(\Omega;E)}\le C'(Q,X_0,p),
\end{align*}
where $C'(Q,X_0,p)>0$.
\end{cor}

Apart from the a priori estimate of $X_{k}^N$,  the Malliavin regularity of the numerical method is needed to control the stochastic integral error term in the weak convergence analysis in Sections \ref{sec-4} and \ref{sec-5}.

\begin{prop}\label{mal-dif}
Let  Assumptions \ref{as-lap}-\ref{as-dri} hold with $\beta\in (0,1]$, $p\ge 1$.
Then  there exist $C(Q,X_0,p)>0$ and $C_1>1$ such that for  $s<t_{k+1}$, $k\in \N$, $z\in U_0$, 
\begin{align*}
\|\mathcal D_s^z X_{k+1}^N\|_{L^p(\Omega, \HH)}
&\le C(Q,X_0,p) \big(1\lor ({1+2C_1(\lambda_1-\lambda_F)\delta t})^{k+1-\lfloor s\rfloor}\big) \\
&\quad \Big(1+t_{k+1-\lfloor s\rfloor }^{\frac {\beta -1}2}\frac 1{(1+\lambda_1\delta t)^{(k+1-\lfloor s\rfloor)\frac{\beta+1}2}}\Big) \|(-A)^{\frac {\beta -1}2}z\|, 
\end{align*}
and 
\begin{align*}
\|(-A)^{\frac {\beta-1}2} \mathcal D_s^z X_{k+1}^N\|_{L^p(\Omega, \HH)}
&\le C(Q,X_0,p)(1\lor ({1+2C_1(\lambda_1-\lambda_F)\delta t})^{k+1-\lfloor s\rfloor})\\
&\quad \|(-A)^{\frac {\beta-1}2}z\|. 
\end{align*}
\end{prop}
\textbf{Proof}
For $s\ge (k+1)\delta t $, $z\in U_0$, we have  $\mathcal D_s^z X^N_{k+1}=0$.
For $0\le s < k\delta t\le T$, $z\in U_0$, we obtain
\begin{align*}
\mathcal D_s^z X^N_{k+1}
&= S_{\delta t}\mathcal D_s^z X^N_k
+\delta tS_{\delta t}P^N(DF(X^N_{k+1})\cdot \mathcal D_s^zX^N_{k+1} )\\
&= \mathcal D_s^z X^N_k+\delta tA \mathcal D_s^z X^N_{k+1}
+\delta tP^N(DF(X^N_{k+1})\cdot \mathcal D_s^zX^N_{k+1} ).
\end{align*}
For $k\delta t \le  s < (k+1)\delta t$, $z\in U_0$, we have
\begin{align*}
\mathcal D_s^z X^N_{k+1}
&= P^NS_{\delta t}z
+\delta tS_{\delta t}P^N(DF(X^N_{k+1})\cdot \mathcal D_s^z X^N_{k+1} ).
\end{align*}
From the above calculations,  it follows that 
for $k \ge [s]_{\delta t}$, $\lfloor s\rfloor\delta t\le s<(\lfloor s\rfloor+1)\delta t$,
\begin{align*}
\mathcal D_s^zX_{k+1}^N
&=P^NS_{\delta t}^{k+1-\lfloor s \rfloor}z
+\delta t\sum_{j= \lfloor s \rfloor}^{k}S_{\delta t}^{k+1-j}
P^N(DF(X_{j+1}^N)\cdot \mathcal D_s^zX_{j+1}^N).
\end{align*}
Then we show the regularity estimate of the Malliavin derivative $\mathcal D_s^zX_{k+1}^N$ by using similar arguments in \cite[Proposition 4.2]{CH18}.
Since in each step, $\mathcal D_s^zX_{k+1}^N$
can be viewed as 
\begin{align*}
\mathcal D_s^z X^N_{k+1}
&= \mathcal D_s^z X^N_k+\delta tA \mathcal D_s^z X^N_{k+1}
+\delta tP^N(DF(X^N_{k+1})\cdot \mathcal D_s^zX^N_{k+1} ),
\end{align*}
it follows that  
\begin{align*}
\|\mathcal D_s^z X^N_{k+1}\|^2
&\le \|\mathcal D_s^z X^N_{k}\|^2
-2\delta t \|\nabla D_s^z X^N_{k+1}\|^2
+2\delta t \<DF(X^N_{k+1})\cdot \mathcal D_s^zX^N_{k+1},\mathcal D_s^zX^N_{k+1} \>\\
&\le \|\mathcal D_s^z X^N_{k}\|^2
+2\delta t(-\lambda_1+\lambda_F)\|\mathcal D_s^z X^N_{k+1}\|^2,
\end{align*} 
which implies 
\begin{align*}
\|\mathcal D_s^z X^N_{k+1}\|^2
&\le \frac {1}{1+2(\lambda_1-\lambda_F)\delta t} \|\mathcal D_s^z X^N_{k}\|^2.
\end{align*}
Next we aim to estimate the regularity of $\mathcal D_s^z X^N_{k+1}$. By defining  
$V^N_z(k+1,s):=\mathcal D_s^z X^N_{k+1}-P^NS_{\delta t}^{k+1-\lfloor s\rfloor}z$, it follows that 
\begin{align*}
V^N_z(k+1,s)&=V^N_z(k,s)-\delta tA V_z^N(k+1,s)
+\delta t P^N(DF(X^N_{k+1}) \cdot V^N_z(k+1,s))\\
&\quad+ \delta t P^N(DF(X^N_{k+1})\cdot P^NS_{\delta t}^{k+1-\lfloor s\rfloor}z).
\end{align*}
By similar  arguments in the proof of  \cite[Proposition 4.1]{CH18}, we obtain 
\begin{align*}
&\|V^N_z(k+1,s)\|\\
&\le  \delta t \sum_{j=\lfloor s\rfloor}^k 
\Big(1+2(\lambda_1-\lambda_F)\delta t \Big)^{-(k+1-j)}\|DF(X^N_{j+1})\cdot S_{\delta t}^{j+1-\lfloor s\rfloor} P^N z\|
\\
&\le C\delta t \sum_{j=\lfloor s\rfloor}^k\Big(1+2(\lambda_1-\lambda_F)\delta t \Big)^{-(k+1-j)}(1+\|X^N_{j+1}\|^2_E) ((j+1-\lfloor s\rfloor)\delta t)^{-\alpha}\\
&\quad \times\frac 1{(1+2\lambda_1\delta t)^{(1-\alpha)(j+1-[s]_{\delta t})}}
\|(-A)^{-\alpha}z\|\\
&\le C\delta t \big(1\lor ({1+2C_1(\lambda_1-\lambda_F)\delta t}\big)^{k+1-\lfloor s\rfloor}) \sum_{j=[s]_{\delta t}}^k(1+\|X^N_{j+1}\|^2_E) ((j+1-\lfloor s\rfloor)\delta t)^{-\alpha}\\
&\quad \times\frac 1{(1+2\lambda_1\delta t)^{(1-\alpha)(j+1-\lfloor s\rfloor)}}
\|(-A)^{-\alpha}z\|,
\end{align*}
where $\delta t\le \frac 1{2(\lambda_F-\lambda_1)} (1-\frac 1{C_1})$ if $\lambda_F>\lambda_1$.
The smoothy effect of $S_{\delta t}$ leads to
\begin{align*}
\|\mathcal D_s^z X_{k+1}^N\|
&\le 
\|V^N_z(k+1,s)\|_{L^p(\Omega, \HH)}
+\|P^NS_{\delta t}^{k+1-\lfloor s\rfloor}z\|\\
&\le C(Q,p)(1+\sup_{j\in \N^+}\|X_{j}\|^2_E) \big(1\lor ({1+2C_1(\lambda_1-\lambda_F)\delta t})^{k+1-\lfloor s\rfloor}\big) \\
&\quad \times \Big(1+t_{k+1-\lfloor s\rfloor}^{-\alpha}\frac 1{(1+2\lambda_1\delta t)^{(k+1-\lfloor s\rfloor)(1-\alpha)}}\Big)\|(-A)^{-\alpha}z\|.
\end{align*}
Based on the above estimate, taking expectation and taking $\alpha = \frac {1 -\beta}2$,  we finish the proof of the first desired estimate.
Similar arguments in the proof of \cite[Proposition 4.2]{CH18} lead to the second desired estimate.
\qed

\begin{rk}
Under the condition of Proposition \ref{mal-dif},
if in addition $\lambda_1\ge \lambda_F$, 
the following time-independent estimates hold, i.e., 
\begin{align*}
\|\mathcal D_s X_{k+1}^N\|_{L^p(\Omega, \LL_2^0)}
&\le C(Q,X_0,p) 
\Big(1+t_{k+1-\lfloor s\rfloor}^{\frac {\beta -1}2}\frac 1{(1+2\lambda_1\delta t)^{(k+1-\lfloor s\rfloor)\frac{\beta+1}2}}\Big) \|(-A)^{\frac {\beta -1}2}\|_{\LL_2^0}
\end{align*}
and 
\begin{align*}
\|(-A)^{\frac {\beta-1}2} \mathcal D_s X_{k+1}^N\|_{L^p(\Omega, \LL_2^0)}
&\le C(Q,X_0,p)\|(-A)^{\frac {\beta-1}2}\|_{\LL_2^0}, 
\end{align*}
for some constant $C(Q,X_0,p)$.
\end{rk}

Based on the strong convergence, the a priori estimate  and the Malliavin regularity of numerical solutions for  \eqref{full-spde}, we are able to deal with the weak convergence of the proposed method in the next sections.

\section{Weak convergence analysis of the full discretization}\label{sec-4}
\label{sec-weak}
In this section, we aim to present the weak error analysis for the considered numerical method approximating Eq. \eqref{spde}. Following the idea of \cite{CH18}, we  introduce the auxiliary regularized stochastic PDE and its corresponding Kolmogorov equation.

Consider the auxiliary problem
\begin{align}\label{rspde}
dX^{\tau}=AX^{\tau}dt+\Psi_{\tau}(X^{\tau})dt
+dW(t),\quad 
X^{\tau}(0)=X_0,
\end{align}
where $\tau$ is regularizing parameter of this splitting approach, $\Psi_{t}(\xi):=\frac {\Phi_t(\xi)-\xi}t$, $t>0$
and $\Psi_{0}(\xi)=F(\xi)$, $\Phi_t$ is the phase flow of
the  differential equation 
\begin{align*}
dx(t)=f(x(t))dt,\quad x(0)=\xi \in \R.
\end{align*}

Next, we give the regularity estimate of Kolmogorov equation with respect to Eq. \eqref{rspde} shown in \cite{CH18},
\begin{align}\label{kol}
\frac {\partial U^{\tau}(t,x)} {\partial t}
&=\<Ax+\Psi_{\tau}(x), D U^{\tau}(t,x)\>+\frac 12tr[Q^{\frac 12}D^2U^{\tau}(t,x)Q^{\frac 12}].
\end{align}

\begin{lm}\label{reg-kol}
For every $\alpha, \theta, \gamma \in [0,1)$, $\theta+\gamma<1$, there exist $\tau_0>0$, $C( T,Q,\tau_0,\alpha)>0$ and 
$C(T, Q,\tau_0,\theta, \gamma )>0$ such that for $\tau \in (0,\tau_0]$, $x\in E, y,z \in \HH$ and $t\in (0,T]$,
\begin{align}\label{kol-d1}
|DU^{\tau}(t,x)\cdot y|&\le \frac {C( T,Q,\tau_0,\alpha)(1+|x|_{E}^{2})}{t^{\alpha}}\|(-A)^{-\alpha}y\|, \\\label{kol-d2}
|D^2U^{\tau}(t,x)\cdot (y,z)|&\le \frac {C(T,Q, \tau_0,\theta, \gamma )(1+|x|_{E}^{9})}{t^{\theta+\gamma}}\|(-A)^{-\theta}y\|\|(-A)^{-\gamma }z\|.
\end{align} 
\end{lm}

Based on the above estimates, now we give the weak error estimate of \eqref{full-spde}.

\textbf{Proof of Theorem \ref{weak0}}
The main idea of deducing the sharp weak convergence rate
lies on the decomposition of $\E\Big[\phi(X(T))-\phi(X^N_{K})\Big]$ into $\E\Big[\phi(X(T))-\phi(X^{\tau}(T))\Big]$ and $\E\Big[\phi(X^{\tau}(T))-\phi(X^N_{K})\Big]$. The first term is estimated by  Lemma \ref{spl} and possesses the strong convergence order 1 with respect to the parameter $\tau$. 
The strong error estimate of the second term is obtained based on Theorem \ref{weak}. Combining these estimations together, we complete the proof of Theorem \ref{weak0}.
\qed

\begin{lm}\label{spl}
Let  Assumptions \ref{as-lap}-\ref{as-dri} hold and  $p\ge 1$.
Then the solution $X^{\tau}$ of Eq. \eqref{rspde} is strongly convergent to the solution
$X$ of Eq. \eqref{spde} and satisfies 
\begin{align*}
&\E\Big[\sup_{t\in [0,T]}\|X^{\tau}(t)\|_E^p\Big]\le C(T,Q,X_0,p),\\
&\Big\|\sup_{t\in [0,T]} \|X^{\tau}(t)-X(t)\|\Big\|_{L^p(\Omega)} \le C(T,Q,X_0,p)\tau,
\end{align*} 
where $C(T,Q,X_0,p)$ is a positive number.
\end{lm}

In the following,  we are devoted to estimating the term $\E\Big[\phi(X^{\tau}(T))-\phi(X^N_{K})\Big]$.
For convenience, we introduce the continuous interpolation of the implicit full discretization.
Similar to \cite{BG18}, we define for $k\in \N^+$, $t\in [t_k,t_{k+1}]$, $\widehat X^N(t_k)=X^N_k$, 
\begin{align*}
d\widehat X^N(t)=(AS_{\delta t}X^N_k 
+S_{\delta t}P^NF(X_{k+1}^N))dt+S_{\delta t}P^NdW(t).
\end{align*}

\begin{tm}\label{weak}
Let  Assumptions \ref{as-lap}-\ref{as-dri} hold with $\beta \in (0,1]$, $\gamma\in (0,\beta)$,
$T>0$ and $\delta t_0\in 
(0,1\land \frac 1{(2\lambda_F-2\lambda_1)\lor 0})$.
Then for any  $\phi\in \mathcal C_b^2(\HH)$, 
there exists $\tau>0$  and $C(T,X_0,Q,\phi)>0$ such that for any $\delta t\in (0, \delta t_0]$, $K\delta t=T$, $K\in\N^+$ and $N\in \N^+$,

\begin{align*}
\Big|\E\Big[\phi(X^{\tau}(T))-\phi(X^N_{K})\Big]\Big|
\le C(T,X_0,Q,\phi)\Big(\delta t^{\gamma}+\lambda_N^{-\gamma}\Big).
\end{align*}
\end{tm}

\textbf{Proof}
We decompose the error  $\E\Big[\phi(X^{\tau}(T))-\phi(X^N_{K})\Big]$  
as
\begin{align*}
&\E\Big[U^{\tau}(T,X_0)\Big]-\E\Big[U^{\tau}(0,X^N_K)\Big]\\
=&\Big(\E\Big[U^{\tau}(T,X_0)\Big]-\E\Big[U^{\tau}(T,X_0^N)\Big]\Big)+\Big(\E\Big[U^{\tau}(T,X^N_0)\Big]-\E\Big[U^{\tau}(0,X^N_K)\Big]\Big).
\end{align*}

The first term is controlled,  by the regularity estimate  of $U^{\tau}$ in Lemma \ref{reg-kol}, as
\begin{align*}
&\Big|\E\Big[U^{\tau}(T,X_0)\Big]-\E\Big[U^{\tau}(T,X_0^N)\Big]\Big|\\
&\le
\int_0^1\Big|\E\Big[DU^{\tau}(T,\theta X_0+(1-\theta)X_0^N)\cdot(I-P^N)X_0\Big]\Big|d\theta\\
&\le C(1+\|X_0\|_E^2+\|X_0^N\|_E^2)\min(T^{-\alpha} \lambda_N^{-\alpha}\|X_0\|,
\lambda_N^{-\frac \beta 2}\|X_0\|_{\HH^{\beta}}).
\end{align*} 
By using the It\^o formula for Skorohod integrals (see e.g. \cite[Chapter 3]{Nua06}), the Kolmogorov equation \eqref{kol} and Malliavin integration by parts, $\E\Big[U^{\tau}(T,X^N_0)\Big]-\E\Big[U^{\tau}(0,X^N_K)\Big]$ is split as
\begin{align*}
&\E\Big[U^{\tau}(T,X^N_0)\Big]-\E\Big[U^{\tau}(0,X^N_K)\Big]\\
&= \sum_{k=0}^{K-1}
\E\Big[U^{\tau}(T-t_k,X^N_k)\Big]-\E\Big[U^{\tau}(T-t_{k+1},X^N_{k+1})\Big]\\
&= 
\E\Big[U^{\tau}(T,X^N_0)\Big]-\E\Big[U^{\tau}(T-\delta t,X^N_1)\Big]\\
&\quad-\sum_{k=1}^{K-1} \E\Big[\int_{t_k}^{t_{k+1}}\sum_{l\in \N^+}D^2U^{\tau}(T-t, \widehat X^N(t))\cdot(\mathcal D_t \widehat  X^N(t)Q^{\frac 12}e_l,S_{\delta t}Q^{\frac 12}e_l)\Big]dt\\
&\quad+
\sum_{k=1}^{K-1}
\Big(\int_{t_k}^{t_{k+1}}\E\Big[\<DU^{\tau}(T-t,\widehat X^N(t)),A\widehat X^N(t)-AS_{\delta t}X_k^N\>\Big]dt\\
&\quad+\int_{t_k}^{t_{k+1}}\E\Big[\<DU^{\tau}(T-t,\widehat X^N(t)), \Psi_\tau(\widehat X^N(t))-S_{\delta t}P^NF(X_{k+1}^N)\>\Big]dt\\
&\quad+\frac 12\int_{t_k}^{t_{k+1}}\sum_{j\in \N^+}\E\Big[D^2U^{\tau}(T-t,\widehat X^N(t))\cdot (Q^{\frac 12}e_j, Q^{\frac 12}e_j)-
(S_{\delta t}P^NQ^{\frac 12}e_j, S_{\delta t}P^NQ^{\frac 12}e_j)\Big)\Big]dt\Big)\\
&=:\E\Big[U^{\tau}(T,X^N_0)\Big]-\E\Big[U^{\tau}(T-\delta t,X^N_1)\Big]+\sum_{k=1}^{K-1} I^k_1+I^k_2+I^k_3+I^k_4.
\end{align*}
The Markov property of $X^{N}_k$, the regularity estimate  \eqref{kol-d1} of $U^{\tau}$ in Lemma \ref{reg-kol} 
and the a priori estimates of $X^N$ and $X^N_k$ in Lemma \ref{pri-xn} and  in Corollary \ref{pri-xnk} yield that for $0<\alpha<1$,
\begin{align*}
&\Big|\E\Big[U^{\tau}(T,X^N_0)\Big]-\E\Big[U^{\tau}(T-\delta t,X^N_1)\Big]\Big|\\
&=\Big|\E\Big[U^{\tau}(T-\delta t,X^\tau(\delta t,X^N_0))-U^{\tau}(T-\delta t,X^N_1)\Big]\Big|\\
&\le C (1+\E[\|X^\tau(\delta t,X^N_0))\|_E^2]+\E[\|X^N_1\|_E^2])(1+(T-\delta t)^{-\alpha })\delta t^{\alpha}\\
&\le C(Q,X_0)(1+(T-\delta t)^{-\alpha })\delta t^{\alpha}.
\end{align*}
For the term $I_1^k,$
the regularity estimate of $U^{\tau}$ and the a priori estimate of $\widehat X^N$ yield that 
\begin{align*}
\Big|\sum_{k=1}^{K-1}I_1^k\Big|
&\le C\sum_{k=1}^{K-1}\int_{t_k}^{t_{k+1}}(T-t)^{\frac {\beta -1}2}
\E\Big[(1+\|\widehat X^N(t)\|^{9}_E)
\|\mathcal D_t\widehat X^N(t)\|_{\LL_2^0}\|S_{\delta t}\|_{\LL}\|(-A)^{\frac {\beta-1}2}\|_{\LL_2^0}\Big]dt
\\
&\le C(T,Q,X_0)\delta t\sum_{k=1}^{K-1}\int_{t_k}^{t_{k+1}}(T-t)^{\frac {\beta-1}2}((t_{k+1}-[t]_{\delta t})^{\frac {\beta-1}2}+1)dt
\le C(T,Q,X_0)\delta t^{\beta}.
\end{align*}
where  we  use Proposition \ref{mal-dif} and the fact that for $t_k\le  t\le s \le t_{k+1}$,
\begin{align*}
\mathcal D_s\widehat X^N(t)
&=S_{\delta t}\mathcal D_sX^N(t_k)
+(t-t_k)P^NS_{\delta t}DF(\widehat X^N(t_{k+1}))\mathcal D_s \widehat X^N_{k+1}
+\mathcal D_s \int_{t_k}^{t}S_{\delta t}dW(s)\\
&=(t-t_k)P^NS_{\delta t}DF(\widehat X^N(t_{k+1}))\mathcal D_s \widehat X^N_{k+1}.
\end{align*}

Next, we  estimate $I_2^k, I_3^k$ and $I_4^k$, $k\ge 1$ separately.  
The definition of $\widehat X$ leads to
\begin{align*}
I^k_2
&= \int_{t_k}^{t_{k+1}}\E\Big[\<DU^{\tau}(T-t,\widehat X^N(t)),A(X_k^N-S_{\delta t}X_k^N)\>\Big]dt\\
 &\quad+\int_{t_k}^{t_{k+1}}\E\Big[\<DU^{\tau}(T-t,\widehat X^N(t)),(t-t_k)(-A)^2S_{\delta t}X_k^N\>\Big]dt\\
 &\quad+ \int_{t_k}^{t_{k+1}}\E\Big[\<DU^{\tau}(T-t,\widehat X^N(t)),
(t-t_k)AS_{\delta t}P^NF(X_{k+1}^N)\>\Big]dt\\
&\quad+ \int_{t_k}^{t_{k+1}}\E\Big[\<DU^{\tau}(T-t,\widehat X^N(t)),A
\int_{t_k}^tS_{\delta t}P^NdW(s)\>\Big]dt
\\
&=:I^k_{21}+I^k_{22}+I^k_{23}+I^k_{24}.
\end{align*}

From 
the property $I-S_{\delta t}=-A\delta t(I-A\delta t)^{-1}$, the mild form  of 
 $X_k^N$, the 
a priori estimate of $\widehat X$, and
the regularity estimate of $U^\tau$ and the smoothing effect of $S_{\delta t}$, it follows that for $k\ge 1$, any small $\epsilon_1>0$,
\begin{align*}
|I^k_{21}|&\le  \Big|\int_{t_k}^{t_{k+1}}\E\Big[\<DU^{\tau}(T-t,\widehat X^N(t)),(-A)^2\delta tS_{\delta t}^{k+1}X_0^N\>\Big]dt\Big|
\\
&\quad+ \Big|\int_{t_k}^{t_{k+1}}\sum_{j=0}^{k-1}\E\Big[\<DU^{\tau}(T-t,\widehat X^N(t)),(-A)^2\delta t^2S_{\delta t}^{k+1-j}P^NF(X^N_{j+1})\>\Big]dt\Big|\\
&\quad+\Big| \int_{t_k}^{t_{k+1}}\E\Big[\<DU^{\tau}(T-t,\widehat X^N(t)),(-A)^2\delta t\sum_{j=0}^{k-1}S_{\delta t}^{k+1-j}P^N\delta W_j\>\Big]dt\Big|\\
&\le C\delta t\int_{t_k}^{t_{k+1}} (T-t)^{-\alpha}
\E \Big[(1+\|\widehat X^N(t)\|_E^2)\|(-A)^{1-\epsilon_1}S_{\delta t}^{k}\|\|(-A)^{1-\alpha+\epsilon_1}S_{\delta t}\|
\|X_0^N\|\Big]dt
\\
&\quad+C\delta t^2\int_{t_k}^{t_{k+1}}(T-t)^{-\alpha}\sum_{j=0}^{k-1}\E\Big[(1+\|\widehat X^N(t)\|_E^2)\|(-A)^{1-\epsilon_1}S_{\delta t}^{k-j}\|\|(-A)^{1-\alpha+\epsilon_1}S_{\delta t}\|
\|F(X^N_{j+1})\|\Big]dt\\
&\quad+
\Big| \int_{t_k}^{t_{k+1}}\E\Big[\<DU^{\tau}(T-t,\widehat X^N(t)),(-A)^2\delta t\sum_{l=0}^{k-1}S_{\delta t}^{k+1-j}P^N\delta W_j\>\Big]dt\Big|\\
&\le C(T,X_0,Q)\delta t^{\alpha-\epsilon_1}\int_{t_k}^{t_{k+1}} (T-t)^{-\alpha} (t_k^{-1+\epsilon_1}+1)dt\\
&\quad+\Big| \int_{t_k}^{t_{k+1}}\E\Big[\<DU^{\tau}(T-t,\widehat X^N(t)),(-A)^2\delta t\sum_{j=0}^{k-1}S_{\delta t}^{k+1-j}P^N\delta W_j\>\Big]dt\Big|.
\end{align*}
By using the Malliavin calculus integration by parts and the Malliavin differentiability of $\widehat X^N$, we have 
\begin{align*}
&\Big| \int_{t_k}^{t_{k+1}}\E\Big[\<DU^{\tau}(T-t,\widehat X^N(t)),(-A)^2\delta t\sum_{j=0}^{k-1}S_{\delta t}^{k+1-j}P^N\delta W_j\>\Big]dt\Big|\\
&= \delta t\int_{t_k}^{t_{k+1}} \sum_{j=0}^{k-1}\int_{t_j}^{t_{j+1}}\sum_{l\in\N^+} \E\Big[\Big| D^2U^{\tau}(T-t,\widehat X^N(t))\cdot (\mathcal D_s^{Q^{\frac12 }e_l}\widehat X^N(t),(-A)^2S_{\delta t}^{k+1-j}P^NQ^{\frac 12}e_l)\Big|\Big]dsdt\\
&\le C\delta t\int_{t_k}^{t_{k+1}}
\sum_{j=0}^{k-1}\int_{t_j}^{t_{j+1}}
\sum_{l\in\N^+}  \E\Big[\Big| \< (-A)^{\frac {1+\beta}2-\epsilon_1}D^2U^{\tau}(T-t,\widehat X^N(t)) (-A)^{\frac {1-\beta}2}\\
&\qquad (-A)^{\frac {\beta-1}2} \mathcal D_s^{Q^{\frac12 }e_l}\widehat X^N(t),(-A)^{2-\frac {\beta+1}2+\epsilon_1}S_{\delta t}^{k+1-j}P^NQ^{\frac 12}e_l)\>\Big|\Big]dsdt
\\
&\le C\delta t\int_{t_k}^{t_{k+1}}(T-t)^{-1+\epsilon_1}
\sum_{j=0}^{k-1}\int_{t_j}^{t_{j+1}}
\E\Big[(1+\|\widehat X^N(t)\|_E^9)\|(-A)^{\frac {\beta-1}2}\mathcal D_s \widehat X^N(t)\|_{\LL_2^0}\\
&\qquad\|(-A)^{1-\epsilon_1} S_{\delta t}^{k-j}\|
\|(-A)^{1-\beta+2\epsilon_1}S_{\delta t}\|\|(-A)^{\frac {\beta-1}2}\|_{\LL_2^0}\Big]dsdt\\
&\le C(T,X_0,Q)\delta t^{\beta-2\epsilon_1}\int_{t_k}^{t_{k+1}}(T-t)^{-1+\epsilon_1}\int_0^{t_{k}}(t_{k}-[s]_{\delta t})^{-1+\epsilon_1}dsdt.
\end{align*} 
The above analysis leads to \begin{align*}
|I^k_{21}|&\le C(T,X_0,Q)\delta t^{\alpha-\epsilon_1}\int_{t_k}^{t_{k+1}} (T-t)^{-\alpha} (t_k^{-1+\epsilon_1}+1)dt\\
&\quad+C(T,X_0,Q)\delta t^{\beta-2\epsilon_1}\int_{t_k}^{t_{k+1}}(T-t)^{-1+\epsilon_1}\int_0^{t_{k}}(t_{k}-[s]_{\delta t})^{-1+\epsilon_1}dsdt,
\end{align*}
for $k\ge 1$.
Since the estimation for $I^k_{22}$ for $k\ge 1$ is similar,
we omit the procedures.
For $I^k_{23}$, by the regularity estimate of  $DU^\tau$, we have 
\begin{align*}
|I^k_{23}|&\le C \delta t\int_{t_k}^{t_{k+1}}(T-t)^{-1+\epsilon_1}\E\Big[(1+\|\widehat X^N(t)\|_E^2)\|(-A)^{\epsilon_1}S_{\delta t}\|\|F(X_{k+1}^N)\|\Big]dt\\
&\le C(T,Q,X_0)\delta t^{1-\epsilon_1}\int_{t_k}^{t_{k+1}}(T-t)^{-1+\epsilon_1}dt.
\end{align*}
Again using  Malliavin calculus integration by parts
yields that
\begin{align*}
|I^k_{24}|&=
\Big|\int_{t_k}^{t_{k+1}}\E\Big[\<DU^{\tau}(T-t,\widehat X^N(t)),A
\int_{t_k}^tP^NS_{\delta t}dW(s)\>\Big]dt\Big|\\
&=\Big|\int_{t_k}^{t_{k+1}}\int_{t_k}^t\E\Big[\<D^2U^{\tau}(T-t,\widehat X^N(t))\mathcal D_s\widehat X^N(t),P^NA
S_{\delta t}\>_{\LL_2^0}\Big]dsdt\Big|\\
&\le \Big|\int_{t_k}^{t_{k+1}}\int_{t_k}^t\E\Big[\<(-A)^{\frac {1+\beta}2-\epsilon_1}D^2U^{\tau}(T-t,\widehat X^N(t))(-A)^{\frac {1-\beta}2} (-A)^{\frac {\beta -1}2}\mathcal D_s\widehat X^N(t), \\
&\quad (-A)^{\epsilon_1}P^N
S_{\delta t}(-A)^{\frac {1-\beta}2}\>_{\LL_2^0}\Big]dsdt\Big|\\
&\le C(T,Q,X_0)\delta t^{\beta-\epsilon_1} 
\int_{t_k}^{t_{k+1}}(T-t)^{-1+\epsilon_1}dt.
\end{align*}
Thus we have 
\begin{align*}
|I_2^k|\le C(T,Q,X_0)\delta t^{\beta-2\epsilon_1}\Big(\int_{t_k}^{t_{k+1}}(T-t)^{-1+\epsilon_1}(1+\int_0^{t_{k}}(t_{k}-[s]_{\delta t})^{-1+\epsilon_1}ds)dt\Big).
\end{align*}

Now, we are in a position to control $I_3^k$.
It follows from the continuity of $\Psi_\tau$ in \cite[Lemma 4.2]{CH18} and  the regularity of $DU^{\tau}$ that
\begin{align*}
|I_{3}^k|&\le \Big|\int_{t_k}^{t_{k+1}}\E\Big[\<DU^{\tau}(T-t,\widehat X^N(t)), \Psi_\tau(\widehat X^N(t))-
F(\widehat X^N(t))\>\Big]dt\Big|\\
&\quad+\Big|\int_{t_k}^{t_{k+1}}\E\Big[\<DU^{\tau}(T-t,\widehat X^N(t)), (I-P^N)F(\widehat X^N(t))\>\Big]dt\Big|\\
&\quad+\Big|\int_{t_k}^{t_{k+1}}\E\Big[\<DU^{\tau}(T-t,\widehat X^N(t)), (I-S_{\delta t})P^NF(X_{k+1}^N)\>\Big]dt\Big|\\
&\quad+ \Big|\int_{t_k}^{t_{k+1}}\E\Big[\<DU^{\tau}(T-t,\widehat X^N(t)),P^N(F(\widehat X^N(t))-F(X_{k+1}^N)\>\Big]dt\Big|\\
&\le C\tau\int_{t_k}^{t_{k+1}}\E\Big[1+\|\widehat X^N(t)\|_E^7\Big]dt+
C(\lambda_N^{-\alpha}+\delta t^{\alpha})
\int_{t_k}^{t_{k+1}}(T-t)^{-\alpha}\\
&\qquad\E\Big[(1+\|\widehat X^N(t)\|_E^2)(1+\|\widehat X^N(t)\|_{L^6}^3+\|\widehat X^N_{k+1}\|_{L^6}^3)\Big]dt\\
&\quad+ \Big|\int_{t_k}^{t_{k+1}}\E\Big[\<DU^{\tau}(T-t,\widehat X^N(t)),P^N(F(\widehat X^N(t))-F(X_{k+1}^N)\>\Big]dt\Big|.
\end{align*} 
Thus it  suffices to estimate the last term in the above inequality.
 From the Taylor expansion of $F$, the regularity estimate of  $DU^{\tau}$ and the a priori estimate of $\widehat X^N$, it follows  that 
\begin{align*}
&\int_{t_k}^{t_{k+1}}\E\Big[\<DU^{\tau}(T-t,\widehat X^N(t)),P^N(F(\widehat X^N(t))-F(X_{k+1}^N)\>\Big]dt\\
&\le \int_{t_k}^{t_{k+1}}(t-t_{k+1})\E\Big[\<DU^{\tau}(T-t,\widehat X^N(t)),P^N(DF(\widehat X^N(t))\cdot (AS_{\delta t}X_{k}^N))\>\Big]dt\\
&\quad+\int_{t_k}^{t_{k+1}}(t-t_{k+1})\E\Big[\<DU^{\tau}(T-t,\widehat X^N(t)),P^N(DF(\widehat X^N(t))\cdot (S_{\delta t}P^NF(X_{k+1}^N))\>\Big]dt\\
&\quad+\int_{t_k}^{t_{k+1}}\E\Big[\<DU^{\tau}(T-t,\widehat X^N(t)),P^N(DF(\widehat X^N(t))\cdot (\int_t^{t_{k+1}}P^NS_{\delta t}dW(s)))\>\Big]dt\\
&\quad+\int_{t_k}^{t_{k+1}}\E\Big[\<DU^{\tau}(T-t,\widehat X^N(t)),P^N(\int_0^1(1-\theta)D^2F(\theta \widehat X^N(t)+(1-\theta) X^N_{k+1})) \\
&\qquad\cdot (\widehat X^N(t)-X^N_{k+1},\widehat X^N(t)-X^N_{k+1})d\theta\>)\Big]dt
=:I^k_{31}+I^k_{32}+I^k_{33}+I^k_{34}.
\end{align*}
The mild form of $X_{k+1}^N$ and Malliavin calculus integration by parts  yield that
\begin{align*}
|I^k_{31}|
&=\Big|\int_{t_k}^{t_{k+1}}(t-t_{k+1})\Big(\E\Big[\<DU^{\tau}(T-t,\widehat X^N(t)),P^N(DF(\widehat X^N(t))\cdot (AS_{\delta t}^{k+1}X_{0}^N))\>\Big]\\
&\quad+\delta t\E\Big[\<DU^{\tau}(T-t,\widehat X^N(t)),P^N(DF(\widehat X^N(t))\cdot (\sum_{j=0}^{k-1}AS_{\delta t}^{k+1-j}P^NF(X_{j+1}^N))\>\Big]
\\
&\quad+\E\Big[\<DU^{\tau}(T-t,\widehat X^N(t)),P^N(DF(\widehat X^N(t))\cdot (\sum_{j=0}^{k-1}\int_{t_j}^{t_{j+1}}AP^NS_{\delta t}^{k+1-j}dW(s))\>\Big]\Big)dt\Big|\\
&\le 
C\delta t^2\sup_{t\in[0,T]}\E\Big[(1+\|\widehat X^N(t)\|^4_E)\|(-A)^{1-\epsilon_1}S_{\delta t}^k\|\|(-A)^{\epsilon_1}S_{\delta t}\|\|X_0\|\Big]]\\
&\quad+
C\delta t^2
\int_{t_k}^{t_{k+1}}\sum_{j=0}^{k-1}\E\Big[(1+\|\widehat X^N(t)\|^4_E)
\|(-A)^{1-\epsilon_1}S_{\delta t}^{k-j}\|(-A)^{\epsilon_1} S_{\delta t}\|\|F(X_{j+1}^N)\|\Big]dt\\
&\quad+
 \Big|\int_{t_k}^{t_{k+1}}(t-t_{k+1})
\sum_{j=0}^{k-1}\int_{t_j}^{t_{j+1}}
\sum_{l\in \N^+}\E\Big[D^2U^\tau(T-t,\widehat X^N(t))\cdot P^N((DF(\widehat X(t))\\
&\qquad \cdot(P^NAS_{\delta t}^{k+1-j}Q^{\frac 12}e_l)), \mathcal D_s^{Q^{\frac 12}e_l}(\widehat X^N(t))\Big]dsdt \Big|\\
&\quad+ \Big|\int_{t_k}^{t_{k+1}}(t-t_{k+1})
\sum_{j=0}^{k-1}\int_{t_j}^{t_{j+1}}
\sum_{l\in \N^+}\E\Big[\<DU^\tau(T-t,\widehat X^N(t)), P^N(D^2F(\widehat X^N(t)\\
&\qquad \cdot(P^NAS_{\delta t}^{k+1-j}Q^{\frac 12}e_l,\mathcal D_s^{Q^{\frac 12}e_l}(\widehat X^N(t)))\>\Big]dsdt\Big|.
\end{align*}

By the a priori estimate of $\widehat X^N$ and the Sobolev embedding theorem $E \hookrightarrow \HH^{\frac d2+\epsilon_1} $, $\epsilon_1>0$,  we have 
\begin{align*}
|I^k_{31}|&\le
C(T,	Q,X_0)\delta t^{1-\epsilon_1}(t_k^{-1+\epsilon_1}\delta t
+\delta t\sum_{j=0}^{k-1}t_{k-j}^{-1+\epsilon_1})\\
&\quad+
 C\delta t\int_{t_k}^{t_{k+1}}
\sum_{j=0}^{k-1}\int_{t_j}^{t_{j+1}}
\sum_{l\in \N^+}\E\Big[(1+\|\widehat X(t)\|^{11}_E)
\|AS_{\delta t}^{k+1-j}Q^{\frac 12}e_l\| \|\mathcal D_s^{Q^{\frac 12}e_l}\widehat X(t)\|\Big]dsdt \\
&\quad+ C\delta t\int_{t_k}^{t_{k+1}}
\sum_{j=0}^{k-1}\int_{t_j}^{t_{j+1}}
\sum_{l\in \N^+}\E\Big[\Big\|(-A)^\eta DU^\tau(T-t,\widehat X(t))\Big\| \Big\|(-A)^{-\eta}D^2F(\widehat X(t)\\
&\qquad \cdot(P^NAS_{\delta t}^{k+1-j}Q^{\frac 12}e_l,\mathcal D_s^{Q^{\frac 12}e_l}(\widehat X(t)) \Big\|\Big]dsdt,
\end{align*}
where $\eta >\frac d4+\frac {\epsilon_1}2$.
And by using the smoothing effect of $S_{\delta t}$, the Malliavin regularity and the a priori estimate of $\widehat X(t)$, we have
\begin{align*}
|I^k_{31}|
&\le
C(T,	Q,X_0)\delta t^{1-\epsilon_1}(t_k^{-1+\epsilon_1}\delta t
+\delta t\sum_{j=0}^{k-1}t_{k-j}^{-1+\epsilon_1})\\
&\quad+
 C(T,Q,X_0)\delta t\int_{t_k}^{t_{k+1}}
\sum_{j=0}^{k-1}\int_{t_j}^{t_{j+1}}
\E\Big[\|(-A)^{\frac {3-\beta}2}S_{\delta t}^{k+1-j}(-A)^{\frac {
\beta -1}2}\|_{\LL_2^0} \|\mathcal D_s\widehat X(t)\|_{\LL_2^0}\Big]dsdt \\
&\quad+ C\delta t\int_{t_k}^{t_{k+1}}
\sum_{j=0}^{k-1}\int_{t_j}^{t_{j+1}}
(T-t)^{-\eta}\E\Big[(1+\|\widehat X(t)\|_E^3)
\|(-A)^{\frac {3-\beta}2}S_{\delta t}^{k+1-j}(-A)^{\frac {\beta-1}2}\|_{\LL_2^0}
\| \mathcal D_s\widehat X(t)\|_{\LL_2^0}\Big]dsdt\\
&\le C(T,	Q,X_0)\delta t^{1-\epsilon_1}t_k^{-1+\epsilon_1}\delta t
+C(T,Q,X_0)\delta t^{\beta-\epsilon_1}
\int_{t_k}^{t_{k+1}}
\int_{0}^{t_{k}}(t_{k}-[s]_{\delta t})^{-\frac{\beta +1}2 +\epsilon_1}dsdt
\\
&\quad+C(T,Q,X_0)\delta t^{\beta-\epsilon_1}\int_{t_k}^{t_{k+1}}(T-t)^{-\eta}\int_{0}^{t_k}(t_k-[s]_{\delta t})^{-1+\epsilon_1}dt.
\end{align*}
Similarly, we get 
\begin{align*}
|I^k_{32}|
&\le C\delta t\int_{t_{k}}^{t_{k+1}}
\E\Big[(1+\|\widehat X^N(t)\|_E^4)(1+\|\widehat X^N(t)\|_{L^6}^3)\Big]dt
\le C\delta t^2
\end{align*}
and
\begin{align*}
|I^k_{33}|&\le
\Big|\int_{t_k}^{t_{k+1}}\int_t^{t_{k+1}}\sum_{l\in \N^+}\E\Big[D^2U^{\tau}(T-t,\widehat X^N(t))\cdot\Big(P^N(DF(\widehat X^N(t))\cdot (P^NS_{\delta t}Q^{\frac 12}e_l), \mathcal D_s^{Q^ {\frac 12e_l}}\widehat X^N(t)\Big)\Big]dsdt
\Big|\\
&+
\Big|\int_{t_k}^{t_{k+1}}\int_t^{t_{k+1}}\sum_{l\in \N^+}\E\Big[\<DU^{\tau}(T-t,\widehat X^N(t)),P^N(D^2F(\widehat X^N(t))\cdot (\mathcal D_s^{Q^{\frac 12}e_l} \widehat X^N(t),P^NS_{\delta t}Q^{\frac 12}e_l )\>\Big]dsdt\Big|
\\
&\le 
C(T,Q,X_0)\delta t^{\beta-\epsilon_1}.
\end{align*}
Combining with the continuity of $\widehat X^N$, we have that for $t\in [t_k,t_{k+1}]$,
\begin{align*}
\|\widehat X^N(t)-\widehat X^N_{k+1}\|_{L^p(\Omega;\HH)}
&\le 
(t_{k+1}-t)\|(-A)^{1-\frac \beta 2}S_{\delta t}\|_{\LL(\HH)}
\|X(t_k)\|_{L^p(\Omega;\HH^{ \beta })}
+\|X(t_k)-X^N_k\|_{L^p(\Omega;\HH^{\beta })}\\
&\quad+C\|F(X_{k+1}^N)\|_{L^p(\Omega;\HH)}(t_{k+1}-t)
+\|\int_{t}^{t_{k+1}}S_{\delta t}dW(s)\|_{L^p(\Omega;\HH)}\\
&\le C(T,X_0,Q)(t_{k+1}-t)^{\frac \beta 2},
\end{align*}
which implies that for $\eta>\frac d4+\epsilon_1$, 
\begin{align*}
|I^k_{34}|&\le 
C\Big|\int_{t_k}^{t_{k+1}}\E\Big[\<DU^{\tau}(T-t,\widehat X^N(t)),P^N(\int_0^1D^2F(\theta \widehat X^N(t)+(1-\theta) X^N_{k+1})) \\
&\qquad\cdot (\widehat X^N(t)-\widehat X^N_{k+1},\widehat X^N(t)-\widehat X^N_{k+1})(1-\theta)d\theta\>)\Big]dt\Big|\\
&\le C\int_{t_k}^{t_{k+1}}(T-t)^{-\eta}\E\Big[(1+\|\widehat X(t)\|^{3}_E)\|\widehat X^N(t)-\widehat X^N_{k+1} \|^2\Big]dt\\
&\le C(T,Q,X_0)\delta t^{\beta}\int_{t_k}^{t_{k+1}}(T-t)^{-\eta}dt.
\end{align*}
From the estimations of $II_{31}^k$-$II_{34}^k$, it is concluded that 
\begin{align*}
|I_3^k|&\le 
C(T,Q,X_0)\delta t^{\beta-\epsilon_1}
\Big(1+t_k^{-1+\epsilon_1}\delta t+
\int_{t_k}^{t_{k+1}}(T-t)^{-\eta}\int_{0}^{t_k}(1+(t_k-[s]_{\delta t})^{-1+\epsilon_1})dt\Big).
\end{align*}

For $I_4^k$, by applying the regularity estimate  of $D^2U^{\tau}$, we obtain
\begin{align*}
|I_4^k|
&\le \Big|\int_{t_k}^{t_{k+1}}\sum_{j\in \N^+}\E\Big[D^2U^{\tau}(T-t,\widehat X^N(t))\cdot\Big((I-P^N)Q^{\frac 12}e_j, (I+P^N)Q^{\frac 12}e_j\Big)\Big]dt\Big|\\
&\quad+ \Big|\int_{t_k}^{t_{k+1}}\sum_{j\in \N^+}\E\Big[D^2U^{\tau}(T-t,\widehat X^N(t))\cdot\Big(P^N(I-S_{\delta t})Q^{\frac 12}e_j, P^N(I+S_{\delta t})Q^{\frac 12}e_j\Big)\Big]dt\Big|\\
&\le C \int_{t_k}^{t_{k+1}}(T-t)^{-1+\epsilon_1}\E\Big[(1+\|\widehat X^N(t)\|^9_E)\|(-A)^{\frac {\beta-1}2}\|_{\LL_2^0}^2\|(-A)^{-\frac {1+\beta}2+\epsilon_1}(I-P^N)(-A)^{\frac {1-\beta}2}\|\Big]dt\\
&\quad+C \int_{t_k}^{t_{k+1}}(T-t)^{-1+\epsilon_1}\E\Big[(1+\|\widehat X^N(t)\|^9_E)\|(-A)^{\frac {\beta-1}2}\|_{\LL_2^0}^2\|(-A)^{-\frac {1+\beta}2+\epsilon_1}(I-S_{\delta t})(-A)^{\frac {1-\beta}2}\|\Big]dt\\
&\le C(T,Q,X_0)(\delta t^{\beta-\epsilon_1}+\lambda_N^{-\beta+\epsilon_1})\int_{t_k}^{t_{k+1}}(T-t)^{-1+\epsilon_1}dt.
\end{align*}
Combining all the estimations of $I^k_1$-$I^k_4$ and summing up over $k$, taking $\tau = \mathscr O(\delta t^{\beta})$ or $\mathscr O(\lambda_N^{-\beta})$, we finish the proof.
\qed

\section{Time-independent weak convergence analysis and approximation of invariant measures}
\label{sec-5}
In this section, we consider whether the proposed method can be used to approximate the invariant measure
of Eq. \eqref{spde}. 
Different from analyzing weak error in Section \ref{sec-4}, we need to give the time-independent regularity estimates of the 
Kolmogorov equation, which are 
more involved.

\subsection{V-uniform ergodicity for the semi-discretization}

To  ensure the existence of  a unique ergodic invariant measure for Eq. \eqref{spde} and to study the invariant measure numerically, the following assumptions are introduced, that is, the dissipative condition in Assumption \ref{erg-as}
and the non-degenerate condition in  Assumption  \ref{erg-as1}.

\begin{as}\label{erg-as}
Let $\lambda_F<\lambda_1$ and $\|(-A)^{\frac {\beta -1}2}\|_{\LL^2_0}<\infty$, $\beta \le 1$.
\end{as}

The above Assumption \ref{erg-as} immediately implies the following  result on  exponential convergence to equilibrium for Eq. \eqref{spde}. 

\begin{prop}\label{exp-erg0}
Under Assumptions \ref{as-lap}-\ref{erg-as}, there exist $c>0, C>0$ such that 
for any $\phi\in \mathcal C_b^1(\HH)$, $t\ge 0$ and  $x_1,x_2\in \HH$,
\begin{align*}
|\E[\phi(X(t,x_1))-\phi(X(t,x_2))]|
&\le C|\phi|_{1}e^{-ct}(1+\|x_1\|^2+\|x_2\|^2).
\end{align*}
\end{prop}

\begin{rk}\label{ind-str}
Based on the proof of Lemma \ref{str} and Corollary \ref{pri-xnk}, together with the strict dissipative condition $\lambda_F< \lambda_1$ in  Assumption \ref{erg-as},  the full discretization 
is strongly convergent and satisfies
\begin{align*}
\sup_{k\in \N^+}\Big\|X_k^N-X(t_k)\Big\|_{L^p(\Omega;\HH)}
&\le C(X_0,Q)(\delta t^{\frac \beta 2}+\lambda_N^{-\frac \beta 2}).
\end{align*}
\end{rk}

In some situations, it may occurs that $\lambda_1\le \lambda_F$, which leads that  Assumption \ref{erg-as} does not hold. In this case, we introduce the following non-degenerate condition.

\begin{as}\label{erg-as1}
Let 
the covariance operator $Q$ be invertible
and commute with $A$, $\|Q^{-\frac 12}(-A)^{-\frac 12}\|<\infty$ and $\|(-A)^{\frac {\beta -1}2}\|_{\LL^2_0}<\infty$, $\beta \le 1$.
\end{as}

Under Assumption \ref{erg-as1}, the existence of the unique invariant measure $\mu$ for Eq. \eqref{spde}, as well as the invariant measure $\mu^N$ for the spatial Galerkin method, can be obtained according to Doob theorem for general $\lambda_F\in \R.$
Besides the ergodicity of the invariant measure, we
also need the following exponential convergence result in Proposition \ref{exp-erg1}. Its proof  lies heavily on the strong Feller property and $V$-uniform ergodicity of the Markov 
semigroup $P_t$ generated by the solution of Eq. \eqref{spde} and Eq. \eqref{semi-spde} (see, e.g., \cite{DMT95,GM06}). In fact, we first follow the proof of \cite{GM06} to show
the a priori estimate of a Lyapunov functional $V$ and to obtain the existence of the invariant measure. Then we prove that the Markov semigroup of the solution is strong Feller and irreducible, which implies the uniqueness and ergodicity of the invariant measure. By using again a priori estimate of $V$, one can obtain the $V$-uniform ergodicity. 
In particular, we choose $\phi\in B_b(\HH)$ to get 
the exponential ergodicity of the invariant measure, which  immediately  implies Proposition \ref{exp-erg1}.

\begin{prop}\label{exp-erg1}
Under Assumptions \ref{as-lap}-\ref{as-dri}  and Assumptions \ref{erg-as1},
there exist $c>0$, $C>0$ such that for any $\phi\in B_b(\HH)$ and for $t\ge0$,  any $x_1,x_2\in \HH$ and $y_1^N,y_2^N \in P^N(\HH)$, we have
\begin{align}\label{exp-erg-exa}
|\E[\phi(X(t,x_1)]-\E[\phi(X(t,x_2))]|&\le 
C\|\phi\|_{0}e^{-ct}(1+\|x_1\|^2+\|x_2\|^2),\\\label{exp-erg-sem}
|\E[\phi(X^N(t,y_1^N)]-\E[\phi(X^N(t,y_2^N))]|&\le 
C\|\phi\|_{0}e^{-ct}(1+\|y_1^N\|^2+\|y_2^N\|^2).
\end{align}
\end{prop}

\textbf{Proof}
For the exponential convergence to equilibrium \eqref{exp-erg-exa} for the original equation, we refer to \cite{GM06}. 
We  focus on the semi-discretization and define
$$P_t^N\phi(x)=\E\phi(X^N(t,x)),\quad\phi \in B_b(P^N(\HH)),\quad t\ge 0.$$
For the sake of simplicity, we omit the index $N$ of $P_t^N$
for convenience.
The Markov property and Feller property of $P_t$ can be obtained by the similar arguments in \cite[Chapter 4]{Dap04}.
The left proof will be divided into three steps. 

Step 1: $P_t$ is strong Feller.
To get the strong Feller property, $P_t(B_b(P^N(\HH))\subset C_b(P^N(\HH))$  for $t>0$, it suffices to show that for any $\phi \in C_b(P^N(\HH))$ and $t>0$,  there exists $C(t)>0$ such that $\sup\limits_{x\in P^N(\HH)}\|DP_t\phi(x)\| \le C(t)\|\phi\|_0$. Indeed, the strong Feller property follows from 
$|P_t\phi(x)-P_t\phi(y)|\le C(t)\|\phi\|_0\|x-y\|, x,y\in P^N(\HH)$ and the density of $C_b(P^N(\HH))$ in $B_b(P^N(\HH))$.

Now, we are in a position to deduce the regularity estimate of $P_t$, i.e., $\|DP_t\phi\|_0\le C(t)\|\phi\|_0$.
Recall that $\eta^h(t,x)=D\E[X^N(t,x)]\cdot h$ satisfies 
\begin{align*}
\frac 12\|\eta^h(t,x)\|^2+\int_0^t\|\nabla \eta^h(t,x)\|^2ds \le \frac 12\|h\|^2+\int_0^t\lambda_F\|\eta^h(s,x)\|^2ds.
\end{align*} 
This, combined with the equivalence of the norms in  $\HH^1$ and $H\cap H^1_0$, implies that 
\begin{align*}
\int_0^t \|(-A)^{\frac 12} \eta^h(t,x)\|^2ds\le C(t)\|h\|^2.
\end{align*}
The Bismut--Elworthy--Li formula
\begin{align*}
\<DP_t\phi(x),h\>=\frac 1t\E[\phi(X^N(t,x))\int_0^t\<
Q^{-\frac 12}\eta^h(s,x),P^NdW(s)\>],
\end{align*}
together with the H\"older inequality, leads to 
\begin{align*}
\|DP_t\phi(x)\|_0^2&\le  \frac 1{t^2}\|\phi\|_0^2
\E\Big[\int_0^t\sup_{\|h\|\le 1}\|P^N(Q^{-\frac 12}\eta^h(s,x))\|^2ds\Big]
\\
&\le C(t)\|\phi\|_0^2\|Q^{-\frac 12}(-A)^{-\frac 12}\|^2,
\end{align*}
for $t>0$,	
which implies the strong Feller property of $P_t$.

Step 2: $P_t$ is irreducible. A basic tool for proving the irreducibility property is using the approximate controllability of 
the following system 
\begin{align}\label{ske}
d\widetilde X^N(t)&= A\widetilde X^N(t)dt
+P^NF(\widetilde X^N(t))dt+P^N(Q^{\frac 12}u(t))dt, \; t>0,\\\nonumber 
\widetilde X^N&=x,
\end{align}
where $x\in P^N(\HH)$ and $u\in L^{2}([0,T];P^N(\HH))$.
Denoting by $\widetilde X^N(t,x,u)$  the mild solution of the above system, it follows that 
\begin{align*}
\widetilde X^N(t)
=e^{tA}x+\int_0^te^{(t-s)A}P^N(F(\widetilde X^N(s)))ds
+\int_0^te^{(t-s)A}P^N(Q^{\frac12}u(s))ds.
\end{align*}

Thus it needs to show that for any fixed time $T>0$, for any $\epsilon >0$, $x_0,x_1\in P^N(\HH)$, there exists $u \in L^2([0,T];P^N(\HH))$ such that 
$\|\widetilde X^N(T,x_0,u)-x_1\|\le \epsilon$.
Now, we denote $\alpha_{x_0,x_1}(t)=\frac {T-t}Tx_0+\frac tTx_1$
and $\beta_{x_0,x_1}(t)=\frac {d}{dt}\alpha_{x_0,x_1}(t)-A\alpha_{x_0,x_1}(t)-P^NF(\alpha_{x_0,x_1}(t)), t\in [0,T]$.
 
Since $x_0,x_1\in P^N(\HH)$ $\subset D(A)$ and $Q$ is invertible, we  choose $u\in C([0,T];P^N(\HH))$ such that 
$\|\beta_{x_0,x_1}(t)-Q^{\frac 12}u(t)\|\le C\epsilon,$ $t\in [0,T]$.  Denote $z(t)=X^N(t,x_0,u)-\alpha_{x_0,x_1}(t)$, $t\in [0,T]$.
By using the monotonicity of $F$ and  the dissipativity     
of $A$, we have 
\begin{align*}
\frac 12\|z(t)\|^2
&\le \int_0^t-\|\nabla z(t)\|^2ds
+\int_0^t \<F(X^N(s,x_0,u))-F(\alpha_{x_0,x_1}(s)),
z(s)\>ds\\
&\quad+\int_0^t\<Q^{\frac 12}u(s)-\beta_{x_0,x_1}(s), z(s)\>ds\\
&\le \int_0^t(\frac12+ \lambda_F-\lambda_1)\|z(s)\|^2ds
+\frac T2 C^2\epsilon^2.
\end{align*}
Then the Gronwall inequality implies that 
$
\|z(T)\| \le C\sqrt{T}e^{(\frac 12+\lambda_F-\lambda_1)T}\epsilon.
$
Choosing a proper $C$ completes the proof of the approximate controllability. 
By applying the approximate controllability of the skeleton equation \eqref{ske}, we deduce that 
for $x_0,x_1\in P^N(\HH)$ and $T>0$,
$\PP(\|X^N(T,x_0)-x_1\|<\epsilon)>0$.
Indeed, the approximate controllability leads to 
the existence of a control $u\in L^2([0,T];P^N(\HH))$ such that $\|\widetilde X^N(T,x_0,u)-x_1\|\le \frac \epsilon 2$.
Then we have 
\begin{align*}
\PP(\|X^N(T,x_0)-x_1\|\ge \epsilon)\le 
\PP\left(\|X^N(T,x_0)-\widetilde X^N(T,x_0,u)\|\ge \frac \epsilon2\right).
\end{align*}
Similar arguments in the proof of the priori estimate of $X^N$ lead to 
\begin{align*}
\|X^N(t,x_0)-\widetilde X^N(t,x_0,u)\|
\le \|Y^N(t,x_0)-\widetilde Y^N(t,x_0,u)\|
+\|Z^N(t)-\widetilde Z^N(t,u)\|,
\end{align*}
where $\widetilde X^N=\widetilde Y^N+\widetilde Z^N$,  $\widetilde Y^N$ and $\widetilde Z^N$
satisfy 
\begin{align*}
\frac d{dt} \widetilde Z^N&=A\widetilde Z^N
+Q^{\frac 12}u(t),\; \widetilde Z^N(0)=0,
\\
\frac d{dt}\widetilde Y^N&=A\widetilde Y^N
+P^NF(\widetilde Y^N+\widetilde Z^N),\;
\widetilde Y^N= x_0,
\end{align*}
By the monotonicity of $F$ and the dissipativity     
of $A$, we have
\begin{align*}
&\frac 12\|Y^N(t,x_0)-\widetilde Y^N(t,x_0,u)\|^2\\
\le& \int_0^t -\lambda_1\|Y^N(s,x_0)-\widetilde Y^N(s,x_0,u)\|^2ds\\
&+\int_0^t\<F(Y^N+ Z^N)-F(\widetilde Y^N+\widetilde Z^N),Y^N- \widetilde Y^N\>ds\\
\le& \int_0^t (-\lambda_1+\lambda_F)\|Y^N(s,x_0)-\widetilde Y^N(s,x_0,u)\|^2ds\\
&+\int_0^t\<F(\widetilde Y^N+ Z^N)-F(\widetilde Y^N+\widetilde Z^N),Y^N- \widetilde Y^N\>ds\\
\le& \int_0^t C\|Y^N(s,x_0)-\widetilde Y^N(s,x_0,u)\|^2ds\\
&+\int_0^tC(1+\|\widetilde Y^N\|_E^4+\|\widetilde Z^N\|_E^4+\| Z^N\|_E^4)\|Z^N-\widetilde Z^N\|^2ds.
\end{align*} 
Then the Gronwall inequality leads to 
\begin{align*}
\|Y^N(T,x_0)-\widetilde Y^N(T,x_0,u)\|
&\le Ce^{CT}\sqrt{T}\|Z^N-\widetilde Z^N\|_{C([0,T];P^N(\HH))}(1+\|\widetilde Y^N\|^2_{C([0,T];E)}\\
&\qquad+\|\widetilde Z^N\|_{C([0,T];E)}^2+\|Z^N\|_{C([0,T];E)}^2).
\end{align*}
The Sobolev embedding theorem $\HH^2\hookrightarrow E$, the inverse inequality $\|x\|_{\HH^2}\le \lambda_N\|x\|, x\in P^N(\HH)$,
and the uniform boundedness of $\widetilde Y^N$, $\widetilde Z^N$ and $Z^N$ imply that 
\begin{align*}
&\|Y^N(T,x_0)-\widetilde Y^N(T,x_0,u)\|
\le C(\lambda_N)e^{CT}\sqrt{T}\|Z^N-\widetilde Z^N\|_{C([0,T];P^N(\HH))}\\
&\qquad\left(1+\|\widetilde Y^N\|^2_{C([0,T];P^N(\HH))}+\|Z^N-\widetilde
 Z^N\|_{C([0,T];P^N(\HH))}^2+\|\widetilde Z^N\|_{C([0,T];P^N(\HH))}^2\right).
\end{align*}
It is concluded that 
\begin{align*}
&\PP\left(\|X^N(T,x_0)-\widetilde X^N(T,x_0,u)\|\ge \frac \epsilon2\right)\\
\le& \PP\left( \|Y^N(T,x_0)-\widetilde Y^N(T,x_0,u)\|
+\|Z^N(T)-\widetilde Z^N(T,u)\| \ge \frac \epsilon2\right)\\
\le& \PP\Big( C(\lambda_N)e^{CT}\sqrt{T}\|Z^N-\widetilde Z^N\|_{C([0,T];P^N(\HH))}(1+\|\widetilde Y^N\|^2_{C([0,T];P^N(\HH))}\\
&\qquad+\|Z^N-\widetilde Z^N\|_{C([0,T];P^N(\HH))}^2+\|\widetilde Z^N\|_{C([0,T];P^N(\HH))}^2)\ge \frac \epsilon 2\Big).
\end{align*}
Since $Z^N$ is full in $C([0,T];P^N(\HH))$  and $C(\lambda_N)$ is polynomially dependent on $\lambda_N$, we have that  there exists $R=R(\lambda_N,\epsilon,T)$ such that 
\begin{align*}
&\PP\left(\|X^N(T,x_0)-\widetilde X^N(T,x_0,u)\|\ge \frac \epsilon 2\right)\\
&\le \PP\left(\|Z^N-\widetilde Z^N\|_{C([0,T];P^N(\HH))} \ge R(\lambda_N,\epsilon,T)\right)<1,
\end{align*} 
which completes the proof of the irreducibility property.

Step 3: Existence of the $V$-uniformly ergodic invariant measure.
Similar arguments in the proof of Lemma \ref{pri-y}
imply  
the uniform estimate of $X^N$ in $\HH^{\beta}, \beta\in (0,1]$. 
The existence of the invariant measure $\mu^N$
of Eq. \eqref{semi-spde} is ensured by the uniform estimate of $X^N$ in $\HH^{\beta}$ and the Sobolev compact embedding theorem. To show the exponential ergodicity of the invariant measure, by \cite[Theorem 12.1]{GM06}, it suffices to show that the $p$-th moment of $X^N(t,x)$ is ultimately bounded, i.e., $\E[\|X^N(t,x)\|^p]\le k|x|^pe^{-\omega t}+c$, $t\ge 0$, $x\in P^N(\HH)$, for some positive constants $k,\omega, c$, and $p$.  

For convenience, we only prove the case that $p=2$. Due to the fact that $X^N(t)=Y^N(t)+Z^N(t)$, we estimate the $\HH$-norm of $Z^N$ and $Y^N$, respectively.
The mild form of $Z^N$ yields that 
\begin{align*}
\E[\|Z^N(t)\|^2]&\le \E[\|Z^N(t)\|_E^2]\le C(Q).
\end{align*}
It follows from the variational approach, Poincare, Young 
and H\"older inequalities that
\begin{align*}
\E[\|Y^N(t)\|^2]&= \E[\|Y^N(0)\|^2]
-2\int_0^t\|\nabla Y^N(s)\|^2ds
+2\int_0^t\<F(Y^N(s)+Z^N(s)),Y^N(s)\>ds
\\
&\le 
\E[\|Y^N(0)\|^2]
-2\lambda_1\int_0^t\|Y^N(s)\|^2ds
+C(\epsilon)\int_0^t(1+\|Z^N\|_{L^4}^4)ds
\\
&\quad
+\int_0^t(-(2a_3-\epsilon)\|Y^N(s)\|_{L^4}^4+C\|Y^N(s)\|^2)ds\\
&\le
\E[\|Y^N(0)\|^2]
-2\lambda_1\int_0^t\|Y^N(s)\|^2ds+
C(\epsilon)\int_0^t (1+
\|Z^N(s)\|_{E}^4)ds.
\end{align*} 
Then the Gronwall inequality and  the a priori estimate of $\|Z^N(s)\|_{E}$ imply that 
\begin{align*}
\E[\|Y^N(t)\|^2]
&\le e^{-2\lambda_1t}\|x\|^2
+\int_0^t e^{-2\lambda_1(t-s)}C(Q)ds
\le e^{-2\lambda_1t}\|x\|^2+C(Q,\lambda_1).
\end{align*}
From \cite[Theorem 12.1]{GM06}, it follows that 
$\{P_t^N\}_{t\ge 0}$, is $V$-uniformly ergodic with $V=1+\|x\|^2$, i.e.,
\begin{align*}
\sup_{\|\phi\|_V\le 1}\Big|P_t\phi(x)-\int_{P^N(\HH)}\phi d\mu^N\Big|
\le CV(x)e^{-\alpha t},\; x\in P^N(\HH),
\end{align*}
where $C=C(\alpha,\lambda_1)$, $\alpha=\alpha(\alpha,\lambda_1)$, and $\phi \in B_b^V(P^N(\HH))$, i.e, $\phi$ is 
Borel-measurable and  $\|\phi\|_{V}:=\sup_{x\in P^N(\HH)}\frac {\|\phi(x)\|}{V(x)}<\infty$. 
Now taking any $\phi\in B_b(P^N(\HH))$, we have 
$\|\phi\|_V\le \|\phi\|_0$. Combining  with $V$-uniformly ergodicity of $P_t$, we deduce that 
\begin{align*}
\Big|P_t\phi(x)-\int_{P^N(\HH)}\phi d\mu^N\Big|
\le C(\|\phi\|_0)V(x)e^{-\alpha t},\; x\in P^N(\HH).
\end{align*}
we obtain the exponential ergodicity of the unique invariant measure.
By the fact that $B_b(\HH)\subset B_b(P^N(\HH))$,  
taking two different initial data $y_1^N$, $y_2^N$, combining with the exponential ergodicity of $\mu^N$, we complete the proof. 
\qed

\begin{rk}
Under the same conditions of Proposition \ref{exp-erg1}, one can obtain the 
uniformly exponentially
ergodicity ($V=1$) of $P_t^N$, $t\ge 0$ (see e.g. \cite{GM06}),
\begin{align*}
\sup_{|\phi|_0\le 1}\Big|P_t^N\phi(x)-\int_{P^N(\HH)}\phi d\mu^N\Big|\le Ce^{-ct},
\end{align*}
which can be used to improve the bound of regularity estimates in Lemma \ref{reg-um1}.
The condition $\|Q^{-\frac 12}(-A)^{-\frac12}\|<\infty$ in Assumption \ref{erg-as1} is necessary for the strong Feller property of $P_t$. However, from the proof of the strong Feller property of $P_t^N$, 
it follows that  
the estimate \eqref{exp-erg-sem} holds even for the case
$\|Q^{-\frac 12}(-A)^{-\frac12}\|=\infty $.
\end{rk}

\subsection{Time-independent regularity estimate of Kolmogorov equation}
In order to obtain the time-independent optimal weak error estimate, we need more refined regularity estimates of the Kolmogorov equation.
However, the regularizing approach by the splitting strategy may be not used directly to get these time-independent regularity estimates.
To overcome this difficulty, we investigate time-independent regularity estimates of the Kolmogorov equation by means of  a finite dimensional approximation. 
Recall the Kolmogorov equation of the Galerkin approximation
\begin{align}\label{kol-um}
\frac {\partial U^{M}(t,x)} {\partial t}
&=\<Ax+P^MF(x), D U^{M}(t,x)\>+\frac 12tr[P^MQ^{\frac 12}(P^MQ^{\frac 12})^*D^2U^{M}(t,x)],
\end{align}
where $M\in \N^+$.
The proofs of the time-independent regularity estimates under Assumption \ref{erg-as} or Assumption \ref{erg-as1}
are totally different. Under the strong dissipative condition $\lambda_1>\lambda_F$, we show the exponential decay of the regularity estimate by 
variational arguments. Once it occurs that 
$\lambda_1\le \lambda_F$, the $V$-uniform ergodicity is used to achieve this type estimate.

\begin{lm}\label{reg-um}
Let Assumptions \ref{as-lap}-\ref{erg-as} hold and 
$\phi\in \mathcal C_b^2(\HH)$.
For  $\alpha, \theta, \gamma \in [0,1)$, $\theta+\gamma<1$ and $\eta \in (\frac d4,1)$, there exist $c>0$, $C(Q,\alpha,\phi,\eta)$ and 
$C(Q,\theta, \gamma,\phi,\eta)$ such that for $x, h,k \in P^M(\HH)$, $M\in \N^+$ and $t>0$,
\begin{align}\label{kol-d3}
|DU^{M}(t,x)\cdot h|&\le C(Q,\alpha,\phi,\eta)(1+\sup_{s\in[0,t]}\E[\|X^M(s,x)\|_{E}^{2}])(1+t^{-\alpha})e^{-c t}\|(-A)^{-\alpha}h\|, 
\\\nonumber
|D^2U^{M}(t,x)\cdot(h,k)|&\le C(Q,\theta, \gamma,\phi,\eta)(1+\sup_{s\in[0,t]}\E[\|X^M(s,x)\|_{E}^{7}])(1+t^{-\eta}+{t^{-\theta-\gamma}})
\\\label{kol-d4}
&\qquad e^{-c t}\|(-A)^{-\theta}h\|\|(-A)^{-\gamma }k\|.
\end{align} 
\end{lm}
\textbf{Proof}
Similar to the proof of \cite[Proposition 4.1]{CH18},  we have 
\begin{align*}
DU^M(t,x)\cdot h&=\E[D\phi(X^M(t,x))\cdot \eta^h(t,x)],\\
DU^M(t,x)\cdot(h,k)
&=\E[D\phi(X^M(t,x))\cdot \zeta^{h,k}(t,x)]\\
&\quad+\E[D^2\phi(X^M(t,x))\cdot (\eta^h(t,x),\eta^k(t,x))]
\end{align*}
for $h,k,x \in P^M(\HH)$, $t\ge 0$,
where $\eta^h$ and $\zeta^{h,k}$ satisfy  
\begin{align*}
\frac {\partial \eta^h(t,x)} {\partial t}
&= A \eta^h(t,x)+P^M(DF(X^M(t,x))\eta^h(t,x)),\;\eta^h(0,x)=h,\\
\frac {\partial \zeta^{h,k}(t,x)}{\partial t}
&=A \zeta^{h,k}(t,x)+P^M(DF(X^M(t,x))\zeta^{h,k}(t,x))\\
&\quad+P^M(D^2F(X^M(t,x))\eta^h(t,x)\eta^k(t,x)),\;\zeta^{h,k}(0,x)=0.
\end{align*}
For convenience, the parameter $M$ is omitted in the notations of $\eta^h$ and $\zeta^{h,k}$.
Consider the following auxiliary equation 
\begin{align*}
\frac {\partial V(t,s)h}{\partial t}
=(A+P^MDF(X^M(t,x)))V(t,s)h,\; V(s,s)h=h.
\end{align*}
The straightforward argument leads to 
$
\|V(t,s)h\|^2\le e^{-2(\lambda_1-\lambda_F)(t-s)}\|h\|^2.
$

Denote $\widetilde \eta^h(t,x):=\eta^h(t,x)-e^{tA}h$. 
It follows from 
the smoothing effect of $e^{tA}$ and the estimate of $V(t,s)$, $0\le s\le t$, that for some $c\in (0,\lambda_1)$,
\begin{align*}
|\E[D\phi(X^M(t,x))\cdot e^{tA}h]|
&\le |\phi|_{1}C_{\alpha}t^{-\alpha}e^{-ct}\|(-A)^{-\alpha}h\|,
\end{align*}
and
\begin{align*}
\|\widetilde \eta^h(t,x)\|
&=\|\int_0^tV(t,s)(P^M(DF(X^M(t,x))e^{sA}h)ds\|\\
&\le
\int_0^t e^{-(\lambda_1-\lambda_F)(t-s)}\|DF(X^M(s,x))\|_E\|e^{sA}h\|ds\\
&\le C|\phi|_{1}\int_0^t (1+\|X^M(s,x)\|_E^2)e^{-(\lambda_1-\lambda_F)(t-s)}s^{-\alpha}e^{-cs}ds\|(-A)^{-\alpha}h\|.
\end{align*}
Taking $c>\lambda_1-\lambda_F$, we obtain
\begin{align*}
&|\E[D\phi(X^M(t,x))\cdot \widetilde \eta^h(t,x)]|\\
&\le 
C|\phi|_{1}\int_0^t (1+\E[\|X^M(s,x)\|_E^2])e^{-(\lambda_1-\lambda_F)(t-s)}s^{-\alpha}e^{-cs}ds\|(-A)^{-\alpha}h\|
\\
&\le Ce^{-(\lambda_1-\lambda_F)t}\int_0^t (1+\E[\|X^M(s,x)\|_E^2]) s^{-\alpha}e^{-(c-\lambda_1+\lambda_F)s}ds\|(-A)^{-\alpha}h\|\\
&\le  C(1+\sup_{s\in[0,t]}\E[\|X^M(s,x)\|_E^2])e^{-(\lambda_1-\lambda_F)t}\|(-A)^{-\alpha}h\|.
\end{align*}
The above two estimates imply that 
\begin{align*}
&|\E[D\phi(X^M(t,x))\cdot \eta^h(t,x)]|\\
&\le C(1+\sup_{s\in[0,t]}\E[\|X^M(s,x)\|_E^2])(1+t^{-\alpha})e^{-ct}\|(-A)^{-\alpha}h\|.
\end{align*}
Similarly, we have 
\begin{align*}
&\E[D^2\phi(X^M(t,x))\cdot (\eta^h(t,x),\eta^k(t,x))]\\
&\le C|\phi|_{2}(1+\sup_{s\in[0,t]}\E[\|X^M(s,x)\|_E^4])
(1+t^{-(\beta+\gamma)})e^{-ct}\|(-A)^{-\beta}h\|\|(-A)^{-\gamma}h\|.
\end{align*}
Denote  $\widetilde V(t,s)h=V(t,s)h-e^{-(t-s)A}h$.
Similar arguments imply that for $t>s\ge 0$, $0\le \alpha<1$,
\begin{align}\label{reg-vts}
\|V(t,s)h\|
&\le Ce^{-(\lambda_1-\lambda_F)(t-s)} \Big((t-s)^{-\alpha}+\int_s^t  (r-s)^{-\alpha}e^{(-c+\lambda_1-\lambda_F)(r-s)}\\\nonumber
&\qquad\|DF(X^M(r,x))\|_E dr\Big)\|(-A)^{-\alpha}h\|.
\end{align}
Based on the representation of $\zeta^{h,k}$ and 
\eqref{reg-vts}, we obtain 
\begin{align*}
\|\zeta^{h,k}(t,x)\|
&=\Big\|\int_0^tV(t,s)P^M\Big(D^2F(X^M(s,x))\eta^y(s,x)\eta^z(s,x)\Big)ds\Big\|\\
&\le  C\int_0^te^{-c(t-s)}(t-s)^{-\eta}\big(1+\int_s^te^{-c_1(r-s)}(r-s)^{-\eta}\|DF(X^M(r,x))\|_Edr\big)\\
&\qquad  \|(-A)^{-\eta}P^M((D^2F(X^M(s,x))\eta^y(s,x)\eta^z(s,x)))\|ds,
\end{align*}
for $\eta>\frac d4$, $c,c_1>0$.
Thus we have 
\begin{align*}
&\E[D\phi(X^M(t,x))\cdot \zeta^{h,k}(t,x)]\\
&\le Ce^{-ct}(1+\sup_{r\in [0,t]}\E[\|(X^{M}(r,x))\|^7_E])\int_0^te^{-cs}(t-s)^{-\eta}\\
&\qquad(1+s^{-\theta-\gamma})\|(-A)^{-\theta}h\|\|(-A)^{-\gamma}k\|ds\\
&\le Ce^{-c_1t}(1+t^{-\eta})(1+\sup_{r\in [0,t]} \E[\|(X^{M}(r,x))\|^7_E])\|(-A)^{-\theta}h\|\|(-A)^{-\gamma}k\|,
\end{align*}
which completes the proof.
\qed

\begin{lm}\label{reg-um1}
Let Assumptions \ref{as-lap}-\ref{as-dri} and  Assumption \ref{erg-as1} hold, $\phi\in \mathcal C_b^2(\HH)$.
For $\alpha, \theta, \gamma \in [0,1)$, $\theta+\gamma<1$ and $\eta\in(\frac d4,1)$, there exist $c>0$, $C(Q,\alpha,\phi,\eta)$ and 
$C(Q,\theta, \gamma,\phi,\eta )$ such that for  $x, y,z \in P^M(\HH)$ and $t\in (0,T]$,
\begin{align}\label{kol-d5}
|DU^{M}(t,x)\cdot y|&\le {C(Q,\alpha)}(1+\sup_{s\in[0,t]}\E[\|X^M(s,x)\|_{E}^{4}])\\\nonumber
&\qquad(1+t^{-\alpha})e^{-c t}\|(-A)^{-\alpha}y\|, \\\label{kol-d6}
|D^2U^{M}(t,x)\cdot(y,z)|&\le {C(Q,\theta, \gamma )|}(1+\sup_{s\in[0,t]}\E[\|X^M(s,x)\|_{E}^{14}])\\\nonumber
&\qquad  
(1+t^{-\eta}+{t^{-\theta-\gamma}})
e^{-c t}\|(-A)^{-\theta}y\|\|(-A)^{-\gamma }z\|.
\end{align} 
\end{lm}

\textbf{Proof}
By the similar arguments in Lemma \ref{reg-um},
we obtain the regularity estimate for $0<t\le T$.
For convenience, we take $T=1$ and get for $0<t\le 1$ 
\begin{align}\label{reg-est-t=1}
\|\eta^{h}(t,x)\|&\le  {C(1+\sup_{s\in[0,t]}\|X^M(s,x)\|_{E}^{2})}{(1+t^{-\alpha})}\|(-A)^{-\alpha}h\|, \\\nonumber 
\|\zeta^{h,k}(t,x)\|&\le  {C(1+\sup_{s\in[0,t]}\|X^M(s,x)\|_{E}^{7})}(1+{t^{-\eta}}+t^{-\theta-\gamma})\|(-A)^{-\theta}h\|\|(-A)^{-\gamma}h\|.
\end{align}
To get the time-independent regularity estimate, we need another a priori estimate of $X^M$ which does not depend on $x$.
According to the fact $X^M=Y^M+Z^M$ and the evolution of $\|Y^M\|^2$, we have
\begin{align*}
\frac {\partial}{\partial t} \|Y^M(t)\|^2
&= -2\|\nabla Y^M(t)\|^2
+2\<F(Y^M(t)+Z^M(t)),Y^N(t)\>\\
&\le 
-c\|Y^M(t)\|^4+C(1+\|Z^M(t)\|_{L^4}^4),
\end{align*}
 for $0<c<2$, $C>0$.
By applying \cite[Lemma 1.2.6]{Cer01}, we have
\begin{align}\label{y-ind}
\sup_{x\in P^M(\HH)}\|Y^M(t)\|^p
&\le C(p,t)(t\land 1)^{-\frac p 2},
\end{align}
 for $p\ge 1$,
where $C(p,t)$ has  finite moments of any order.
Combining with the equivalence of  norms in finite dimensional space, we have 
\begin{align*}
\|\eta^{h}(t,x)\|&\le C(M,t,|\phi|_{1})\|h\|, \\
\|\zeta^{h,k}(t,x)\|&\le C(M,t,|\phi|_{2})\|h\|\|k\|,
\end{align*}
for $t>0$.
Indeed, by the chain rule, we have 
\begin{align*}
\|\eta^h(t,x)\|^2+\int_0^t\|\nabla \eta^h(s,x)\|^2ds
\le \|h\|^2+\int_0^t \<DF(X^M(s))\eta^h(s,x),\eta^h(s,x)\>ds.
\end{align*}
Therefore
\begin{align*}
\|\eta^h(t,x)\|^2+\int_0^t\|\nabla \eta^h(s,x)\|^2ds
\le \|h\|^2+C\int_0^t\|\eta^h(s,x)\|^2ds.
\end{align*}
The Gronwall inequality leads that 
$\|\eta^{h}(t,x)\|\le C(M,t,|\phi|_{1})\|h\|$ 
for $ t>0$.
The same arguments, together with \eqref{y-ind}, the Sobolev embedding theorem and the inverse inequality, yield that 
\begin{align*}
\|\zeta^{h,k}(t,x)\|
&\le \|\int_0^t V(t-s) D^2F(X^M(s))\eta^h(s,x)\eta^k(s,x)ds\|\\
&\le e^{Ct}\int_0^t \|X^M(s)\|_{L^6}\|\eta^h(s,x)\|_{L^6}
\|\eta^k(s,x)\|_{L^6}ds\\
&\le  e^{Ct}\int_0^t \|X^M(s)\|_{\HH^1}\|\eta^h(s,x)\|_{\HH^1}
\|\eta^k(s,x)\|_{\HH^1}ds\\
&\le e^{Ct}\lambda_M^{\frac 32}\int_0^t \|X^M(s)\|\|\eta^h(s,x)\|
\|\eta^k(s,x)\|ds\\
&\le e^{Ct}\lambda_M^{\frac 32}\sqrt{t}\|h\|\|k\|.
\end{align*}
Thus, we get 
\begin{align}\label{x-ind}
|DU^{M}(t,x)\cdot h|&\le e^{Ct}\|h\|,\\\nonumber
|D^2U^{M}(t,x)\cdot (h,k)|&\le C(p,t)e^{Ct}\lambda_M^{\frac 32}\sqrt{t}\|h\|\|k\|,
\end{align}
which implies that for any $\phi \in \mathcal C_b^2(\HH)$, $U^M(t)=P_t\phi \in \mathcal C_b^2(\HH)$, $t>0$.

The Bismut--Elworthy--Li formula states that if $\Phi: P^M(\HH) \to \R$ belongs to $\mathcal C_b^2(P^M(\HH) )$ and $|\Phi(x)|\le M(\Phi)(1+|x|^q), q\ge 1$, then we can calculate the first and second order derivatives of $U_{\Phi}^M(t,x):=\E [\Phi(X^M(t,x))]$
with respect to $x$.
Indeed, we have
\begin{align*}
DU_{\Phi}^M(t,x)\cdot h
=\frac 1t\E\Big[\int_0^t\<Q^{-\frac 12}\eta^{h}(s,x), d\widetilde W(s)\>\Phi(X^M(t,x))\Big],
\end{align*}
  for any $x\in P^M(\HH)$, $h\in P^M(\HH)$.
The Markov property of $P_t$ implies that 
\begin{align*}
DU_{\Phi}^M(t,x)\cdot h
=\frac 2t\E\Big[\int_0^{\frac t2}\<Q^{-\frac 12}\eta^{h}(s,x), d\widetilde W(s)\>U_{\Phi}^N(\frac t2, X^M(\frac t2, x))\Big],
\end{align*}
where $\widetilde W$ is the cylindrical Wiener process. 
By applying again the Bismut--Elworthy--Li formula, 
we get  a formula of the second derivative
\begin{align*}
D^2U_{\Phi}^M(t,x)\cdot (h,k)
&=\frac 2t \E\Big[\int_0^{\frac t2}\<Q^{-\frac 12}\zeta ^{h,k}(t,x), d\widetilde W(s)\>U_{\Phi}^M(\frac t2, X^M(\frac t2, x))\Big]\\
&\quad+\frac 2t \E\Big[\int_0^{\frac t2}\<Q^{-\frac 12}\eta ^{h}(t,x), d\widetilde W(s)\>DU_{\Phi}^M(\frac t2, X^M(\frac t2, x))\cdot \eta^{k}(\frac t2,x)\Big],
\end{align*}
for $x,h,k\in P^M(\HH)$.
By using a priori estimates of $\eta^h$ and $\zeta^{h,k}$, we obtain 
\begin{align*}
|DU_{\Phi}^M(t,x)\cdot h|
&\le 
\frac 1t \sqrt{\E[ |\Phi(X^M(t,x))|^2]}
\sqrt{\E[\int_0^t \|Q^{-\frac 12}\eta^h(s,x)\|^2ds]}\\
&\le \frac 1t C(t)M(\Phi)(1+\sqrt{\E[ \|X^M(t,x)\|^{2q}]})\|h\|,
\end{align*}
and 
\begin{align*}
|D^2U_{\Phi}^M(t,x)\cdot (h,k)|
&\le 
\frac 2t \sqrt{\E[ |U^M_{\Phi}(\frac t2,X^M(\frac t2,x))|^2]}
\sqrt{\E[\int_0^t \|Q^{-\frac 12}\zeta^{h,k}(s,x)\|^2ds]}\\
&\quad+ \frac 1t \sqrt{\E[\int_0^t \|Q^{-\frac 12}\eta^h(s,x)\|^2ds]}\sqrt{\E[|DU_{\Phi}^M(\frac t2, X^M(\frac t2, x))\cdot \eta^{k}(\frac t2,x)|^2]}\\
&\le \frac 1{t}  C(t)M(\Phi)(1+\sqrt{\E[ \|X^M(\frac t2,x)\|^{2q}]})
\sqrt{\E[\int_0^t \|Q^{-\frac 12}\zeta^{h,k}(s,x)\|^2ds]}\\
&\quad+
 \frac 1{t^2} C(t,Q)M(\Phi)(1+\sqrt{\E[ \|X^M(\frac t2,x)\|^{2q}]})\|h\|\|k\|
\end{align*}
for $0<t\le 1$.
To estimate $\E[\int_0^t \|Q^{-\frac 12}\zeta^{h,k}(s,x)\|^2ds]$, we consider the $M$-independent  estimation of $\zeta^{h,k}$
and get
\begin{align*}
\|\zeta^{h,k}(t,x)\|
&\le \|\int_0^t V(t-s) D^2F(X^M(s))\eta^h(s,x)\eta^k(s,x)ds\|\\
&\le e^{Ct}\int_0^t (1+\|X^M(s)\|_{E})\|\eta^h(s,x)\|_{L^4}
\|\eta^k(s,x)\|_{L^4}ds\\
&\le e^{Ct}(1+\sup_{s\in [0,t]}\|X^M(s)\|_{E})
\int_0^t  \|(-A)^{\frac 12}\eta^h(s,x)\|\|(-A)^{\frac 12}\eta^h(s,x)\|ds\\
&\le C(t)(1+\sup_{s\in [0,t]}\|X^M(s)\|_{E})\|h\|\|k\|,
\end{align*}
which implies that 
\begin{align*}
\int_0^t \|Q^{-\frac 12}\zeta^{h,k}(s,x)\|^2ds
&\le 
C(t)\|\zeta^{h,k}(t,x)\|^2
+\int_0^t\<D^2F(X^M(s))\eta^h(s,x)\eta^k(s,x),\zeta^{h,k}(s,x)\>ds\\
&\le 
e^{Ct}C(t)(1+\sup_{s\in [0,t]}\|X^M(s)\|^2_{E})\|h\|^2\|k\|^2+C\sup_{s\in[0,t]}\|\zeta^{h,k}(s,x)\|\\
&\quad
(1+\sup_{s\in [0,t]}\|X^M(s)\|_{E})
\int_0^t \|(-A)^{\frac 12}\eta^h(s,x)\|\|(-A)^{\frac 12}\eta^k(s,x)\|ds\\
&\le e^{Ct}C(t)(1+\sup_{s\in [0,t]}\|X^M(s)\|^2_{E})\|h\|^2\|k\|^2.
\end{align*}
It is concluded that 
\begin{align*}
|D^2U_{\Phi}^M(t,x)\cdot (h,k)|
&\le 
(\frac 1t+\frac 1{t^2})C(t)M(\Phi)
(1+\sqrt{\E[ \|X^M(\frac t2,x)\|^{2q}]})\\
&\quad(1+\sup_{s\in [0,t]}\sqrt{\E[\|X^M(s,x)\|^2_{E}}])
\|h\|\|k\|.
\end{align*}

For any $t\ge 1$, we have $U^M(t,x)=\E[U^M(t-1,X^M(1,x))]$. The exponential convergence estimate \eqref{exp-erg-sem} yields that 
\begin{align*}
|U^M(t-1,x)-\int_{P^M(\HH)}\phi d\mu^M|
\le Ce^{-c(t-1)}(1+\|x\|^2).
\end{align*}
Inspired by \cite{Bre14}, we choose $\Phi(x)=U^M(t-1,x)-\int_{P^M(\HH)}\phi d\mu^M$, we have that $D\Phi$ and $D^2\Phi$ are uniformly bounded by \eqref{x-ind}. From the above estimate
\eqref{reg-est-t=1} in $0<t\le 1$ and the fact that $U^M(t,x)=\E[\Phi(X^M(1,x))]+\int_{P^M(\HH)}\phi d\mu^M$, it follows that 
\begin{align*}
|DU^M(t,x)\cdot h|&\le Ce^{-c(t-1)}(1+\sqrt{\E[ \|X^M(1,x)\|^{4}]})\|h\|\\
|D^2U^M(t,x)\cdot (h,k)|&\le Ce^{-c(t-1)}(1+\sqrt{\E[ \|X^M(\frac 12,x)\|^{4}]})\\
&\quad(1+\sup_{s\in [0,1]}\sqrt{\E[\|X^M(s)\|^2_{E}}])
\|h\|\|k\|,
\end{align*}
 for $t\ge 1$.
Thus we conclude that 
\begin{align*}
\|DU^M(t,x)\|&\le Ce^{-ct}(1+\sup_{s\in[0,1]}\sqrt{\E[ \|X^M(s,x)\|^{4}]}),\\
\|D^2U^M(t,x)\|_{\LL(\HH)}&\le 
Ce^{-ct}(1+\sup_{s\in[0,1]}\sqrt{\E[\|X^M(s,x)\|^{4}]})(1+\sup_{s\in[0,1]}\sqrt{\E[ \|X^M(s,x)\|_E^{2}]}).
\end{align*}
Combining with the Markov property of $X^M$, 
we have 
\begin{align*}
|DU^M(t,x)\cdot h|
&\le Ce^{-c(t-1)}\Big(1+\sqrt{\sup_{s\in[0,t]}\E[\|X^M(s,x)\|^4]}\Big)\sqrt{\E[\|\eta^h(1,x)\|^2]}\\
&\le Ce^{-c(t-1)}\Big(1+\sup_{s\in[0,t]}\E[\|X^M(s,x)\|_E^4]\Big)\|h\|_{\HH^{-\alpha}}\\
|D^2U^M(t,x)\cdot(h,k)|
&= |\E[D^2(U^M(t-1,X^M(1,x)))\cdot(h,k)]|\\
&\le  |\E[D^2U^M(t-1,X^M(1,x))\cdot(\eta^h(1,x),\eta^k(1,x))]|\\
&\quad+ |\E[DU^M(t-1,X^M(1,x))\cdot \zeta^{h,k}(1,x)]|\\
&\le Ce^{-c(t-1)}(1+\sup_{s\in[0,t]}\E[\|X^M(s,x)\|^{14}_E])\|h\|_{\HH^{-\theta}}\|k\|_{\HH^{-\gamma}}.
\end{align*}
The above estimate, together with \eqref{reg-est-t=1} completes the proof.
\qed

\subsection{Time-independent weak error estimate
and approximation of the invariant measure}

To get the time-independent weak convergence analysis, we introduce the another solution  $X^M$, $M\gg N$, of spectral Galerkin method.  The regularity estimates in Lemmas \ref{reg-um} and \ref{reg-um1} of $U^M$ are crucial. Before that, we first give a useful estimate to deal 
with the conditional expectation appeared in the regularity estimate of $U^M$ in Lemmas \ref{reg-um} and \ref{reg-um1}. 
 For convenience, we use the notation $\E_x$ as the conditional expectation at $x\in E$.

\begin{lm}\label{pri-con}
Under Assumptions \ref{as-lap}-\ref{as-dri},
for any $T>t\ge 0$, there exists a constant $C(Q,X_0,p)$ such that  for any $p\ge 2$,
\begin{align*}
\E\Big[\sup_{s\in [0,T-t]}\E_{\widehat X^N(t)}\Big[\|X^M(s, \widehat X^N(t))\|_E^p\Big]\Big]\le C(Q,X_0,p).
\end{align*}
\end{lm}
\textbf{Proof}
Without loss of generality, we assume that $t\in [t_k,t_{k+1})$, $0\le k\le K-1$.
By the procedures in proving Corollary \ref{pri-xnk}, 
we have, for $0\le s\le T-t$,
\begin{align*}
X^M(s,\widehat X^N(t))=Y^M(s,\widehat Y^N(t))+Z^M(s,\widehat Z^N(t)).
\end{align*}
Here $Y^M$ and $Z^M$ satisfy
\begin{align*}
dY^M&=AY^Mds+P^MF(Y^M+Z^M)ds,\\
dZ^M&=AZ^Mds+P^Md\widehat W(s),
\end{align*}
where $\widehat W(s)$, $s\ge 0$ is another Wiener process independent of $\{W(r)\}_{r\in [0,t_{k+1}]}$
and has the same distribution as $W(t+s)-W(t)$. Here $Y^M(0)=\widehat Y^N(t)$ and 
 $Z^M(0)=\widehat Z^N(t)$,
where 
\begin{align*}
\widehat Y^N(t)&=Y^N_k
+AS_{\delta t}Y^N_k(t-k\delta t)
+S_{\delta t}P^NF(Y^N_{k+1}+Z^N_{k+1})(t-k\delta t)\\
\widehat Z^N(t)&=Z^N_k+AS_{\delta t} Z^N_k(t-k\delta t)+S_{\delta t}P^N(W(t)-W(t_k)).
\end{align*}
Now we show that $\|\widehat Y^N(t)\|_{\HH^1}$ and $\|\widehat Z^N(t)\|_E$ have any finite $q$th moment, $q\ge 2$.
Indeed, we have 
\begin{align*}
\|\widehat Y^N(t)\|_{\HH^1}&\le C\|Y^N_k\|_{\HH^1}
+C(1+\|Y^N_{k+1}\|_{\HH^1}^3+\|Z^N_{k+1}\|_{E}^3)\\
\|\widehat Z^N(t)\|_{E}&\le C\|Z^N_k\|_E+\|S_{\delta t}P^N(W(t)-W(t_k))\|_E,
\end{align*}
which, together with the estimations in Lemmas \ref{pri-z} and \ref{pri-y}, implies the boundedness of any $q$th moment
of $\|\widehat Y^N(t)\|_{\HH^1}$ and $\|\widehat Z^N(t)\|_E$. 
Similar arguments in Lemmas \ref{pri-z} and \ref{pri-y} yield that for $s\in [0,T-t]$
\begin{align*}
\|Z^M(s,\widehat Z^N(t))\|_E&\le 
\|\widehat Z^N(t)\|_E+\|\int_{t}^{t+s}S(s-r)d\widehat W(r)\|_E,\\
\|Y^M(s,\widehat Y^N(t))\|_E&\le \|\widehat Y^N(t)\|_E+C_d(\|\widehat Y^N(t)\|_{\HH^1},Q,\{\|Z^M(r,\widehat Z^N(t))\|_E\}_{r\in[0,s]}),
\end{align*}
where $C_d(\|\widehat Y^N(t)\|_{\HH^1},Q,\{\|Z^M(r,\widehat Z^N(t))\|_E\}_{r\in [0,s]})$ is a random variable and polynomially depends on $\|\widehat Y^N(t)\|_{\HH^1}$ and $\|\widehat Z^N(t)\|_E$. Similar to the proof of Lemma 
\ref{pri-y}, in $d=1$, we do not need the a priori estimate of $\|Y^M\|_{\HH^1}$.
Combining with the a priori estimate of stochastic convolution of $Z^M$, we deduce that $C_d$ has finite $q$th moment, $q\ge 2$, which leads to the desired result.
\qed

\textbf{Proof of Theorem \ref{weak-dis0}}
Let $K\delta t=T>0$.
We transform the error estimate from $\HH$
into $P^M (\HH)$,
\begin{align*}
|\E[\phi(X(K\delta t,X_0))-\phi(X_k^N)]|
&\le |\E[\phi(X(K\delta t,X_0))-\phi(X^M(K\delta t,X_0^M))]|\\
&\quad+|\E[\phi(X^M(K\delta t,X_0^M))-\phi(X_K^N)]|.
\end{align*}
From the strong convergence analysis in Lemma \ref{str} and Remark \ref{rk-str}, it follows that for $M\in \N^+$,
\begin{align*}
|\E[\phi(X(K\delta t,X_0))-\phi(X^M(K\delta t,X_0^M))]|
&\le C(K\delta t,X_0)\lambda_M^{-\frac \beta 2}.
\end{align*}
Then, after taking $M\to \infty$, it suffices to estimate $|\E[\phi(X^M(K\delta t,X_0^M))-\phi(X_K^N)]|$. 
We decompose $\E[\phi(X^M(K\delta t,X_0^M))-\phi(X_K^N)]$ as
\begin{align*}
\E\Big[U^{M}(K\delta t,X_0^M)\Big]-\E\Big[U^{M}(0,X^N_K)\Big]
&=\Big(\E\Big[U^{M}(K\delta t,X_0^M)\Big]-\E\Big[U^{M}(K\delta t,X_0^N)\Big]\Big)\\
&\quad+\Big(\E\Big[U^{M}(K\delta t,X^N_0)\Big]-\E\Big[U^{M}(0,X^N_K)\Big]\Big).
\end{align*}
The regularity estimate  of $U^{M}$ in Lemma \ref{reg-um} leads to
\begin{align*}
&\Big|\E\Big[U^{M}(K\delta t,X_0^M)\Big]-\E\Big[U^{M}(K\delta t,X_0^N)\Big]\Big|\\
&\le
\int_0^1\Big|\E\Big[DU^{M}(K\delta t,\theta X_0^M+(1-\theta)X_0^N)\cdot(I-P^N)X_0^M\Big]\Big|d\theta\\
&\le C(1+\|X_0^M\|_E^2+\|X_0^N\|_E^2)\min((1+{(K\delta t)}^{-\alpha})e^{-cK\delta t}\lambda_N^{-\alpha}\|X_0\|,
\lambda_N^{-\frac \beta 2}\|X_0\|_{\HH^{\beta}}).
\end{align*} 
From the It\^o formula for Skorohod integrals, the Kolmogorov equation \eqref{kol-um} and Malliavin integration by parts, it follows that 
\begin{align*}
&\E\Big[U^{M}(K \delta t,X^N_0)\Big]-\E\Big[U^{M}(0,X^N_K)\Big]\\
&= \sum_{k=0}^{K-1}
\E\Big[U^{M}(K\delta t-k\delta t,X^N_k)\Big]-\E\Big[U^{M}(K\delta t-(k+1)\delta t,X^N_{k+1})\Big]\\
&= 
\E\Big[U^{M}(K\delta t,X^N_0)\Big]-\E\Big[U^{M}(K\delta t-\delta t,X^N_{1})\Big]
\\
&\quad-\sum_{k=0}^{K-1}\int_{t_k}^{t_{k+1}}\sum_{l\in \N^+}\E\Big[D^2U^{M}(T-t, \widehat X^N(t))\cdot(\mathcal D_t \widehat X^N(t)P^NQ^{\frac 12}e_l,S_{\delta t}P^NQ^{\frac 12}e_l)\Big]dt\\
&\quad+
\sum_{k=0}^{K-1}
\Big(\int_{t_k}^{t_{k+1}}\E\Big[\<DU^{M}(T-t,\widehat X^N(t)),A\widehat X^N(t)-AS_{\delta t}X_k^N\>\Big]dt\Big)\\
&\quad+\Big(\int_{t_k}^{t_{k+1}}\E\Big[\<DU^{M}(T-t,\widehat X^N(t)), P^M F(\widehat X^N(t))-S_{\delta t}P^NF(X_{k+1}^N)\>\Big]dt\Big)\\
&\quad+\frac 12\Big(\int_{t_k}^{t_{k+1}}\sum_{l\in \N^+}\E\Big[D^2U^{M}(T-t,\widehat X^N(t))\cdot\Big((P^MQ^{\frac 12}e_l, P^MQ^{\frac 12}e_l)-
(S_{\delta t}P^NQ^{\frac 12}e_l, S_{\delta t}P^NQ^{\frac 12}e_l)\Big)\Big]dt\Big)\\
&=:\E\Big[U^{M}(K\delta t,X^N_0)\Big]-\E\Big[U^{M}(K\delta t-\delta t,X^N_{1})\Big]+\sum_{k=1}^{K-1} II_1^k+II_2^k+II_3^k+II_4^k.
\end{align*}
The estimation for the first term in the above equation can be easily obtained by
the similar arguments in the proof of Theorem \ref{weak} and thus we focus on the estimations of $II_1^k$-$II_4^k$, $k\ge1$.
By the regularity estimate of $U^{M}$ in Lemmas \ref{reg-um} and \ref{reg-um1} and the a priori estimate of $\widehat X^N$, we have
\begin{align*}
|\sum_{k=1}^{K-1} II_1^k|
&\le \Big|\sum_{k=0}^{K-1}\int_{t_k}^{t_{k+1}}\sum_{l\in \N^+}\E\Big[D^2U^{M}(T-t, \widehat X^N(t))\cdot(\mathcal D_t\widehat X^N(t)P^NQ^{\frac 12}e_l,S_{\delta t}P^NQ^{\frac 12}e_l)\Big]dt\Big|\\
&\le C\sum_{k=1}^{K-1}\int_{t_k}^{t_{k+1}}(1+(T-t)^{\frac {\beta -1}2})e^{-c(T-t)}
\E\Big[(1+\sup_{s\in [0,T-t]}\E_{\widehat  X^N(t)}[\|X^M(s,\widehat X^N(t))\|^{14}_E])\\
&\qquad\|\mathcal D_t\widehat X^N(t)\|_{\LL_2^0}\|(-A)^{\frac {1-\beta}2}S_{\delta t}\|_{\LL}\|(-A)^{\frac {\beta-1}2}\|_{\LL_2^0}\Big]dt
\\
&\le C(Q,X_0)\delta t^{\frac {\beta+1}2}\sum_{k=1}^{K-1}\int_{t_k}^{t_{k+1}}(1+(T-t)^{\frac {\beta-1}2})e^{-c(T-t)}(t_{k+1}-[t]_{\delta t})^{\frac {\beta-1}2}dt\\
&\le C(Q,X_0)\delta t^{\beta},
\end{align*}
where  we use the a priori estimate in Proposition \ref{mal-dif} and the  fact that for $t_k\le  t\le s \le t_{k+1}$,
\begin{align*}
\mathcal D_s\widehat X^N(t)
&=S_{\delta t}\mathcal D_sX^N(t_k)
+(t-t_k)P^NS_{\delta t}DF(\widehat X^N(t_{k+1}))\mathcal D_s \widehat X^N_{k+1}
+\mathcal D_s \int_{t_k}^{t}S_{\delta t}P^NdW(s)\\
&=(t-t_k)P^NS_{\delta t}DF(\widehat X^N(t_{k+1}))\mathcal D_s \widehat X^N_{k+1}.
\end{align*}
Then we estimate $II_2^k, II_3^k$ and $II_4^k$, $k\ge 1$ separately.  
The definition of $\widehat X$ leads to
\begin{align*}
II_2^k
&= \int_{t_k}^{t_{k+1}}\E\Big[\<DU^{M}(T-t,\widehat X^N(t)),A( X^N_k-S_{\delta t}X_k^N)\>\Big]dt\\
 &\quad+\int_{t_k}^{t_{k+1}}\E\Big[\<DU^{M}(T-t,\widehat X^N(t)),(t-t_k)A^2S_{\delta t}X_k^N\>\Big]dt\\
 &\quad+ \int_{t_k}^{t_{k+1}}\E\Big[\<DU^{M}(T-t,\widehat X^N(t)),
(t-t_k)AS_{\delta t}P^NF(X_{k+1}^N)\>\Big]dt\\
&\quad+ \int_{t_k}^{t_{k+1}}\E\Big[\<DU^{M}(T-t,\widehat X^N(t)),A
\int_{t_k}^tS_{\delta t}dW(s)\>\Big]dt
\\
&:=II_{21}^k+II^k_{22}+II^k_{23}+II^k_{24}.
\end{align*}
From $I-S_{\delta t}=-A\delta t(I-A\delta t)^{-1}$, the mild form of $X_k^N$ \eqref{full-spde}, Lemma \ref{pri-y} and  Lemmas \ref{reg-um}-\ref{pri-con},
it follows that for $k\ge 1$ and any small $\epsilon_1>0$
\begin{align*}
|II^k_{21}|
&\le C\delta t\int_{t_k}^{t_{k+1}} (1+(T-t)^{-\alpha})e^{-c(T-t)}
\E \Big[(1+\sup_{s\in [0,T-t]}\E_{\widehat X^N(t)}[\|X^M(s,\widehat X^N(t))\|_E^4])\\
&\qquad \|(-A)^{1-\epsilon_1}S_{\delta t}^{k}\|\|(-A)^{1-\alpha+\epsilon_1}S_{\delta t}\|
\|X_0^N\|\Big]dt
\\
&\quad+C\delta t\int_{t_k}^{t_{k+1}}(1+(T-t)^{-\alpha})e^{-c(T-t)}\sum_{j=0}^{k-1}\E\Big[(1+\sup_{s\in [0,T-t]}\E_{\widehat X^N(t)}[\|X^M(s,\widehat X^N(t))\|_E^4])\\
&\qquad \|(-A)^{1-\epsilon_1}S_{\delta t}^{k-j}\|\|(-A)^{1-\alpha+\epsilon_1}S_{\delta t}\|
\|F(X^N_{j+1})\|\Big]dt\\
&\quad+
\Big| \int_{t_k}^{t_{k+1}}\E\Big[\<DU^{M}(T-t,\widehat X^N(t)),A^2\delta t\sum_{j=0}^{k-1}S_{\delta t}^{k+1-j}P^N\delta W_j\>\Big]dt\Big|\\
&\le C(X_0,Q)\delta t^{\alpha-\epsilon_1}\int_{t_k}^{t_{k+1}} (1+(T-t)^{-\alpha})e^{-c(T-t)} (1+(t_k)^{-1+\epsilon_1}e^{-c_1 t_k})dt+II_{211}^k,
\end{align*}
where $II_{211}^k:=\Big| \int_{t_k}^{t_{k+1}}\E\Big[\<DU^{M}(T-t,\widehat X^N(t)),A^2\delta t\sum_{j=0}^{k-1}S_{\delta t}^{k+1-j}P^N\delta W_j\>\Big]dt\Big|$.
Notice that the lack of regularity and bad time behavior  do not happen at the same time.
We split $II_{211}^k$  as 
\begin{align*}
|II_{211}^k|
&\le \Big| \int_{t_k}^{t_{k+1}}\E\Big[\<DU^M(T-t,\widehat X^N(t)),A^2\delta t\int_0^{\max(t_k-1,0)}P^NS_{\delta t}^{k+1-\lfloor s \rfloor}dW(s)\>\Big]dt\Big|\\
&\quad+ \Big| \int_{t_k}^{t_{k+1}}\E\Big[\<DU^M(T-t,\widehat X^N(t)),A^2\delta t\int_{\max(t_k-1,0)}^{t_k}P^NS_{\delta t}^{k+1-\lfloor s \rfloor}dW(s)\>\Big]dt\Big|.
\end{align*}
 The Cauchy--Schwarz inequality, the regularity estimate of $U^M$, a priori estimate of $X^N$ and the smoothy effect of  $S_{\delta t}$  and Lemma \ref{pri-con}  yield that 
\begin{align*}
&\Big| \int_{t_k}^{t_{k+1}}\E\Big[\<DU^M(T-t,\widehat X^N(t)),A^2\delta t\int_0^{\max(t_k-1,0)}P^NS_{\delta t}^{k+1-\lfloor s \rfloor}dW(s)\>\Big]dt\Big|\\
&\le 
C \delta t \int_{t_k}^{t_{k+1}}
\sqrt{\E\Big[\|\int_0^{\max(t_k-1,0)}A^2P^NS_{\delta t}^{k+1-\lfloor s \rfloor}dW(s)\|^2\Big]}\sqrt{\E\Big[\|DU^M(T-t,\widehat X^N(t)\|^2\Big]}dt\\
&\le 
C(Q,X_0)\delta t\int_{t_k}^{t_{k+1}}\sqrt{\E\Big[\|\int_0^{\max(t_k-1,0)}A^2P^NS_{\delta t}^{k+1-\lfloor s \rfloor}dW(s)\|^2\Big]}e^{-c(T-t)}dt\\
&\le 
C(Q,X_0)\delta t\int_{t_k}^{t_{k+1}}\sqrt{\E\Big[\int_0^{\max(t_k-1,0)}\|(-A)^{2+\frac {1-\beta}2}P^NS_{\delta t}^{k+1-\lfloor s \rfloor}\|^2_{\mathcal L(\HH)}\|(-A)^{\frac {\beta-1}2}\|^2_{\LL_2^0}ds\Big]}e^{-c(T-t)}dt\\
&\le C(Q,X_0)\delta t\int_{t_k}^{t_{k+1}}\sqrt{\int_0^{\max(t_k-1,0)}
\frac 1{(1+\lambda_1\delta t)^{k-\lfloor s \rfloor}}(t_k-[s]_{\delta t})^{-5+\beta}ds}e^{-c(T-t)}dt\\
&\le C(Q,X_0)\delta t\int_{t_k}^{t_{k+1}}\sqrt{\int_0^{\infty}
\frac 1{(1+\lambda_1\delta t)^{\lfloor s \rfloor}}ds}e^{-c(T-t)}dt\le C(Q,X_0)\delta t\int_{t_k}^{t_{k+1}}e^{-c(T-t)}dt.
\end{align*}
Applying Malliavin calculus integration by parts, Malliavin differentiability of $\widehat X^N$ and the regularity estimate of $U^M$ and Lemma \ref{pri-con}, we have 
\begin{align*}
&\Big| \int_{t_k}^{t_{k+1}}\E\Big[\<DU^{M}(T-t,\widehat X^N(t)),A^2\delta t\int_{\max(t_k-1,0)}^{t_k}P^NS_{\delta t}^{k+1-\lfloor s \rfloor}dW(s)\>\Big]dt\Big|\\
&= \delta t\int_{t_k}^{t_{k+1}}\int_{\max(t_k-1,0)}^{t_k} \sum_{l\in\N^+} \E\Big[\Big| D^2U^{M}(T-t,\widehat X^N(t))\cdot (\mathcal D_s^{Q^{\frac12 }e_l}\widehat X^N(t),A^2S_{\delta t}^{k+1-\lfloor s \rfloor}P^NQ^{\frac 12}e_l)\Big|\Big]dsdt\\
&\le C\delta t\int_{t_k}^{t_{k+1}}
\int_{\max(t_k-1,0)}^{t_k} 
\sum_{l\in\N^+}  \E\Big[\Big| \< (-A)^{\frac {1+\beta}2-\epsilon_1}D^2U^{M}(T-t,\widehat X^N(t)) (-A)^{\frac {1-\beta}2}\\
&\qquad (-A)^{\frac {\beta-1}2} \mathcal D_s^{Q^{\frac12 }e_l}\widehat X^N(t),(-A)^{2-\frac {\beta+1}2+\epsilon_1}S_{\delta t}^{k+1-\lfloor s \rfloor}P^NQ^{\frac 12}e_l)\>\Big|\Big]dsdt
\\
&\le C\delta t\int_{t_k}^{t_{k+1}}e^{-c(T-t)}(1+(T-t)^{-1+\epsilon_1})
\int_{\max(t_k-1,0)}^{t_k} 
\E\Big[(1+\sup_{r\in [0,T-t]}\E_{\widehat X^N(t)}\|X^M(r,\widehat X^N(t))\|_E^{14}])\\
&\qquad\|(-A)^{\frac {\beta-1}2}\mathcal D_s \widehat X^N(t)\|_{\LL_2^0}
\|(-A)^{1-\epsilon_1} S_{\delta t}^{k-\lfloor s \rfloor}\|
\|(-A)^{1-\beta+2\epsilon_1}S_{\delta t}\|\|(-A)^{\frac {\beta-1}2}\|_{\LL_2^0}\Big]dsdt\\
&\le C(X_0,Q)\delta t^{\beta-2\epsilon_1}\int_{t_k}^{t_{k+1}}e^{-c(T-t)}(1+(T-t)^{-1+\epsilon_1})\int_{\max(t_k-1,0)}^{t_k} (t_{k}-[s]_{\delta t})^{-1+\epsilon_1}e^{-c(t_{k}-[s]_{\delta t})}\\
&\qquad\|(-A)^{\frac {\beta-1}2}\mathcal D_s \widehat X^N(t)\|_{L^2(\Omega;\LL_2^0)} dsdt\\
&\le C(X_0,Q)\delta t^{\beta-2\epsilon_1}\int_{t_k}^{t_{k+1}}e^{-c(T-t)}(1+(T-t)^{-1+\epsilon_1})dt.
\end{align*} 
The above analysis leads to 
\begin{align*}
|II^k_{21}|&\le C(X_0,Q)\delta t^{\alpha-\epsilon_1}\int_{t_k}^{t_{k+1}} (1+(T-t)^{-\alpha})e^{-c(T-t)}dt\\
&\quad+C(X_0,Q)\delta t^{\beta-2\epsilon_1}\int_{t_k}^{t_{k+1}}(1+(T-t)^{-1+\epsilon_1})e^{-c(T-t)}dt,
\end{align*}
for $k\ge 1$.
Since the estimations for $II^k_{22}$ and $II^k_{23} $ for $k\ge 1$ are similar,
we omit the procedures.
 Malliavin  integration by parts
yields that
\begin{align*}
|II^k_{24}|&=
\Big|\int_{t_k}^{t_{k+1}}\E\Big[\<DU^{M}(T-t,\widehat X^N(t)),A
\int_{t_k}^tS_{\delta t}dW(s)\>\Big]dt\Big|\\
&=\Big|\int_{t_k}^{t_{k+1}}\int_{t_k}^t\E\Big[\<D^2U^{M}(T-t,\widehat X^N(t))\mathcal D_s\widehat X^N(t),A
S_{\delta t}\>_{\LL_2^0}\Big]dsdt\Big|\\
&\le \Big|\int_{t_k}^{t_{k+1}}\int_{t_k}^t\E\Big[\<(-A)^{\frac {1+\beta}2-\epsilon_1}D^2U^{M}(T-t,\widehat X^N(t))(-A)^{\frac {1-\beta}2} (-A)^{\frac {\beta -1}2}\mathcal D_s\widehat X^N(t), \\
&\quad (-A)^{\epsilon_1}
S_{\delta t}(-A)^{\frac {1-\beta}2}\>_{\LL_2^0}\Big]dsdt\Big|\\
&\le C(Q,X_0)\delta t^{\beta-\epsilon_1} 
\int_{t_k}^{t_{k+1}}(1+(T-t)^{-1+\epsilon_1})e^{-c(T-t)}dt.
\end{align*}
It follows that 
\begin{align*}
|II^k_2|&\le C(X_0,Q)\delta t^{\beta-2\epsilon_1}\int_{t_k}^{t_{k+1}}e^{-c(T-t)}(1+(T-t)^{-1+\epsilon_1})dt\\
&\quad +C(X_0,Q)\delta t^{\alpha-\epsilon_1}\int_{t_k}^{t_{k+1}} (1+(T-t)^{-\alpha})e^{-c(T-t)} (1+(t_k)^{-1+\epsilon_1}e^{-c_1 t_k})dt.
\end{align*}
Now, we are in a position to estimate $II_3^k$.
By the  regularity of $DU^{M}$  and a priori estimate of $\widehat X^N$, we have 
\begin{align*}
|II_3^k|&\le
\Big|\int_{t_k}^{t_{k+1}}\E\Big[\<DU^{M}(T-t,\widehat X^N(t)), P^M(I-P^N)F(\widehat X^N(t))\>\Big]dt\Big|\\
&\quad+\Big|\int_{t_k}^{t_{k+1}}\E\Big[\<DU^{M}(T-t,\widehat X^N(t)), (I-S_{\delta t})P^NF(X_{k+1}^N)\>\Big]dt\Big|\\
&\quad+ \Big|\int_{t_k}^{t_{k+1}}\E\Big[\<DU^{M}(T-t,\widehat X^N(t)),P^N(F(\widehat X^N(t))-F(X_{k+1}^N)\>\Big]dt\Big|\\
&\le 
C(\lambda_N^{-\alpha}+\delta t^{\alpha})
\int_{t_k}^{t_{k+1}}(1+(T-t)^{-\alpha})e^{-c(T-t)}\\&\qquad\E\Big[(1+\sup_{r\in [0,T-t]}\E_{\widehat X^N(t)}[\|X^M(r,\widehat X^N(t))\|_E^4])(1+\|\widehat X(t)\|_{L^6}^3)\Big]dt\\
&\quad+ \Big|\int_{t_k}^{t_{k+1}}\E\Big[\<DU^{M}(T-t,\widehat X^N(t)),P^N(F(\widehat X^N(t))-F(X_{k+1}^N)\>\Big]dt\Big|.
\end{align*} 
Thus it  sufficient to estimate the last term in the above inequality.
It follows from 
Taylor expansion of $F$, the regularity of $\widehat X^N$ and $DU^{M}$, and the a priori estimate of $\widehat X^N$  that 
\begin{align*}
&\int_{t_k}^{t_{k+1}}\E\Big[\<DU^{M}(T-t,\widehat X^N(t)),P^N(F(\widehat X^N(t))-F(X_{k+1}^N)\>\Big]dt\\
&\le \int_{t_k}^{t_{k+1}}(t-t_{k+1})\E\Big[\<DU^{M}(T-t,\widehat X^N(t)),P^N(DF(\widehat X^N(t))\cdot (AS_{\delta t}X_{k}^N))\>\Big]dt\\
&\quad+\int_{t_k}^{t_{k+1}}(t-t_{k+1})\E\Big[\<DU^{M}(T-t,\widehat X^N(t)),P^N(DF(\widehat X^N(t))\cdot (S_{\delta t}P^NF(X_{k+1}^N))\>\Big]dt\\
&\quad+\int_{t_k}^{t_{k+1}}\E\Big[\<DU^{M}(T-t,\widehat X^N(t)),P^N(DF(\widehat X^N(t))\cdot P^N(\int_t^{t_{k+1}}S_{\delta t}dW(s)))\>\Big]dt\\
&\quad+\int_{t_k}^{t_{k+1}}\E\Big[\<DU^{M}(T-t,\widehat X^N(t)),P^N(\int_0^1D^2F(\theta \widehat X^N(t)+(1-\theta )X^N_{k+1})) \\
&\qquad\cdot (\widehat X^N(t)-X^N_{k+1},\widehat  X^N(t)-X^N_{k+1})(1-\theta)d\theta\>)\Big]dt
=:II^k_{31}+II^k_{32}+II^k_{33}+II^k_{34}.
\end{align*}
The estimation of $II_{31}^k$ is similar to the estimation of $II_{21}^k$ and we need to use a proper decomposition of the stochastic integral.
The mild form of $X_{k+1}^N$, Malliavin  integration by parts and Lemma \ref{pri-con}  yield that
\begin{align*}
|II^k_{31}|
&:=\Big|\int_{t_k}^{t_{k+1}}(t-t_{k+1})\E\Big[\<DU^{M}(T-t,\widehat X^N(t)),P^N(DF(\widehat X^N(t))\cdot (AS_{\delta t}^{k+1}X_{0}^N))\>\Big]dt\\
&\quad+\int_{t_k}^{t_{k+1}}(t-t_{k+1})\delta t\E\Big[\<DU^{M}(T-t,\widehat X^N(t)),P^N(DF(\widehat X^N(t))\cdot (\sum_{j=0}^{k-1}A\\
&\qquad S_{\delta t}^{k+1-j}P^NF(X_{j+1}^N))\>\Big]dt+II_{311}^k\Big|\\
&\le 
C\delta t \int_{t_k}^{t_{k+1}}e^{-c(T-t)}\E\Big[(1+\sup_{r\in [0,T-t]}\E_{\widehat X^N(t)}[\|X^M(r,\widehat X^N(t))\|^4_E])(1+\|\widehat X^N(t)\|^2_E)\\
&\qquad\|(-A)^{1-\epsilon_1}S_{\delta t}^k\|\|(-A)^{\epsilon_1}S_{\delta t}\|\|X_0^N\|\Big]dt
\\
&\quad+C\delta t^2
\int_{t_k}^{t_{k+1}} e^{-c(T-t)}\sum_{j=0}^{k-1}\E\Big[(1+\sup_{r\in [0,T-t]}\E_{\widehat X^N(t)}[\|X^M(r,\widehat X^N(t))\|^4_E])(1+\|\widehat X^N(t)\|^5_E)\\
&\qquad
\|(-A)^{1-\epsilon_1}S_{\delta t}^{k-j}\|(-A)^{\epsilon_1} S_{\delta t}\|\|X_0^N\|\Big]dt+|II_{311}^k|\\
&\le
C(Q,X_0) \int_{t_k}^{t_{k+1}}e^{-c(T-t)}dt\delta t^{1-\epsilon_1}(t_k^{-1+\epsilon_1}e^{-c_1t_k}
+\delta t\sum_{j=0}^{k-1}t_{k-j}^{-1+\epsilon_1}e^{-c_1{(t_k-t_j)}})+|II_{311}^k|.
\end{align*}
Similar  arguments on estimating  $II_{211}^k$ leads to
\begin{align*}
|II_{311}^k|
&\le \Big| \int_{t_k}^{t_{k+1}}(t_{k+1}-t)\E\Big[\<DU^M(T-t,\widehat X^N(t)),A \int_0^{\max(t_k-1,0)}P^NS_{\delta t}^{k+1-\lfloor s \rfloor}dW(s)\>\Big]dt\Big|\\
&\quad+ \Big| \int_{t_k}^{t_{k+1}}(t_{k+1}-t)\E\Big[\<DU^M(T-t,\widehat X^N(t)),A\int_{\max(t_k-1,0)}^{t_k}P^NS_{\delta t}^{k+1-\lfloor s \rfloor}dW(s)\>\Big]dt\Big|.
\end{align*}
The regularity estimate of  $DU^M$ yields that
\begin{align*}
 &\Big|\int_{t_k}^{t_{k+1}}(t_{k+1}-t)\E\Big[\<DU^M(T-t,\widehat X^N(t)),A \int_0^{\max(t_k-1,0)}P^NS_{\delta t}^{k+1-\lfloor s \rfloor}dW(s)\>\Big]dt\Big|\\
 &\le 
 C\delta t
 \int_{t_k}^{t_{k+1}}\sqrt{\E[\|DU^M(T-t,\widehat X^N(t))\|]}
 \sqrt{\E[\|\int_0^{\max(t_k-1,0)}AP^NS_{\delta t}^{k+1-\lfloor s \rfloor}dW(s)\|^2]}dt\\
 &\le C(X_0,Q) \delta t\int_{t_k}^{t_{k+1}}e^{-c(T-t)}dt.
\end{align*}
From  the smoothing effect  of $S_{\delta t}$, the Malliavin regularity and the a priori estimate of $\widehat X(t)$, and the  Sobolev embedding theorem $E \hookrightarrow \HH^{\frac d2+\epsilon_1}, \epsilon_1>0 $, it follows  that 
\begin{align*}
&\Big| \int_{t_k}^{t_{k+1}}(t_{k+1}-t)\E\Big[\<DU^M(T-t,\widehat X^N(t)),A\int_{\max(t_k-1,0)}^{t_k}P^NS_{\delta t}^{k+1-\lfloor s \rfloor}dW(s)\>\Big]dt\Big|\\
&\le
 C(Q,X_0)\delta t\int_{t_k}^{t_{k+1}}e^{-c(T-t)}
\int_{\max(t_k-1,0)}^{t_k}
\sum_{l\in \N^+}\|AS_{\delta t}^{k+1-\lfloor s \rfloor}Q^{\frac 12}e_l\|  \|\mathcal D_s^{Q^{\frac 12}e_l}\widehat X(t)\|_{L^2(\Omega;\HH)}dsdt \\
&\quad+ C(Q,X_0)\delta t\int_{t_k}^{t_{k+1}}e^{-c(T-t)}
\int_{\max(t_k-1,0)}^{t_k}
\sum_{l\in \N^+}\Big\|(-A)^\eta DU^M(T-t,\widehat X(t))\Big\|_{L^2(\Omega;\HH)} \Big\|(-A)^{-\eta}D^2F(\widehat X(t)\\
&\qquad \cdot(P^NAS_{\delta t}^{k+1-\lfloor s \rfloor}Q^{\frac 12}e_l,\mathcal D_s^{Q^{\frac 12}e_l}(\widehat X(t)))\Big\|_{L^2(\Omega;\HH)}dsdt\\
&\le
 C(Q,X_0)\delta t\int_{t_k}^{t_{k+1}}e^{-c(T-t)}
\int_{\max(t_k-1,0)}^{t_k}
\|(-A)^{\frac {3-\beta}2}S_{\delta t}^{k+1-\lfloor s \rfloor}(-A)^{\frac {
\beta -1}2}\|_{\LL_2^0} 
\|\mathcal D_s\widehat X(t)\|_{L^2(\Omega; \LL_2^0)}dsdt \\
&\quad+ C(Q,X_0)\delta t\int_{t_k}^{t_{k+1}}e^{-c(T-t)}(1+(T-t)^{-\eta})
\int_{\max(t_k-1,0)}^{t_k}
\|(-A)^{\frac {3-\beta}2}S_{\delta t}^{k+1-\lfloor s \rfloor}(-A)^{\frac {\beta-1}2}\|_{\LL_2^0}\\
&\qquad \sqrt{\sum_{l\in \N^+}\|\mathcal D_s^{Q^{\frac 12}e_l}\widehat X(t)\|^2_{L^4(\Omega;\HH)}}dsdt\\
&\le
C(Q,X_0)\delta t^{\beta-\epsilon_1}
\int_{t_k}^{t_{k+1}}e^{-c(T-t)}(1+(T-t)^{-\eta})
\int_{0}^{t_{k}}(t_{k}-s)^{-1 +\epsilon_1}e^{-c_1(t_k-s)}dsdt,
\end{align*}
for $\eta >\frac d4+\frac {\epsilon_1}2$.
Similarly, we have 
\begin{align*}
|II^k_{32}|
&\le C\delta t\int_{t_{k}}^{t_{k+1}}e^{-c(T-t)}
\E\Big[(1+\sup_{s\in [0,T-t]}\E_{\widehat X^N(t)}\|X^M(s,\widehat X^N(t))\|_E^4)(1+\|\widehat X^N(t)\|_{E}^5)\Big]dt
\end{align*}
and
\begin{align*}
|II^k_{33}|&\le
\Big|\int_{t_k}^{t_{k+1}}\int_t^{t_{k+1}}\sum_{l\in \N^+}\E\Big[D^2U^{M}(T-t,\widehat X^N(t))\cdot\Big(P^N(DF(\widehat X^N(t))\cdot (P^NS_{\delta t}Q^{\frac 12}e_l), \mathcal D_s^{Q^ {\frac 12e_l}}\widehat X^N(t)\Big)\Big]dsdt
\Big|\\
&\quad+
\Big|\int_{t_k}^{t_{k+1}}\int_t^{t_{k+1}}\sum_{l\in \N^+}\E\Big[\<DU^{M}(T-t,\widehat X^N(t)),P^N(D^2F(\widehat X^N(t))\cdot (\mathcal D_s^{Q^{\frac 12}e_l} \widehat X^N(t),P^NS_{\delta t}Q^{\frac 12}e_l )\>\Big]dsdt\Big|
\\
&\le 
C(Q,X_0)\delta t\int_{t_k}^{t_{k+1}} (1+(T-t)^{-\eta})e^{-c(T-t)}dt,
\end{align*}
where  we utilize the fact that for $t_k\le  t\le s \le t_{k+1}$,
\begin{align*}
\mathcal D_s\widehat X^N(t)
&=S_{\delta t}\mathcal D_sX^N(t_k)
+(t-t_k)P^NS_{\delta t}F(\widehat X^N(t_{k+1}))\mathcal D_s \widehat X^N_{k+1}
+\mathcal D_s \int_{t_k}^{t}P^NS_{\delta t}dW(s)\\
&=(t-t_k)P^NS_{\delta t}F(\widehat X^N(t_{k+1}))\mathcal D_s \widehat X^N_{k+1}.
\end{align*}
Combining with the continuity of $\widehat X^N$, for $t\in [t_k,t_{k+1}]$,
\begin{align*}
\|\widehat X^N(t)-\widehat X^N_{k+1}\|_{L^p(\Omega;\HH)}
&\le 
(t_{k+1}-t)\|(-A)^{1-\frac \beta 2}S_{\delta t}\|_{\LL}
\|X(t_k)\|_{L^p(\Omega;\HH^{ \beta })}
+C\|X(t_k)-X^N_k\|_{L^p(\Omega;\HH)}\\
&\quad+C\|F(X_{k+1}^N)\|_{L^p(\Omega;\HH)}(t_{k+1}-t)
+\|\int_{t}^{t_{k+1}}P^NS_{\delta t}dW(s)\|_{L^p(\Omega;\HH)}\\
&\le C(X_0,Q)(t_{k+1}-t)^{\frac \beta 2},
\end{align*}
we deduce that
\begin{align*}
|II^k_{34}|&\le 
C\Big|\int_{t_k}^{t_{k+1}}\E\Big[\<DU^{M}(T-t,\widehat X^N(t)),P^N(\int_0^1D^2F(\theta \widehat X^N(t)+(1-\theta) X^N_{k+1})) \\
&\qquad\cdot (\widehat X^N(t)-\widehat X^N_{k+1},\widehat X^N(t)-\widehat X^N_{k+1})d\theta\>)\Big]dt\Big|\\
&\le C(Q,X_0)\int_{t_k}^{t_{k+1}}e^{-c(T-t)}(T-t)^{-\eta}
\|\widehat X^N(t)-\widehat X^N_{k+1}\|_{L^4(\Omega;\HH)}^2dt\\
&\le C(Q,X_0)\delta t^{\beta}\int_{t_k}^{t_{k+1}}e^{-c(T-t)}(T-t)^{-\eta}dt.
\end{align*}
Thus we conclude that
\begin{align*}
|II^k_3|&\le C(Q,X_0) \int_{t_k}^{t_{k+1}}e^{-c(T-t)}dt\delta t^{1-\epsilon_1}(t_k^{-1+\epsilon_1}e^{-c_1t_k}
+\delta t\sum_{j=0}^{k-1}t_{k-j}^{-1+\epsilon_1}e^{-c_1{(t_k-t_j)}})\\
&\quad +C(Q,X_0)\delta t^{\beta-\epsilon_1}
\int_{t_k}^{t_{k+1}}e^{-c(T-t)}(1+(T-t)^{-\eta})
\int_{0}^{t_{k}}(t_{k}-s)^{-1 +\epsilon_1}e^{-c_1(t_k-s)}dsdt.
\end{align*}
For $II_{4}^k$, by applying the regularity estimate of $D^2U^{M}$, we obtain
\begin{align*}
|II_4^k|
&\le  \Big|\int_{t_k}^{t_{k+1}}\sum_{j\in \N^+}\E\Big[D^2U^{M}(T-t,\widehat X^N(t))\cdot\Big(P^M(I-P^N)Q^{\frac 12}e_j, (P^M+P^NS_{\delta t})Q^{\frac 12}e_j\Big)\Big]dt\Big|\\
&\quad+ \Big|\int_{t_k}^{t_{k+1}}\sum_{j\in \N^+}\E\Big[D^2U^{M}(T-t,\widehat X^N(t))\cdot\Big(P^N(I-S_{\delta t})Q^{\frac 12}e_j, (P^M+P^NS_{\delta t})Q^{\frac 12}e_j\Big)\Big]dt\Big|\\
&\le C(Q,X_0) \int_{t_k}^{t_{k+1}}e^{-c(T-t)}(1+(T-t)^{-1+\epsilon_1})\|(-A)^{\frac {\beta-1}2}\|_{\LL_2^0}^2\|(-A)^{-\frac {1+\beta}2+\epsilon_1}(I-P^N)(-A)^{\frac {1-\beta}2}\|dt\\
&\quad+C(Q,X_0) \int_{t_k}^{t_{k+1}}e^{-c(T-t)}(1+(T-t)^{-1+\epsilon_1})\|(-A)^{\frac {\beta-1}2}\|_{\LL_2^0}^2\|(-A)^{-\frac {1+\beta}2+\epsilon_1}(I-S_{\delta t})(-A)^{\frac {1-\beta}2}\|dt\\
&\le C(Q,X_0)(\delta t^{\beta-\epsilon_1}+\lambda_N^{-\beta+2\epsilon_1})\int_{t_k}^{t_{k+1}}(1+(T-t)^{-1+\epsilon_1})dt.
\end{align*}
Combining all the estimation of $II_{1}^k$-$II_{4}^k$
together and summing up over $k$, we finish the proof.\\
\qed

The above time-independent error estimate, together with the $V$-uniformly ergodicity of Eq.\eqref{spde} in Proposition \ref{exp-erg1}, immediately yields the result of Corollary \ref{erg-cor}. We remark that one can first take $\delta t\to0$, and get the weak error between  $\mu^{N}$ and 
$\mu$. However, it is still unknown  the invariant measure of the 
implicit method is unique or not when
$N\to \infty$ firstly. This will be studied further.
This weak convergence approach to approximating  the invariant measure is available 
for other type numerical methods since we have given the time-independent regularity estimates of Kolmogorov equation in Lemmas \ref{reg-um} and \ref{reg-um1}. The key requirement lies on the time-independent a priori estimates of numerical solutions in $E$. 
In particular, if $d=1$,  according  to the a priori estimate in \cite{CH18} and the arguments in Lemma \ref{pri-y},
we  get the sharp weak convergence rate 
of the full discretization $\{X^h_k\}_{k\in \N^+,h\in (0,1]}$ given by the temporal implicit Euler method and the spatial linear finite element method. For convenience, denoting $V_h$ the finite element space and using the notations of  the finite element method in \cite{CH18}, we have the following result.

\begin{cor}\label{weak-fem}
Let Assumptions \ref{as-lap}-\ref{as-dri} hold with $d=1$ and $\beta\in (0,1]$, $\gamma\in (0,\beta)$,
$X_0\in E$, $T>0$, $\delta t_0\in 
(0,1\land \frac 1{(2\lambda_F-2\lambda_1)\lor 0})$. 
Then for any  $\phi\in \mathcal C_b^2(\HH)$, 
there exists  $C(X_0,Q,T)>0$ such that for any $\delta t\in (0, \delta t_0]$, $K\delta t=T $, $K\in\N^+$ and $h\in (0,1]$, 
\begin{align*}
\Big|\E\Big[\phi(X(T))-\phi(X^h_{K})\Big]\Big|
\le C(X_0,T,Q)\Big(\delta t^{\gamma}+h^{2\gamma}\Big).
\end{align*}
In addition, under Assumption \ref{erg-as} or \ref{erg-as1},  for  
any $\phi \in \mathcal C_b^2(\HH)$, there exist constants $c>0$, $C(X_0,Q)>0$ such that for any large $K$, $\delta t\in (0, \delta t_0]$ and $h\in (0,1]$,
\begin{align*}
\Big|\E\Big[\phi(X_K^h(X^h_0))-\int_{\HH}\phi d\mu\Big]\Big|
&\le C(X_0,Q)(\delta t^{\gamma}+h^{2\gamma}+e^{-cK\delta t}).
\end{align*}
Furthermore, if $\mu^{h,\delta t}$ is an ergodic invariant measure of the numerical  
solution $\{X^h_k\}_{k\in \N^+}$,
we have 
\begin{align*}
\Big|\E\Big[\int_{V_h}\phi d\mu^{h,\delta t} -\int_{\HH} \phi d\mu \Big]\Big|
&\le C(X_0,Q)(\delta t^{\gamma}+h^{2\gamma}).
\end{align*}
\end{cor}

\begin{rk}
The weak convergence analysis can be extended to the functional space $\mathcal C_p^2(\HH)$, i.e., for $\phi \in \mathcal C_p^2(\HH)$, the first and second derivatives of $\phi$ grow polynomially. 
For instance, under  Assumption \ref{erg-as1}, 
one can first use the arguments in the proof of Lemma \ref{reg-um} to get the regularity estimate of 
Kolmogorov equation in a finite time $T$.
Then similar arguments in Lemma \ref{reg-um1}
yield the exponential decay estimate for $t\ge T$ by using
the $x$-independent uniform boundedness of $X^N$ (see the estimate \eqref{y-ind}) and the Bismut--Elworthy--Li formula. 
Combining with the proof of Theorem \ref{weak-dis0}, we can obtain the similar convergence rate of 
the proposed method for  $\mathcal C_p^2(\HH)$. 
\end{rk}

\section{Numerical experiments}

In this section,  several numerical tests are presented  to verify the temporal weak convergent rates and the ergodicity of \eqref{full-spde}.
Consider
$f(\xi)=-\xi^3+\lambda_F \xi$, $\lambda_F\in \R$ and $W(t,\xi) = \sum_{j=1}^{\infty}
\frac{1}{1 + j^{\rm\kappa}}\sqrt{2}sin(j\pi \xi)
\beta_j(t)$
with $\kappa$ characterizing the smoothness of the driving  noise.
In our numerical tests,  we truncate the series by the first $M$ terms, $M\in \N^+$. 

\begin{figure}
\centering
\subfigure[$\rm \kappa=0$]{
\begin{minipage}{0.31\linewidth}
\includegraphics[width=4cm,height=4cm]{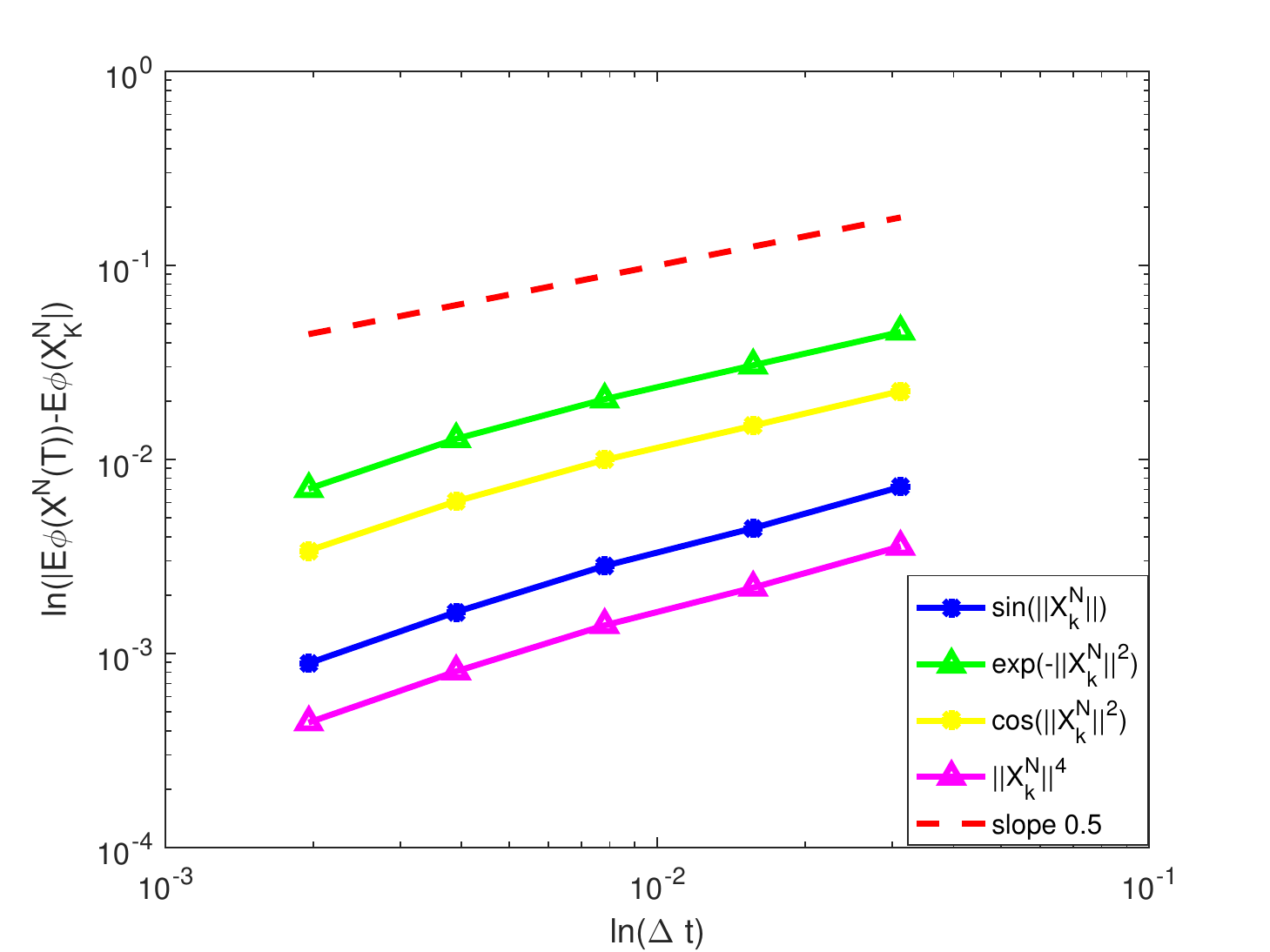}
\end{minipage}}
\subfigure[$\rm\kappa=0.5$]{
\begin{minipage}{0.31\linewidth}
\includegraphics[width=4cm,height=4cm]{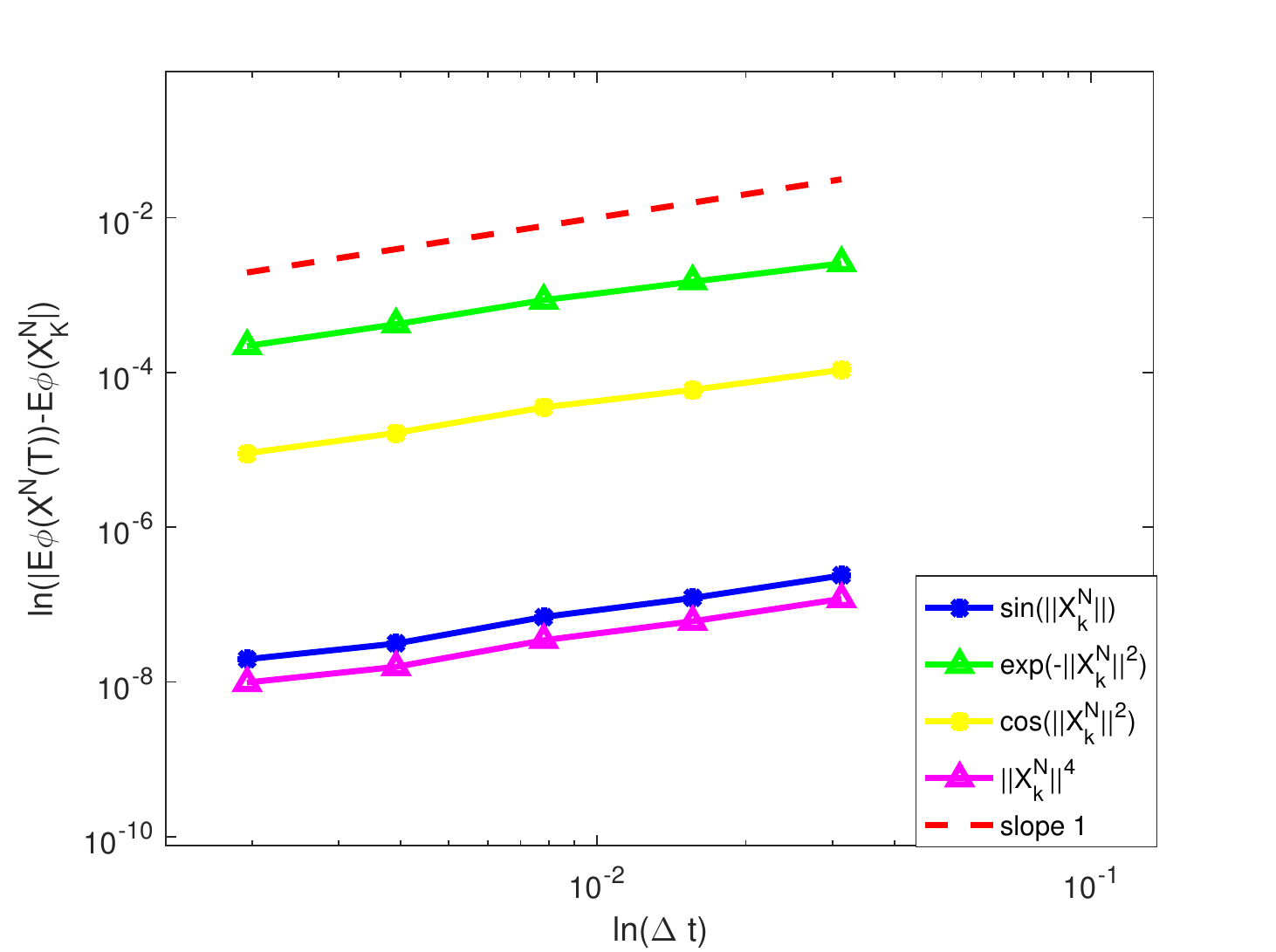}
\end{minipage}}
\subfigure[$\rm\kappa=2$]{
\begin{minipage}{0.31\linewidth}
\includegraphics[width=4cm,height=4cm]{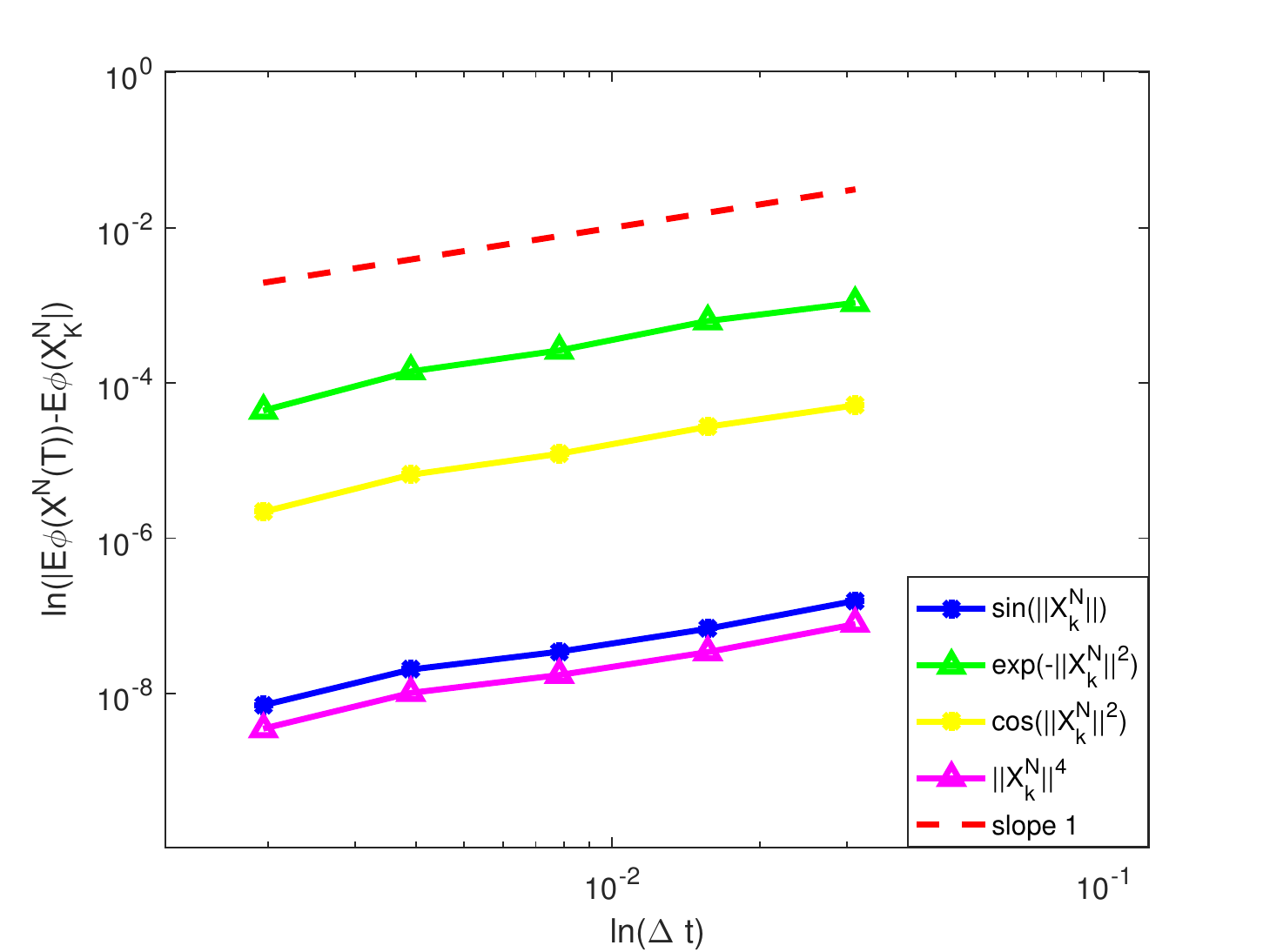}
\end{minipage}}
\caption{Rate of weak convergence in temproal direction $(u(0,x) = sin(\pi x),$  and $M=2^6)$}\label{pp1}
\end{figure}

We first investigate the weak convergence order in temporal direction of
the proposed method \eqref{full-spde}. 
In order to show the rate of weak convergence, we fix $N=2^{6}$ and  take $\delta t^{ref}=2^{-11}$ as the reference solution.
Moreover, we choose four different kinds of functionals (a) $\phi(u)=\cos(\|u\|^2),$ (b) $\phi(u)=\exp{(-\|u\|^2)},$ (c) $\phi(u)=\sin(\|u\|)$ and (d) $\phi(u)=\|u\|^4$, where $u\in \HH$, as the test functions for weak convergence.
Fig \ref{pp1} plots the value $\ln|\mathbb{E}\phi ({X^N}(T))-\mathbb{E}\phi(X^N_K)|,$ against $\ln(\delta t)$ for five different step sizes $\delta t=[2^{-5},2^{-6},2^{-7},2^{-8},2^{-9}]$ at $T=1$, where ${X^N}(T)$ and $X^N_K$ represent the exact and numerical solutions at the terminal time $T$, respectively. 
Here, the expectation $\mathbb{E}$ is approximated by taking average over 2000 realizations.
It can be seen that \eqref{full-spde} is of weak order 0.5 for cylindrical Wiener process, i.e. $\kappa=0$,
and of  weak order 1 for Q-Wiener process with $\kappa=0.5,2,$ which are indicated by the reference lines.
These coincide with the theoretical analysis.

\begin{figure}
\centering
\subfigure[$\phi(u)=\sin(\|u\|^2)$]{
\begin{minipage}{0.31\linewidth}
\includegraphics[width=4cm,height=4cm]{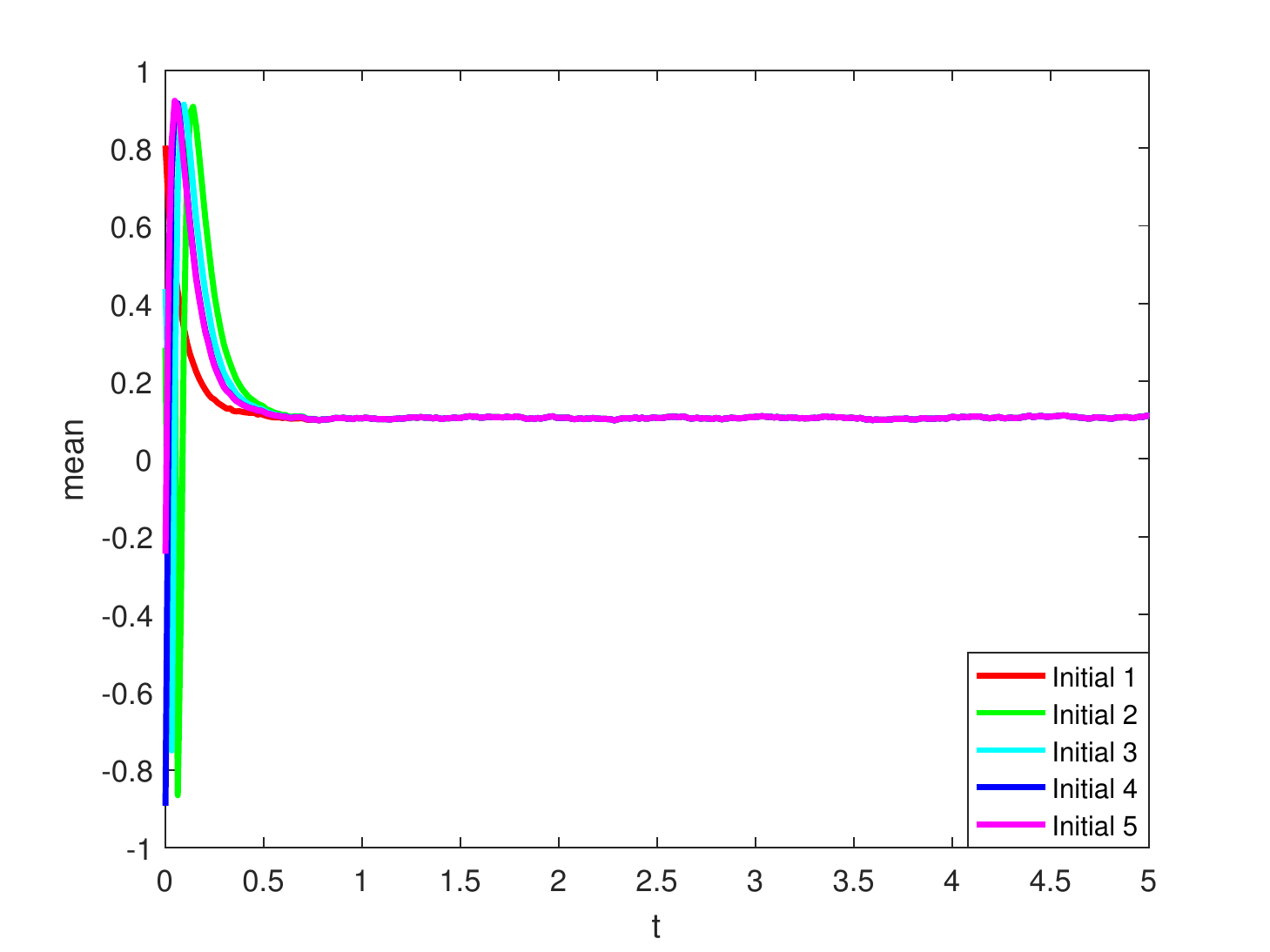}
\end{minipage}}
\subfigure[$\phi(u)=\sqrt{2}\cos(\|u\|^2-\frac \pi 4)$]{
\begin{minipage}{0.31\linewidth}
\includegraphics[width=4cm,height=4cm]{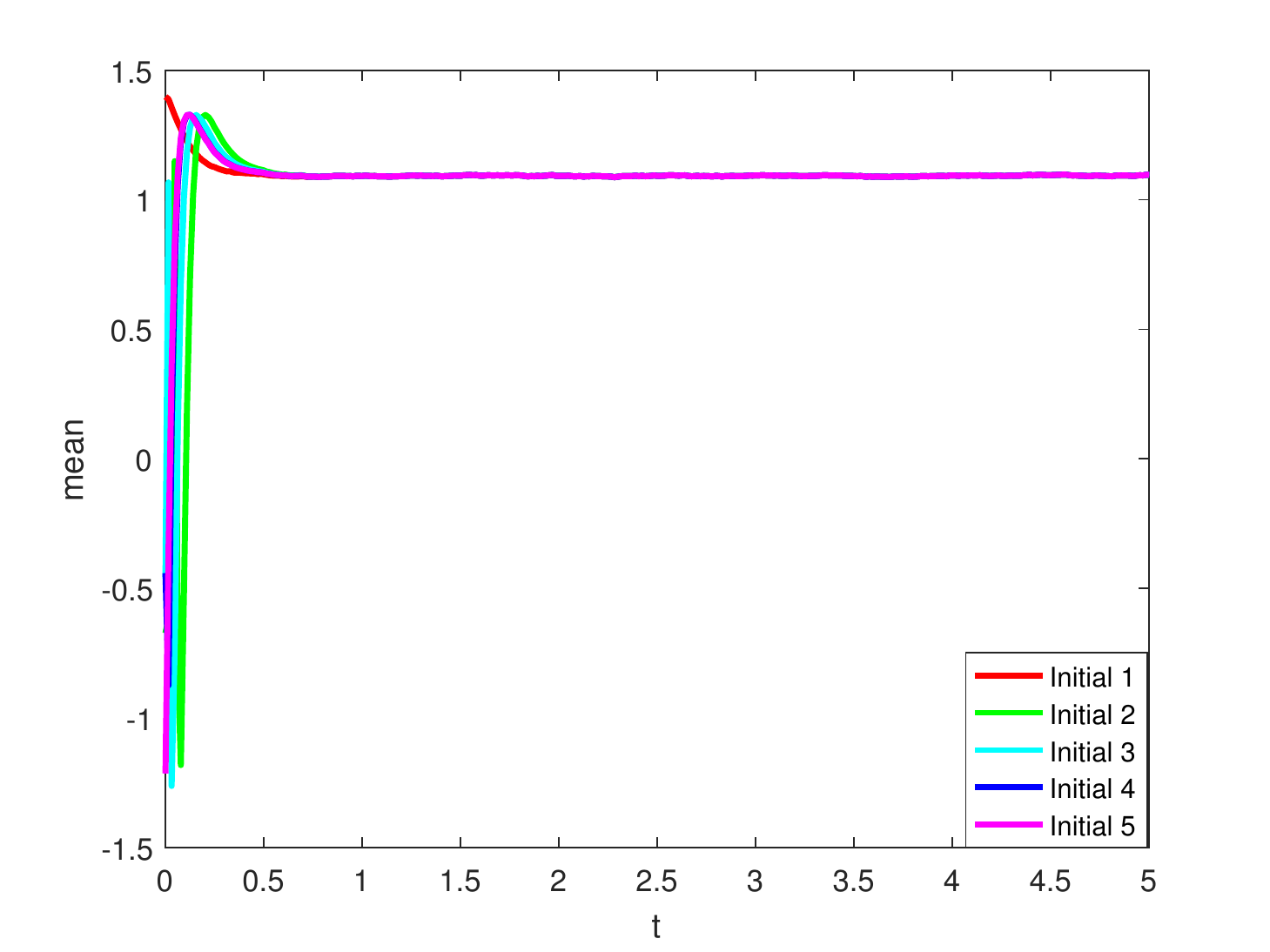}
\end{minipage}}
\subfigure[$\phi(u)=\exp(-\|u\|^2)$]{
\begin{minipage}{0.31\linewidth}
\includegraphics[width=4cm,height=4cm]{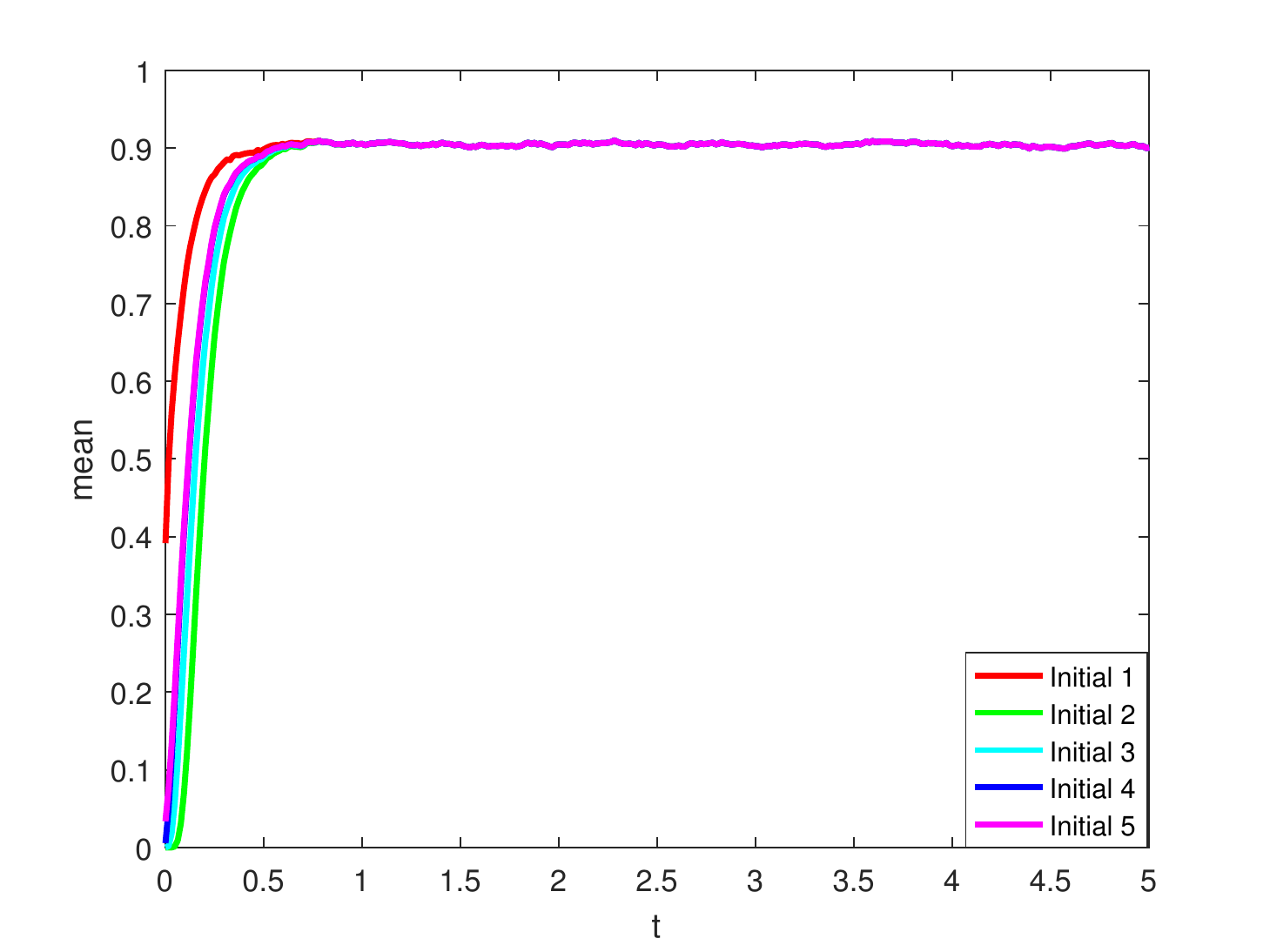}
\end{minipage}}
\caption{The averages $\E[\phi(X_k^N)]$ started from different initial values $(\lambda_F=5,$ $\delta t=2^{-6},$ $T=5)$}\label{pp3}
\end{figure}

\begin{figure}
\centering
\subfigure[$\phi(u)=\sin(\|u\|^2)$]{
\begin{minipage}{0.31\linewidth}
\includegraphics[width=4cm,height=4cm]{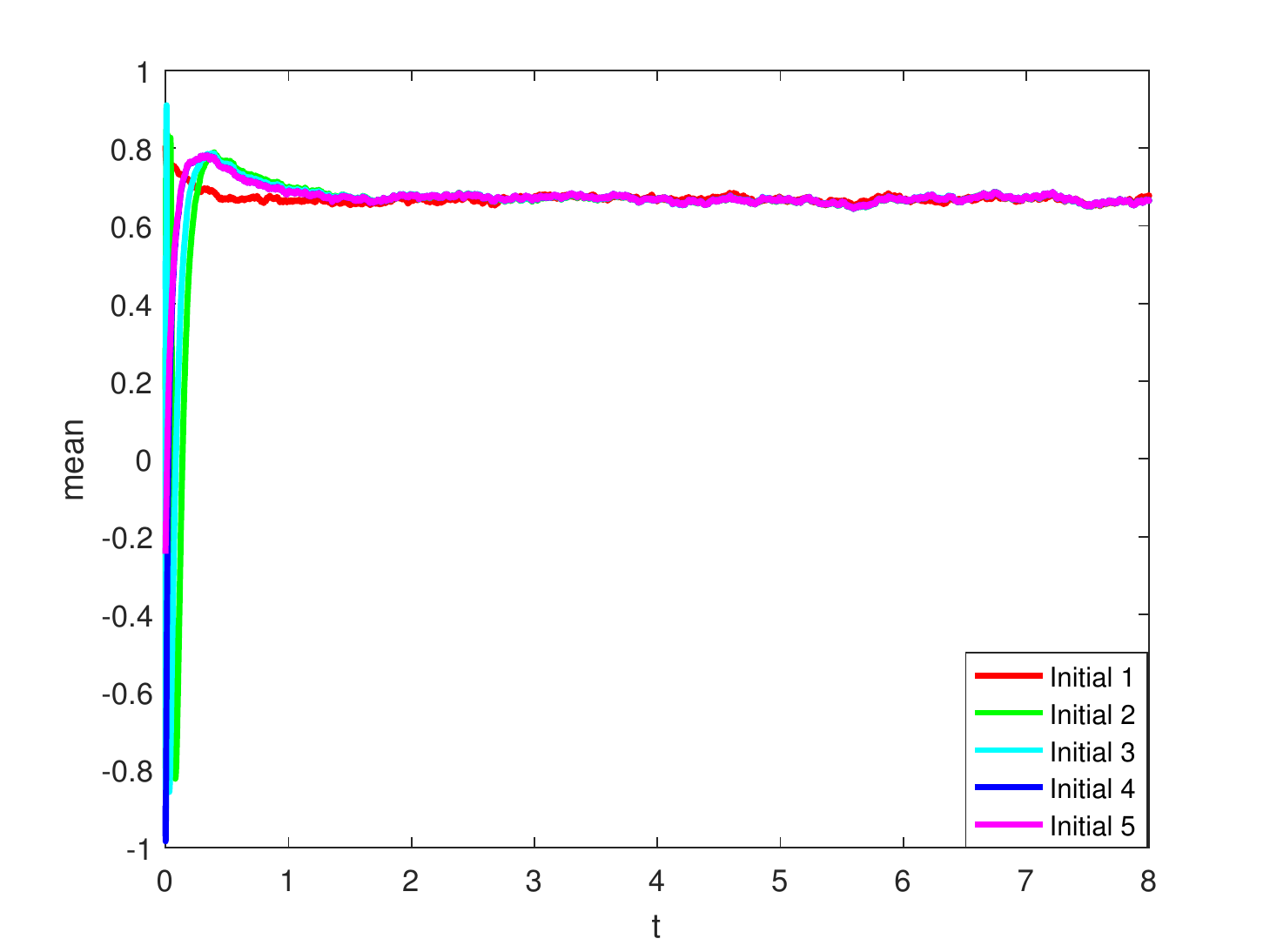}
\end{minipage}}
\subfigure[$\phi(u)=\sqrt{2}\cos(\|u\|^2-\frac \pi 4)$]{
\begin{minipage}{0.31\linewidth}
\includegraphics[width=4cm,height=4cm]{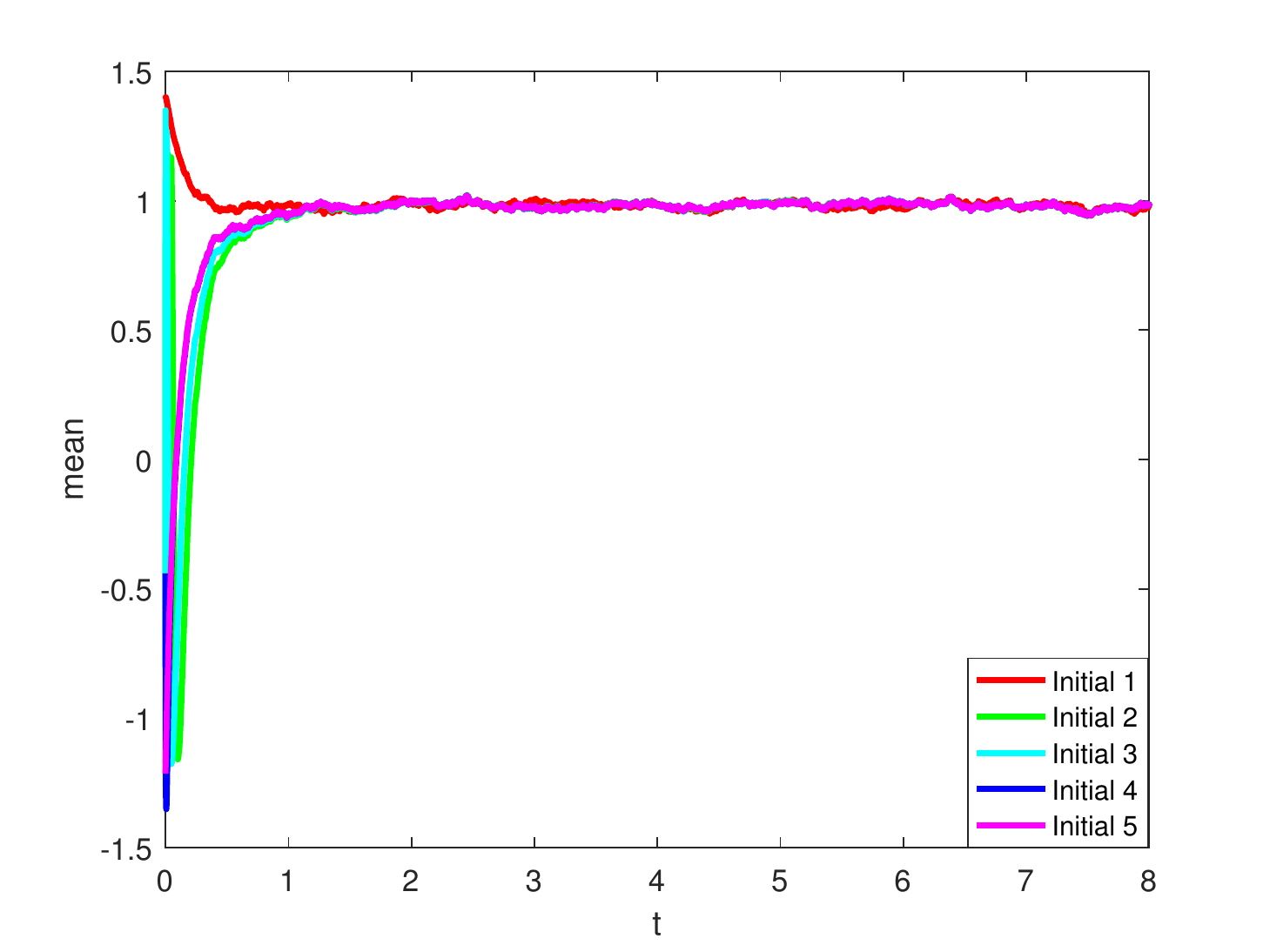}
\end{minipage}}
\subfigure[$\phi(u)=\exp(-\|u\|^2)$]{
\begin{minipage}{0.31\linewidth}
\includegraphics[width=4cm,height=4cm]{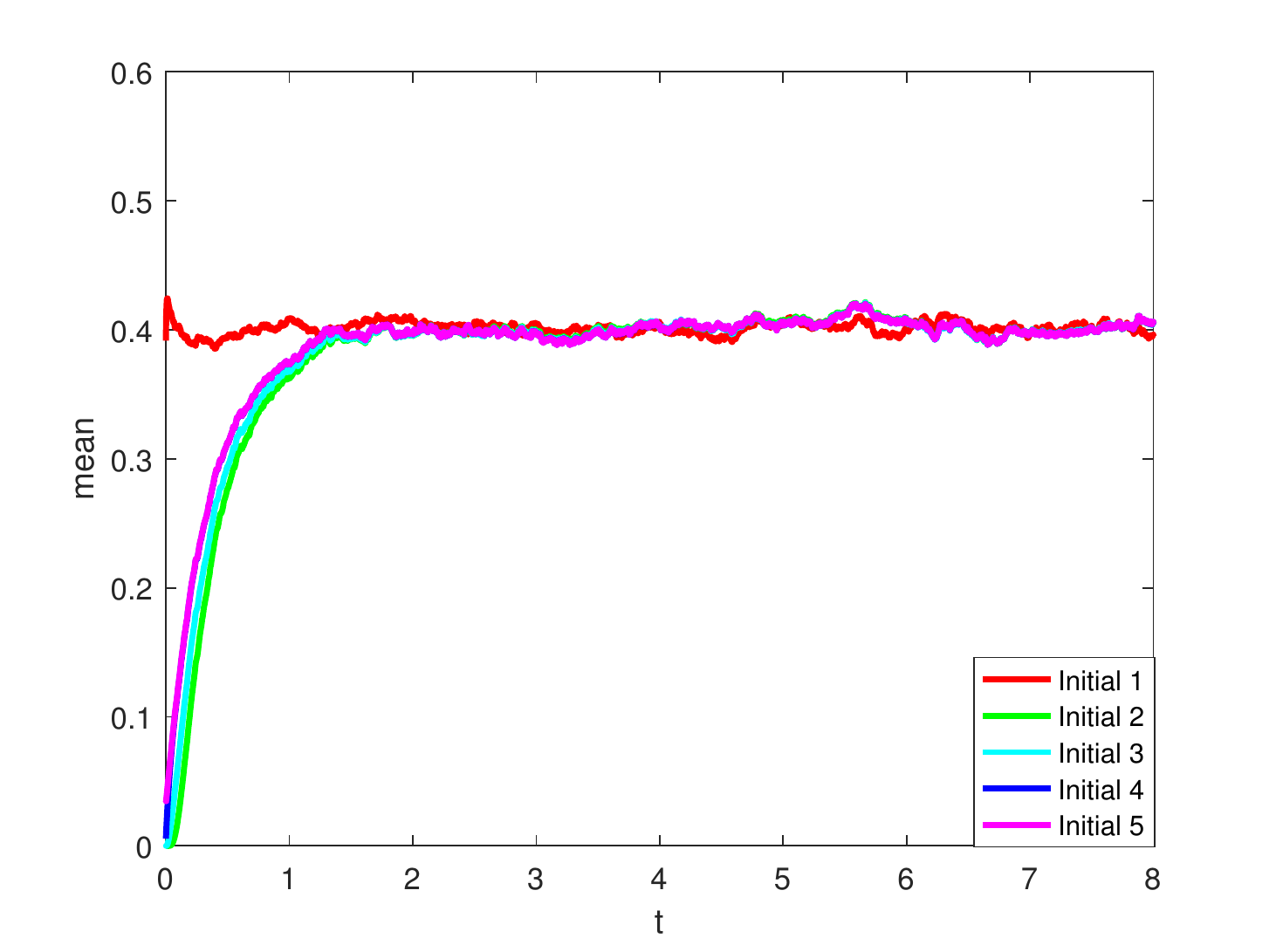}
\end{minipage}}
\caption{The averages $\E[\phi(X_k^N)]$ started from different initial values $(\lambda_F=12,$ $\delta t=2^{-10},$ $T=8)$}\label{pp4}
\end{figure}

Then we consider the longtime behaviors of \eqref{full-spde}.
Based on the definition of ergodicity, if numerical solution \eqref{full-spde} is strongly mixing, the average 
$\E[\phi(X_k^N)]$, $k>0$, started from different initial values will converge to the spatial average for almost every path. 
To verify this property and to make clear how the average value changes when time $t$ goes, Fig. \ref{pp3} shows the average of the fully discrete method started from five different initial values with the terminal time $T$ being 5 and $\kappa=0.$
It indicates that $\E[\phi(X_k^N(X_0))]$ started from different initial values converge to the same value in a short time for three  different kinds of continuous and bounded functions $\phi.$ 
Due to the exponential convergence to equilibrium,
the terminal time chosen here is not very large.
Moreover, aiming at
verifying that the mixed ergodicity does not need the condition $\lambda_F < \lambda_1,$ we also show the case $\lambda_F=12$ which implies 
$\lambda_F>\lambda_1$ in Fig. \ref{pp4}.
It can be seen that for different test functions, the averages will converge to the same value.
Numerical tests confirm theoretical findings.
Besides, the averages started from different initial values will also converge for Eq. \ref{spde} driven by other Q-Wiener processes.
For simplicity, we do not show those figures here.\\

\section{References}

\bibliographystyle{plain} 
\bibliography{bib}
\end{document}